\documentclass[leqno, 12pt]{amsart}
\usepackage{amsmath}
\usepackage{amsxtra}
\usepackage{amscd}
\usepackage{amsfonts}
\usepackage{eucal}
\usepackage{scalerel}
\usepackage{verbatim}
\usepackage{hyperref}
\usepackage{amssymb}
\usepackage{tikz-cd}
\usepackage{relsize}
\usepackage[matrix,arrow,curve]{xy}
\usepackage{stmaryrd}

\newcommand{\nc}{\newcommand}
\newcommand{\rc}{\renewcommand}

\nc{\ra}{{      \rightarrow     }}
\nc{\laa}{{     \leftarrow      }}      
\nc{\lra}{{\longrightarrow}}
\nc{\lr}{{\leftrightarrow}}             
\nc{\lrs}{{\rightleftarrows}}           

\nc{\imp}{{\Rightarrow}}                
\nc{\eq}{{\Leftrightarrow}}             

\nc{\inj}{{\protect  \hookrightarrow }}              
\nc{\injj}{{\protect \hookleftarrow  }}              

\nc{\sur}{{     \twoheadrightarrow      }}      
\nc{\surr}{{    \twoheadleftarrow       }}      

\nc{\mm}{{\mapsto}}                             

\nc{\va}{{\uparrow}}                            

\nc{\ul}{\underline}
\nc{\sub}{\subseteq}

\nc{\se}{       \section                }
\nc{\sus}{      \subsection             }
\nc{\sss}{      \subsubsection          }

\nc{\Lemm}{     \subsection{Lemma}              }
\nc{\lemm}{     \subsubsection{Lemma}           }
\nc{\slemm}{    \subsubsection*{Lemma}          }
\nc{\sublemm}{  \subsubsection{\bf Sublemma}            }
\nc{\ssublemm}{         \subsubsection*{\bf Sublemma}           }

\nc{\Pro}{      \subsection{Proposition}        }
\nc{\pro}{      \subsubsection{Proposition}     }
\nc{\spro}{     \subsubsection*{Proposition}    }

\nc{\Corr}{     \subsection{Corollary}          }
\nc{\corr}{     \subsubsection{\bf Corollary}       }
\nc{\scorr}{    \subsubsection*{\bf Corollary}      }

\nc{\Theo}{     \subsection{Theorem}            }
\nc{\theo}{     \subsubsection{\rm{\bf Theorem}} }
\nc{\stheo}{    \subsubsection*{\rm{\bf Theorem}}        }

\nc{\rem}{      \subsubsection{Remark}          }
\nc{\srem}{     \subsubsection*{Remark} }

\nc{\rems}{     \subsubsection{Remarks}         }
\nc{\srems}{    \subsubsection*{Remarks}        }

\nc{\conj}{     \subsubsection{\bf Conjecture}      }
\nc{\sconj}{    \subsubsection*{\bf Conjecture}     }

\nc{\ex}{       \subsubsection{Example}         }
\nc{\sex}{      \subsubsection*{Example}        }
\nc{\exs}{      \subsubsection{Examples}        }
\nc{\sexs}{     \subsubsection*{Examples}       }

\nc{\que}{      \subsubsection{Question}        }
\nc{\ques}{     \subsubsection{Questions}       }
\nc{\sque}{     \subsubsection*{Question}       }
\nc{\sques}{    \subsubsection*{Questions}      }

\nc{\pl}{{\oplus}}                              
\nc{\tim}{{\times}}
\nc{\btim}{{\boxtimes}}
\nc{\ltim}{\ltimes}                     %
\nc{\rtim}{\rtimes}                     %
\nc{\ltr}{\triangleleft}        %
\nc{\rtr}{\triangleright}       %

\nc{\ten}{{     \otimes         }}
\nc{\Lten}{{    \aa{L}\otimes   }}            
\nc{\Ltim}{{    \aa{L}\times    }}            
\nc{\bten}{{\boxtimes}}                         

\nc{\conl}{{    \begin{CD} @>{\cong}>>\end{CD}     }}
\nc{\cD} [2]{{   \begin{CD} @>{#1}>{#2}>\end{CD}     }}
\nc{\cDD}[2]{{   \begin{CD} @<{#1}<{#2}<\end{CD}     }}
\nc{\conn}{{     \begin{CD} @<{\cong}<<  \end{CD}   }}
\nc{\Con}{{     \equiv          }}      
\nc{\appr}{{    \sim            }}      
\nc{\eqr}{{     \sim            }}      

\nc{\ha}{{ \frac{1}{2} }}               
        \nc{\half}{{ \frac{1}{2} }}

\nc{\ci}{{\circ}}               
\nc{\cd }{{\cdot}}              
\nc{\cddd}{{\cdots}}

\nc{\cupp}{\bigcup}             
\nc{\capp}{\bigcap}
\nc{\pll}{\bigoplus}

\nc{\pii}{\prod}                
\nc{\ppii}{\bigprod}            
\nc{\cci}{\sqcup}              
\nc{\ccii}{\bigsqcup}

\nc{\wwe}{\bigwedge}            
\nc{\cce}{\bigcoprod}           

\nc{\aaa}{      \stackerel      }       

\nc{\Ker}{{    \operatorname{Ker}      }}
\rc{\Im}{{      \operatorname{Im}       }}
\nc{\rank}{{    \ \operatorname{rank}\  }}
\nc{\Res}{{     \  \operatorname{Res}   }}
\nc{\Hom}{{    \operatorname{Hom}      }}
\nc{\End}{{     \operatorname{End}      }}
\nc{\RHom}{{    \operatorname{RHom}     }}
\nc{\HHom}{{    \operatorname{$\HH$om}  }}
\nc{\EEnd}{{    \operatorname{$\EE nd$} }}
\nc{\AAut}{{    \operatorname{$\AA ut$} }}
\nc{\RHHom}{{   \text{R}\HH om	    }}
\nc{\Ext}{\operatorname{Ext}}
\nc{\Der}{{     \operatorname{Der}      }}

\nc{\RGa}{{	\text{R$\Ga$}	}}

\nc{\ord        }{{ \operatorname{ord} }}                       
\nc{\divv       }{{ \operatorname{div} }}                       
\nc{\Lie        }{{ \operatorname{Lie} }}

\nc{\bb}{	\pr\underset 	}           
\rc{\aa}{ 	\pr\overset 	}            

\rc{\AA}{{\mathcal A}}
\nc{\BB}{{\mathcal B}}
\nc{\CC}{{\mathcal C}}
\nc{\DD}{{\mathcal D}}
\nc{\EE}{{\mathcal E}}
\nc{\FF}{{\mathcal F}}
\nc{\GG}{{\mathcal G}}
\nc{\HH}{{\mathcal H}}
\nc{\II}{{\mathcal I}}
\nc{\JJ}{{\mathcal J}}
\nc{\KK}{{\mathcal K}}
\nc{\LL}{{\mathcal L}}
\nc{\MM}{{\mathcal M}}
\nc{\NN}{{\mathcal N}}
\nc{\OO}{{\mathcal O}}
\nc{\PP}{{\mathcal P}}
\nc{\QQ}{{\mathcal Q}}
\nc{\RR}{{\mathcal R}}
\rc{\SS}{{\mathcal S}}
\nc{\TT}{{\mathcal T}}
\nc{\UU}{{\mathcal U}}
\nc{\VV}{{\mathcal V}}
\nc{\WW}{{\mathcal W}}
\nc{\ZZ}{{\mathcal Z}}
\nc{\XX}{{\mathcal X}}
\nc{\YY}{{\mathcal Y}}

\nc{\A}{{\mathbb A }}
\nc{\B}{{\mathbb B}}
\nc{\C}{{\mathbb C}}
\nc{\D}{{\mathbb D}}
\nc{\E}{{\mathbb E}}
\nc{\F}{{\mathbb F}}
\nc{\G}{{\mathbb G}}
\nc{\hH}{{\mathbb H}}
\nc{\I}{{\mathbb I}}
\nc{\J}{{\mathbb J}}
\nc{\K}{{\mathbb K}}
\nc{\lL}{{\mathbb L}}
\rc{\L}{{\mathbb L}}
\nc{\M}{{\mathbb M}}
\nc{\N}{{\mathbb N}}
\nc{\oO}{{\mathbb O}}

\rc{\P}{{\mathbb P}}
\nc{\pP}{{\mathbb P}}
\nc{\Q}{{\mathbb Q}}
\nc{\R}{{\mathbb R}}
\nc{\sS}{{\mathbb S}}
\nc{\T}{{\mathbb T}}
\nc{\U}{{\mathbb U}}
\nc{\V}{{\mathbb V}}
\nc{\W}{{\mathbb W}}
\nc{\Z}{{\mathbb Z}}
\nc{\X}{{\mathbb X}}
\nc{\Y}{{\mathbb Y}}

\rc{\k}{\Bbbk}

\nc{\fA}{{\mathfrak A}}
\nc{\fB}{{\mathfrak B}}
\nc{\fC}{{\mathfrak C}}
\nc{\fD}{{\mathfrak D}}
\nc{\fE}{{\mathfrak E}}
\nc{\fF}{{\mathfrak F}}
\nc{\fG}{{\mathfrak G}}
\nc{\fH}{{\mathfrak H}}
\nc{\fI}{{\mathfrak I}}
\nc{\fJ}{{\mathfrak J}}
\nc{\fK}{{\mathfrak K}}
\nc{\fL}{{\mathfrak L}}
\nc{\fM}{{\mathfrak M}}
\nc{\fN}{{\mathfrak N}}
\nc{\fO}{{\mathfrak O}}
\nc{\fP}{{\mathfrak P}}
\nc{\fQ}{{\mathfrak Q}}
\nc{\fR}{{\mathfrak R}}
\nc{\fS}{{\mathfrak S}}
\nc{\fT}{{\mathfrak T}}
\nc{\fU}{{\mathfrak U}}
\nc{\fV}{{\mathfrak V}}
\nc{\fW}{{\mathfrak W}}
\nc{\fZ}{{\mathfrak Z}}
\nc{\fX}{{\mathfrak X}}
\nc{\fY}{{\mathfrak Y}}
\nc{\fa}{{\mathfrak a}}
\nc{\fb}{{\mathfrak b}}
\nc{\fc}{{\mathfrak c}}
\nc{\fd}{{\mathfrak d}}
\nc{\fe}{{\mathfrak e}}
\nc{\ff}{{\mathfrak f}}
\nc{\fg}{{\mathfrak g}}
\nc{\fh}{{\mathfrak h}}
\nc{\fj}{{\mathfrak j}}
\nc{\fk}{{\mathfrak k}}
\nc{\fl}{{\mathfrak{l}}}
\nc{\fm}{{\mathfrak m}}
\nc{\fn}{{\mathfrak n}}
\nc{\fo}{{\mathfrak o}}
\nc{\fp}{{\mathfrak p}}
\nc{\fq}{{\mathfrak q}}
\nc{\fr}{{\mathfrak r}}
\nc{\fs}{{\mathfrak s}}
\nc{\ft}{{\mathfrak t}}
\nc{\fu}{{\mathfrak u}}
\nc{\fv}{{\mathfrak v}}
\nc{\fw}{{\mathfrak w}}
\nc{\fz}{{\mathfrak z}}
\nc{\fx}{{\mathfrak x}}
\nc{\fy}{{\mathfrak y}}

\nc{\al}{{\alpha }}
\nc{\be}{{\beta }}
\nc{\ga}{{\gamma }}
\nc{\de}{{\delta }}
\nc{\del}{{\partial }}
\nc{\ep}{{\varepsilon }}
\nc{\vap}{{\epsilon }}

\nc{\ze}{{\zeta }}
\nc{\et}{{\eta }}
\rc{\th}{{\theta }}
\nc{\vth}{{\vartheta }}

\nc{\io}{{\iota }}
\nc{\ka}{{\kappa }}
\nc{\la}{{\lambda }}
\nc{\vrho}{{\varrho}}
\nc{\si}{{\sigma }}
\nc{\ups}{{\upsilon }}
\nc{\vphi}{{\varphi }}
\nc{\om}{{\omega }}

\nc{\Ga}{{\Gamma }}
\nc{\De}{{\Delta }}
\nc{\nab}{{\nabla}}
\nc{\Th}{{\Theta }}
\nc{\La}{{\Lambda }}
\nc{\Si}{{\Sigma }}
\nc{\Ups}{{\Upsilon }}
\nc{\Om}{{\Omega }}

\nc{\toc}{\tableofcontents}
\nc{\addl}{     \addcontentsline{toc}{subsection}       }

\nc{\zero}{{\text{\bf{0}}}}
\nc{\one}{{\text{\bf{1}}}}

\nc{\ie}{{,\ \     \text{i.e.,}\ \ 	}}
\nc{\iif}{{\ \     \text{if}\ \ 	}}
\nc{\aand}{{\ \ \  \text{and}\ \ \ 	}}
\nc{\hence}{{\ \ \ \text{hence}\ \ \ 	}}
\nc{\while}{{\ \ \ \text{while}\ \ \ 	}}
\nc{\with}{{\ \ \  \text{with}\ \ \ 	}}
\nc{\oor}{{\ \     \text{or}\ \ 	}}
\nc{\foor}{{\ \     \text{for}\ \ 	}}
\nc{\suchthat}{{\ \     \text{such that}\ \ 	}}

\nc{\subb}{{	\supseteq	}}         
\nc{\nsub}{{	\nsubseteq	}}         
\nc{\nsubb}{{	\nsupseteq	}}         %

\nc{\nin}{{	\notin	}}	

\let\b\big

\nc{\df}{{ \protect\overset{ \text{def}}= 	}}		
\nc{\dff}{{ \ \df\				}}		
\nc{\inv}{{ {}^{-1}      }}			

\nc{\tx}{{	\tii x				     }}

\nc{\vd}{{	(v,d)				     }}
\nc{\Mvd}{	M(v,d)		}
\rc{\Mc}{	M^c		}
\nc{\fMc}{	\fM^c		}
\nc{\fMz}{	\fM^0		}
\nc{\Mcvd}{	M^c(v,d)	}
\nc{\Mzvd}{	M^0(v,d)	}
\nc{\fMvd}{	\fM(v,d)	}
\nc{\fMvdz}{	\fM_\zero(v,d)	}
\nc{\Lacvd}{{	\La^c(v,d)	}}

\nc{\sA}{{\mathsf A}}
\nc{\sB}{{\mathsf B}}
\nc{\sC}{{\mathsf C}}
\nc{\sD}{{\mathsf D}}
\nc{\sE}{{\mathsf E}}
\nc{\sF}{{\mathsf F}}
\nc{\sG}{{\mathsf G}}
\nc{\sH}{{\mathsf H}}
\nc{\sI}{{\mathsf I}}
\nc{\sJ}{{\mathsf J}}
\nc{\sK}{{\mathsf K}}
\nc{\sL}{{\mathsf L}}
\nc{\sM}{{\mathsf M}}
\nc{\sN}{{\mathsf N}}
\nc{\sO}{{\mathsf O}}
\nc{\sP}{{\mathsf P}}
\nc{\sQ}{{\mathsf Q}}
\nc{\sR}{{\mathsf R}}
\rc{\sS}{{\mathsf S}}
\nc{\sT}{{\mathsf T}}
\nc{\sU}{{\mathsf U}}
\nc{\sV}{{\mathsf V}}
\nc{\sW}{{\mathsf W}}
\nc{\sX}{{\mathsf X}}
\nc{\sY}{{\mathsf Y}}
\nc{\sZ}{{\mathsf R}}
\nc{\sa}{{\mathsf a}}
\rc{\sb}{{\mathsf b}}
\rc{\sc}{{\mathsf c}}
\nc{\sd}{{\mathsf d}}

\nc{\es}{{\mathsf e}}
\nc{\sg}{{\mathsf g}}
\nc{\sh}{{\mathsf h}}
\nc{\sj}{{\mathsf j}}
\nc{\sk}{{\mathsf k}}
\nc{\sn}{{\mathsf n}}
\nc{\so}{{\mathsf o}}
\nc{\sr}{{\mathsf r}}
\nc{\su}{{\mathsf u}}
\nc{\sv}{{\mathsf v}}
\nc{\sw}{{\mathsf w}}
\nc{\sx}{{\mathsf x}}
\nc{\sy}{{\mathsf y}}
\nc{\sz}{{\mathsf z}}

\nc{\sq}{{\mathsf q}}

\nc{\pf}{ \begin{proof}}
\nc{\epf}{ \end{proof}}

\nc{\bbb}{ 	\boldsymbol 	}

\nc{\ftt}[1]{{\footnote{#1}}}
\nc{\fttt}[1]{{$^($\footnote{#1}$^)$}}
\nc{\bftt}[1]{\footnote{#1}}

\nc{\Vn}{	V^{\gl_n}	}
\nc{\Vm}{	V^{\gl_m}	}

\nc{\syp}[1]{	^{ (#1) }		} 	
\nc{\up}[1]{	^{ (#1) }		} 	
\nc{\lp}[1]{	_{ (#1) }		}	
\nc{\hp}[1]{	^{ [#1] }		}	

\nc{\bOO}{{	 \sO }}
\nc{\bKK}{{	 \sK }}
\nc{\bOOm}{{	 \sO^- }}
\nc{\OOm}{{	 \sO^- }}

\nc{\Ao}{{	\A^1	}}
\nc{\Po}{{	\P^1	}}
\nc{\So}{{	S^1	}}

\nc{\h}{{	\hslash	}}	

\nc{\An}{{	\A\syp{n}	}}

\nc{\ub}{{	\underline{b}		}}

\nc{\cG}{{	\ch G		}}
\nc{\cB}{{	\ch B		}}
\nc{\cN}{{	\ch N		}}
\nc{\cT}{{	\ch T		}}
\nc{\cH}{{	\ch H		}}

\nc{\bss}{{\backslash}}           		
\nc{\barr}{ 	\overline 	}      		
\nc{\ud}{	\underline	}		

\nc{\ti}{\tilde}              
\nc{\tii}{\widetilde}         

\nc{\hatt}{\widehat}				
\nc{\hata}{{	\bbb{ \hat{} }		}}	

\nc{\ch}{\check}              			

\nc{\tnu}{{	\underline{\nu}			}}
\nc{\un}{{	\underline{\nu}		}}
\nc{\um}{{	\underline{\mu}		}}
\nc{\tmu}{{	\tii\mu			}}

\nc{\glm}{{		  \gl_m	}}
\nc{\gln}{{		  \gl_n	}}

\nc{\bi}{	\begin{itemize}\item		}
\rc{\i}{	\item			}
\nc{\ei}{ \end{itemize}	}
\nc{\ben}{	\begin{enumerate}\item		}
\nc{\een}{	\end{enumerate}			}

\nc{\tgab}{	    \tfg^{a,\sb}			}
\nc{\tgas}{	    \tfg^{a,\bu}			}

\nc{\nGG}{{	 \syp{n}\GG	  }}

\nc{\slem}{ 	\subsubsection*{Lemma}		}
\nc{\scor}{ 	\subsubsection*{Corollary}	}

\nc{\lap}{{	\si				}}

\nc{\yy}{\infty}
\nc{\ys}{{  \frac{\infty}{2}  }}

\nc{\bu}{ \bullet         }  			
\nc{\Gm}{{		  G_m			}}

\nc{\lb}{\langle}             				
\nc{\rb}{\rangle}

\nc{\pmo}{{ 	\pm 1		}}
\nc{\mpo}{{ 	\mp 1		}}

\nc{\Xn}{{	X\syp{n}	}}

\nc{\Coh}{{	\CC oh		}}               %

\nc{\lab}{\label}

\nc{\bububu}{{ \bb{\bu} { \aa{\bu}\bu }                }}

\nc{\ftx}[1]{{${\bububu}{}^($\footnote{#1}$^)\bububu $}}

\newcounter{suggestions}
\nc{\adc}[1]{{\addtocounter{#1}{1}}}
\nc{\fff}[1]{\ftx{ \fbox{S$\thesuggestions.$}\adc{suggestions} {#1} }}

\nc{\sugg}{\bf SUGGESTION}
\nc{\SU}[1]{\fff{\sugg$@>>>$}}
	\nc{\dv}{\tx{\bf A dying version}}
\nc{\DV}[1]{\fff{{\dv$@>>>$}}}
	\nc{\cracra}{\tx{\bf{END}}}
\nc{\ENDD}{\begin{flushright}{\fff{$@<<<$\cracra}}\end{flushright}}
\nc{\END}{{\fff{$@<<<$\cracra}}}
\nc{\ENDDV}{\begin{flushright}{\fff{$@<<<$\tx{\bf END of \dv}}}\end{flushright}}

\nc{\IN}[1]{	    \fbox{IN?:}$\leftarrow$   \fbox{#1}   }
\nc{\OUT}[1]{	    \fbox{OUT?:}$\rightarrow$ \fbox{#1}   }

\nc{\GK}{{      G(\sK)          }}
\nc{\GO}{{      G(\sO)          }}

\nc{\Gd}{{  {\check G}          }}

\nc{\Gh}{{  \hat{\GG}           }}
\nc{\GA}{{  G(\AA)              }}
\nc{\BK}{{  B(\mathcal K)           }}
\nc{\BO}{{  B(\mathcal O)           }}
\nc{\BKz}{{  B(\mathcal K)_0        }}
\nc{\NK}{{  N(\mathcal K)           }}
\nc{\NO}{{  N(\mathcal O)           }}
\nc{\TK}{{  B(\mathcal K)           }}
\nc{\TO}{{  T(\mathcal O)           }}
\nc{\LO}{{  L(\mathcal O)           }}

\nc{\gk}{{      \mathfrak g_\KK             }}
\nc{\bk}{{      \mathfrak b_\KK             }}
\nc{\tk}{{      \mathfrak b_\KK             }}
\nc{\nk}{{      \mathfrak n_\KK             }}
\nc{\go}{{      \mathfrak g_\OO             }}

\nc{\Ll}{{  L_\lambda   }}
\nc{\Lm}{{  L_\mu       }}
\rc{\Lm}{{  \sL^-       }}
\nc{\sLm}{ \sL^-       }

\nc{\Gl}{{      {\mathcal G_\lambda}        }}
\nc{\Glb}{{     \barr{\Gl}              }}
\nc{\Glm}{{     \GG_{\lambda+\mu}       }}
\nc{\Ge}{{      \GG_\eta                }}
\nc{\Geb}{{     \barr{\Ge}              }}
\nc{\Gwl}{{     \GG_{W \lambda}         }}
\nc{\Gn}{{      \GG_\nu                 }}
\nc{\Gnb}{      \barr{\GG_\nu           }}
\nc{\Gum}{{     \GG^\mu         }}
\nc{\Gul}{{     \GG^\lambda     }}
\nc{\Gun}{{     \GG^\nu         }}
\nc{\Pl}{{  \PP_\lambda         }}
\nc{\Pwl}{{  {\PP}_{W \lambda}  }}
\nc{\FFl}{{  {\FF}_\lambda      }}

\nc{\Cn}{{  C\syp{n}       }}
\nc{\Ce}{{  C_\eta      }}

\nc{\Sn }{{     S_\nu           }}
\nc{\Snb}{{  \barr{S_\nu}       }}
\nc{\Sl }{{  S_\lambda          }}
\nc{\Sm }{{  S_\mu              }}

\nc{\GGG}{{  G(\KK)\underset{G(\OO)}\times      \GG                     }}
\nc{\GGlm}{{  (G(\KK)\underset{G(\OO)}\times    \GG)_{(\lambda,\mu)}    }}

\nc{\dgg}{{     \ddot\GG                }}
\nc{\dlm}{{     \ddot\GG_{\la,\mu}      }}

\nc{\db}{{      \hat \Delta     }}
\nc{\dr}{{      \Delta_\R       }}
\nc{\Dp}{{      \Delta^+        }}
\nc{\Dm}{{      \Delta^-        }}
\nc{\da}{{      \Delta_a        }}
\nc{\dap}{{     \Delta_a^+      }}
\nc{\di}{{      \Delta_I        }}
\nc{\cro}{{     \check \rho     }}      
\nc{\dc}{{      \check \delta   }}
\nc{\rac}{{     \check \rho_a   }}

\nc{\Irr}{\operatorname {Irr}}
\nc{\td}{\widetilde d}
\nc{\tg}{\widetilde g}
\nc{\tv}{\widetilde v}
\nc{\tp}{\tilde p}
\nc{\tmm}{\widetilde{\bf m}}

\nc{\TD}{\widetilde D}
\nc{\TV}{\widetilde V}
\nc{\TU}{\widetilde U}
\nc{\tT}{\widetilde T}
\nc{\TA}{\widetilde A}
\nc{\TB}{\widetilde B}
\nc{\TN}{\widetilde{\mathcal N}}
\nc{\TF}{\widetilde{\mathcal F}}
\nc{\TGG}{\widetilde{\mathcal G}}

\nc{\tga}{\widetilde\ga}
\nc{\tde}{\widetilde\de}
\nc{\tphi}{\widetilde\phi}
\nc{\tpsi}{\widetilde\psi}

\nc{\cla}{{\check\lambda}}
\nc{\cmu}{{\check\mu}}
\nc{\cnu}{{\check\nu}}
\nc{\cth}{{\check\theta}}

\nc{\mmm}{{\bf m}}

\nc{\bx}{\barr x}

\nc{\Sym}{\operatorname {Sym}}
\nc{\Spec}{\operatorname {Spec}}
\nc{\Id}{\operatorname {Id}}

\nc{\gl}{\mathfrak{gl}}
\nc{\tfg}{{\tii\fg}}

\nc{\BDGG}{{\bf\mathfrak G}}

\nc{\tGG}{{\tii\GG}}
\nc{\tGGr}{{\tii\GG^{rat}}}
\nc{\tGGra}{{\tii\GG^{rat,a}}}
\nc{\tGGrab}{{\tii\GG^{rat,a,\sb}}}

\nc{\TBDGG}{\widetilde{\GG}}

\nc{\BDG}{{		 \GG			}}
\nc{\TBDG}{{			    \widetilde\GG	}}

\nc{\Perv}{\operatorname {Perv}}
\nc{\Rep}{\operatorname {Rep}}
\nc{\Mlt}{\operatorname {Mlt}}
\nc{\spr}{\operatorname {Spr}}

\nc{\IC}{\operatorname {IC}}
\nc{\Tr}{\operatorname {Tr}}
\nc{\sign}{\operatorname{sgn}}

\nc{\ad}{\operatorname{ad}}
\nc{\proj}{\operatorname{pr}}
\nc{\pr}{\protect}

\nc{\bfb}{{\bf b}}
\nc{\bfg}{{\bf g}}
\nc{\bfe}{{\bf e}}
\nc{\eb}{{\bf e}}

\nc{\mat} {		\left(		\matrix	}	
\nc{\emat}{		\endmatrix	\right)	}
\nc{\sm} {		\left(		\smallmatrix	}	
\nc{\esm}{		\endsmallmatrix	\right)	}
\nc{\smat} {		\left(		\smallmatrix	}	
\nc{\esmat}{		\endsmallmatrix	\right)	}

\nc{\matr} {		\left[		\matrix	}	
\nc{\ematr}{		\endmatrix	\right]	}
\nc{\smr} {		\left[		\smallmatrix	}	
\nc{\esmr}{		\endsmallmatrix	\right]	}
\nc{\smatr} {		\left[		\smallmatrix	}	
\nc{\esmatr}{		\endsmallmatrix	\right]	}

\nc{\imat} {		\left.		\matrix	}	
\nc{\eimat}{		\endmatrix	\right.	}
\nc{\ism} {		\left.		\smallmatrix	}	
\nc{\eism}{		\endsmallmatrix	\right.	}

\nc{\ca}{		\left\{		\smallmatrix	}	
\nc{\eca}{		\endsmallmatrix	\right\}	}
\nc{\Ca}{		\left\{		\matrix		}	
\nc{\Eca}{		\endmatrix	\right.		}	
\nc{\eCa}{		\endmatrix	\right\}	}	

\begin{document}

\title[]{
Comparison of quiver varieties,
loop Grassmannians and nilpotent cones in type A
}
\author{       Ivan Mirkovi\'c                 }
\address{Dept. of Mathematics and Statistics, University
of Massachusetts at Amherst, Amherst MA 01003-4515, USA}
\email{                mirkovic@math.umass.edu        }

\author{
Maxim Vybornov}

\author{ {\tiny with an appendix by} Vasily Krylov}
\address{Department of Mathematics
Massachusetts Institute of Technology
\newline
77 Massachusetts Avenue,
Cambridge, MA 02139,
USA;\newline
National Research University Higher School of Economics, Russian Federation\newline
Department of Mathematics, 6 Usacheva st., Moscow 119048;\newline
Skolkovo Institute of Science and Technology}
\email{krvas@mit.edu,
kr-vas57@yandex.ru}

\begin{abstract}

In type $A$ we find equivalences of geometries
arising
in three settings:\
 Nakajima's (``framed'') quiver varieties,
conjugacy classes of matrices
and loop Grassmannians.
These are all given by explicit formulas.
In particular,  we embedd the framed  quiver varieties
into Beilinson-Drinfeld Grassmannians. This
provides a compactification of Nakajima
varieties and a decomposition of  affine Grassmannians
into Nakajima varieties.
As an application we
provide a geometric version of  symmetric and skew $(GL(m), GL(n))$
dualities.
\end{abstract}



\maketitle

\nc{\con}{ \aa{\cong}\longrightarrow		   }
\nc{\emp}{ \emptyset		   }
\rc{\empty}{ \emptyset		   }

\nc{\edd}{ \end{document}	}
\nc{\f}[1]{ \fbox{$ $}\footnote{ \fbox{!}#1 }\fbox{$ $}		}

\rc{\sn}{N}
\nc{\Asn}{A\syp{n}}
\nc{\Csn}{C\syp{n}}


\nc{\bbp}{{ \bbb p	 }}
\nc{\limi}{{	\pr\underset {\rightarrow}\lim		}}      
\nc{\limp}{{	\pr\underset {\leftarrow}\lim		}}      

\nc{\Am}{{\A^-}}
\nc{\bp}{\begin{proof}}
\nc{\bep}{\begin{proof}}
\nc{\enp}{\end{proof}}
\nc{\fra}{\frac}

\nc{\sip}{^{\si+}}
\nc{\UsA}{\UU\sip_\An}

\nc{\ox}{\OO_\A}
\nc{\psiz}{\Psi_0}
\nc{\Psiz}{\Psi_0}


\nc{\thh}{	^{\text{th}}	}                     	
\nc{\hb}{{\hatt b}}
\nc{\hc}{{\hatt c}}
\nc{\hz}{{\hatt 0}}

\nc{\np}{\newpage}
\nc{\tvphi}{\tii\vphi}
\setcounter{tocdepth}{1}
\nc{\setse}{		\setcounter{section}	}

\nc{\hy}{{\hatt\yy}}
\nc{\ty}{{\tii\yy}}

\nc{\GGr}{\GG^{rat}}
\nc{\GGrat}{\GG^{rat}}
\nc{\UUratI}{ \UU^{rat,I} }

\nc{\tcUUratI}{ \tii\UU^{rat,I}_{\A,\Gm} }
\nc{\cUUratI}{ \UU^{rat,I}_{\A,\Gm} }

\nc{\cUU}{\UU_{\A,\Gm}}

\nc{\UUrs}{\UU^{rat,\si}_\A}

\nc{\bz}{{\barr z}}

\nc{\UtI}{U_\ty/I}

\nc{\uLL}{{\ud\LL}}
\nc{\uPhi}{{\ud\Phi}}

\nc{\vPhi}{\Phi}

\toc



\newcommand\iso{\,\vphantom{j^{X^2}}\smash{\overset{\sim}{\vphantom{\rule{0pt}{0.20em}}\smash{\longrightarrow}}}\,}

\se{Introduction}

In type A (only)
we relate three geometric objects
related to representations:\
Nakajima's quiver varieties, conjugacy classes of matrices
and Beilinson-Drinfeld Grassmannians.
The main result is an identification of certain normal slices to standard
stratifications in the three settings, in particular  the singularities are the
same in all three cases.
We also identify natural resolutions of slices.\fttt{
These observations clearly do  not literally extend
beyond type A. For instance, the closures of orbits in the loop Grassmannian
are normal and this is not true for
nilpotent orbits.
} An affine version of this paper  has been constructed in \cite{BF1}.

We  embed Nakajima's quiver varieties into the Beilinson-Drinfeld
Grassmannians and  this provides
a compactification of Nakajima's quiver varieties
and a decomposition of  loop Grassmannians
into a disjoint union of quiver varieties.
For  nilpotent orbits we construct new transverse slices
naturally related to loop Grassmannians.
As an application we construct a geometric version of both the skew and the symmetric
version of the $(GL_m, GL_n)$ duality.

\sss{
Loop Grassmannians and matrices
}
\lab{Loop Grassmannians and matrices}
To go from a loop Grassmannian $\GG$ of $G=GL_m$
to matrices in $\gl(N)$ we consider
a subvariety $\TT$ of $\GG$
such that all lattices $L\in\TT$, contain the standard lattice
$\sL_0\dff\C[[z]]^{\pl m}$.
Then a choice of
vector space trivializations  $\io_L:\C^N\con L/\sL_0$ for $L\in\TT$,
translates the operator $z$ on lattices  $L\in\TT$,
into a family of matrices
$x_L$ of size $N$.
It turns out that when $\TT$ is a standard normal slice
to a $G[[z]]$-orbit in $\GG$ then there is a natural trivialization
$\io$ and the corresponding family of
matrices is a (nonstandard) normal slice to a nilpotent orbit in $\gl(N)$.
So, one exchanges a setting with a fixed operators $z$ on $
\C((z))^{\pl m}/\sL_0$ and a variable invariant subspace $L/\sL_0$
to a setting of a fixed space $\C^N$ and a variable operator $x_L$.
The origin of our constructions is Lusztig's embedding of
the nilpotent cone for $GL_n$ into the loop Grassmannian of $GL_n$~\cite{L81}.

\sss{
Quiver varieties and nilpotent orbits
}
\lab{Quiver varieties and nilpotent orbits}
The relationship between these was conjectured by Nakajima \cite{N94} and
proved by Maffei \cite{M}. Our treatment of this relation
is close to
Maffei's work.
However, while he uses
Slodowy's normal slices to nilpotent orbits we use different
slices suggested by the relation to the loop Grassmannians, and
this makes the construction canonical and explicit while Maffei's approach is
based on an existence result. Explicit formulas for Maffei isomorphism were supplied by  \cite{IMW}
after the original version of this paper
appeared \cite{MVyb}.

\sss{
Isomorphism of quiver varieties and loop Grassmannians
}
These two objects are  in principle
related by combining
\ref{Loop Grassmannians and matrices}
and \ref{Quiver varieties and nilpotent orbits}.
Having an  explicit formula
is a deeper problem which is
resolved in the appendix~\ref{Appendix}
by Vasily Krylov
through a  completely new approach
using moduli of bundles over $\P^2$ (following ideas of Nakajima).

\sss{
}
Some of the results of this paper have appeared in an earlier version
\cite{MV}. The new results are the 
explicit isomorphism of quiver varieties and loop Grassmannians
(see ~\ref{Appendix}) and an isomorphism 
normal slices for conjugacy classes in matrices and for the Beilinson-Drinfeld deformation of loop
Grassmannians (previously we only had inclusions). Beyond this the paper has been
extensively rewritten.
 
\sus{
The Isomorphism Theorem
}
Here we formulate the main result more precisely.
We work over the field of complex numbers $\C$.

We consider
a quiver $(I,H)$ of type $A_{n-1}$.
Any
dimension vector
$d=(d_1,\dots ,d_{n-1})\in \N^{n-1}$
for this quiver
determines
integers $m\dff d_1+\dots +d_{n-1}$
and $N \dff \sum_{j = 1}^{n-1} jd_j$.
These give rise 
to the
loop Grassmannian $\GG$
for $GL_m$
and the nilpotent cone $\NN$ for
$\gl_N$.

In the quiver setting, using a choice of
$d,v$ in $\N^{n-1}$
and a central element $c=(c_1,\dots, c_{n-1})$
of the Lie algebra $\gl(v)\dff\ \pl_{i=1}^{n-1} \gl(v_i)$,
Nakajima \cite{N94, N98} constructs
a map of quiver varieties
$\fMc(v,d)\to \fMc_\zero(v,d)$.
These quiver data $(d,v)$ are also reformulated as a pair of
dominant coweights
$\la$ and $\mu$ for $G=GL_m$
which we can think of as a pair of
partitions $\la,\mu$  of $N$
(see \ref{qdatatogl} for the combinatorics of data).

For the loop Grassmannian setting
the loop group
$LG=G\b(\C((z)))=GL\b(m,\C((z))\b)$ contains
the ``disc group''
$L^{\geq 0}G=G(\C[[z]])$ and the negative congruence subgroup
$L^{<0}G$ (the kernel of the evaluation of $G(\C[[z\inv]]\to $ at $z\inv=0$).

Any  coweight $\eta$ of $GL_m$
defines a point
in the loop Grassmannian
$\GG$, denoted $\sL_\eta\dff \eta\inv \GO$.
It generates the   ``disc group'' orbit $\GG_\eta=L^{\ge 0}G
\sL_eta\sub\GG$.
Then
$L^{<0}G\cdot \sL_\eta
$ is
a normal slice to
$L^{\geq 0}G\sL_\eta$.

For the coweights $\la,\mu$ above one has 
$L^{\geq 0}G\sL_\la\sub\ \barr{ L^{\geq 0}G\sL_\mu }$
In the setting of matrices we view $\la,\mu$ as  partitions  of $N$,
so they provide  nilpotent orbits
$\OO_\la\sub\barr{\OO_\mu}$ in $\gl_N$.
We will denote by
$T_\la$
certain   ``{\em regular}'' normal slice
to $\OO_\la$.

The following theorem (announced in
\cite{MVyb}) is a common generalization of
(some of) the results of Kraft-Procesi \cite{KP}, Lusztig \cite{L81},
and Nakajima \cite{N94}.
The lower row contains an affine quiver variety and  intersections of
normal slices
(in nilpotent matrices and  the loop Grassmannian)
with larger  orbits.
The three   vertical maps are resolutions.
Here, the map $\widetilde{\barr{\OO_\mu}}
\CD @>{\bf m}>>\ \barr{\OO_\mu}\endCD$
is a certain Springer resolution
of the closure
of the nilpotent orbit $\OO_\mu$
and
$\pi: \widetilde\GG_\mu\to \barr{L^{\geq 0}G\cdot {\sL_\mu}}$ is a certain
``convolution map''.

\sss{\bf Theorem}\label{Theorem_iso!}
There exist natural algebraic isomorphisms
$\phi, \widetilde\phi, \psi, \widetilde\psi$
such that the following diagram commutes:
\begin{equation}\nonumber
\begin{CD}
{\text{\tiny $\begin{matrix}
Quiver\ varieties
\\
\ for \
A_{n-1}
\end{matrix}
$}}
@.
{\text{\tiny $\begin{matrix}
Nilpotent\ orbits
\\
\ for \
\gl_N
\end{matrix}
$}}
@.
{\text{\tiny $\begin{matrix}
Loop\ Grassmannians
\\
\ for\ GL_m
\end{matrix}
$}}
@.
\\
\\
\fMz(v,d)
@>{\widetilde\phi}>{\simeq}>
\tii{\barr{\OO_\mu}}|_{T_{\la}\cap\barr{\OO_\mu}}
@>{\widetilde\psi}>{\simeq}>
\tii\GG_\mu \tim_\GG\  L^{<0}G\cdot \la
\\
@V{p}VV
@V{\bf m}VV
@V{\pi}VV
\\
\fMz_\zero(v,d)
@>{\phi}>{\simeq}>
T_{\la}\cap\barr{\OO_\mu}
@>{\psi}>{\simeq}>
\barr{L^{\geq 0}G\cdot {\sL_\mu}}
\cap
L^{<0}G\sL_\la
.
\end{CD}
\end{equation}

\sss{
Deformation
along $Z[\gl(V)]$
}
The above theorem
corresponds to the case when the central element $c$ of
$\gl(V)$ is zero.
For arbitrary $c$
the nilpotent orbits deform to general
conjugacy classes, and the loop Grassmannian
deforms to the fiber
of the Beilinson-Drinfeld Grassmannian
$\BDG_\An$  at the point
$(0,c_1,c_1+c_2,\dots,c_1+\dots+c_{n-1})\in {\mathbb A}^{(n)}$
In this generality the theorem  is formulated as Theorem \ref{main}.

\sus{
The ``regular'' normal slice $T$ to a nilpotent orbit
}
The existence of isomorphisms
such as $\phi,\ \tphi$ -- using Slodowy slices -- was conjectured by Nakajima \cite{N94}
and established by Maffei \cite{M}.
In our  modifications  the normal slice
$T$ is \emph{not} the Slodowy's  slice -- it originates
from  the relation to
loop Grassmannian given by $\psi$.\fttt{
As we will see in
\ref{maximalpartition}
$T$ is more  generic or ``regular'' than the Slodowy slice.
}

Our transverse slice $T$ allows several new results.
First, the isomorphism $\phi$ is given by simple explicit formulas,
cf. \ref{constructphi}
as opposed to an inductive existence result used in \cite{M}.
This leads to an explicit formula for the isomorphism
$\psi\ci\phi$ between quiver varieties and pieces of the loop Grassmannians
in the appendix.
Second, we are able to
decompose a loop Grassmannian into a disjoint union of quiver varieties,
cf. \ref{decompositionAffGrass}.
Finally, our construction provides a natural environment for
the geometric $(GL_m, GL_n)$ duality,
cf. section \ref{applicationsRepTheory}.

\sex
In order to illustrate the difference, let us
give an example for $N = 5$ and a nilpotent element $x$ with Jordan blocks of sizes
$3$ and $2$.  If we fix the basis in which the matrix of $x$ has the Jordan canonical form,
then
the two
slices
are affine spaces
\begin{equation}\notag
\text{\rm Slodowy's slice} = \left (
\begin{smallmatrix}
a_1 & 1 & 0 & | & 0 & 0\\
a_2 & a_1 & 1 & | & b_1 & 0 \\
a_3 & a_2 & a_1 & | & b_2 & b_1 \\
\hline
c_1 & 0 & 0 & | & d_1 & 1 \\
c_2 & c_1 & 0 & | & d_2 & d_1 \\
\end{smallmatrix}
\right ),
\qquad
\text{\rm our slice} = \left (
\begin{smallmatrix}
0 & 1 & 0 & | & 0 & 0\\
0 & 0 & 1 & | & 0 & 0 \\
a_3 & a_2 & a_1 & | & b_2 & b_1 \\
\hline
0 & 0 & 0 & | & 0 & 1 \\
c_2 & c_1 & 0 & | & d_2 & d_1 \\
\end{smallmatrix}
\right ).
\end{equation}
In terms of a corresponding  $sl(2)$-triple
$\{x, h, y\}$
both slices are of the form $x+C$
for an $h,y$-invariant subspace $C$ complementary to $[\gl_N, x]$ in $\gl_N$
on which
$h$ acts with non-positive integral eigenvalues.
For Slodowy's slice
$C=Z_{\gl_N}(y)$,
so here $y$ acts by zero.
For our slice  the action of $y$ on $C$ is
\emph{``as close to regular nilpotent as possible"}, cf. \ref{maximalpartition}.


\sss{
Organization of the paper
}
Section \ref{quiversAsection} recalls quiver varieties of type A.
In section \ref{orbitSection} we define
the ``regular'' normal slices $T$ to nilpotent orbits.
In section \ref{grassmannSection} our slice
$T$ is related to the
Beilinson-Drinfeld family of loop Grassmannians.
In section \ref{mainresults}
we formulate the  Isomorphism Theorem and its consequences
in full generality
(with a deformation parameter $c$).
Section
\ref{quiverConjClassSection} proves relation between quiver varieties
and  matrices in a special case.
The general case is related to this special case
in section \ref{proofMainLemma}.
This allows us to finish
the proof of the Isomorphism Theorem in section \ref{proofIsomorphismTheorem}.
In section \ref{applicationsRepTheory} we give a geometric interpretation of $(GL_m,GL_n)$
dualities.
Finally, the appendix (section \ref{Appendix})
provides an explicit formula for the isomorphism
$\psi\circ\phi$  between quiver varieties and loop Grassmannians.

\sss{
Glossary of notation
}
As the paper compares three settings we provide here
some basic notations. The  relation between the
three kind of data is explained in \ref{data}.

{\bf (I)}
\ud{\em Quiver data} $n,v,d,c$.
Here $n$ ``means'' the quiver $A_{n-1}$.
Vectors $v,d\in\N^{n-1}$ are data
for two Nakajima quiver varieties related by a map
$p:\fMc(v,d)\to \fMc_\zero(v,d)$
(with image $\fMc_\one (v,d)$).
Here we view
$c\in Z[\pl\ \gl(V_i)]$ as a deformation parameter.

{\bf (II)}
\ud{\em Nilpotent orbits data $N,\la,\mu$.}
The partitions $\la,\mu$ of $N$
give nilpotent orbits
$\OO_\la\sub\barr{\OO_\mu}$ in
$\gl(D)$ for a vector space $D$ of dimension $N$.
By $T_\la=T_x$ we denote the ``regular''
normal slice to $\OO_\la$ at $x\in\OO_\la$.

{\bf (III)}
\ud{\em Loop Grassmannian data $m,\la,\mu$.}
Here $m$ ``means''
the group $G=GL_m$ with the loop Grassmannian $\GG=\GG(GL_m)$.
We will realize $GL_m$ as $GL(U)$ for a vector space $U$ with a basis
$\bfe_1,...,\bfe_m$. Then the loop Grassmannian of $G$ consists of
lattices in $U_\sK=U\ten \sK$ where
$\sK\dff\C((z))\ \subb\
\sO\dff\C[[z]]$ (see \ref{Loop Grassmannian and lattices}).
Any coweight $\eta\in\Z^m$ of $GL_m$ defines a lattice
$\sL_\eta=\pl\ z^{-\eta_i}\sO\bfe_i\in\GG$.

For a coweight $\si\in\Z^m$			
the orbit $\TT_\la=L^{<0}G\cd\sL_\la$
of the negative congruence subgroup is the standard  normal slice to $\GG_\la$ in $\GG$
at $\sL_\la$.

($\bu$) \ud{\em Resolution datum $a$.}
In the quiver setting
$\fMc(v,d)$ is a resolution of the image
$\fM^c_1(v,d)$ of
$
p^c:\fMc(v,d)\to\fMc_\zero(v,d)
$.

A decomposition $a\in\N^n$ of $N=\sum_1^n a_i$\
defines resolutions in matrices and loop Grassmannians.
It gives a partial flag variety $\FF^a$ of $GL(D)$
consisting of systems of subspaces $0=F_0\sub F_1\sub\cddd \sub F_n=D$
with $dim(F_i/F_{i-1})=a_i$. Its cotangent bundle
$\TN^a\dff T^*\FF^a$
is (see \ref{ginzburgfiber})\
a resolution of the closure $\barr{\OO_\mu}$
of a nilpotent orbit $\OO_\mu$ where
$\mu$ is the dual of the partition $\cmu$ obtained by
ordering the terms of $a=(a_1,...,a_n)$.

In the loop Grassmannian context, $a$ defines a convolution space
$\TGG^a_\mu
\dff \GG_{\om_{a_1}}\ast\cddd\ast\GG_{\om_{a_n}}$ for the fundamental coweights
$\om_i$. This is a
partial resolution
of $\barr{\GG_\mu}$.

($\bu$)
\ud{\em Deformation data $c,\sb$.}
A quiver datum $c\in Z[\pl\ \gl(V_i)]$  deforms
quiver varieties from $\fM^0(v,d)$
to
$\fMc(v,d)$.

The corresponding ``spectral datum'' $\sb\in\An$
deforms resolutions in matrices and in loop Grassmannians.
In matrices, $\sb$ deforms  resolution
$\TN^a\to\barr{\OO_\mu}$
of
the closure of a  nilpotent orbit
to a resolution of the closure of a  general
conjugacy class $\tgab\to\barr{\OO_{\EE,\tmu}}$
(see \ref{Deformation to a nilpotent orbit}).
Here, $\EE$ is the set of eigenvalues of any element $\sx\in\OO_{\EE,\tmu}$
and $\tmu$ is a family of partitions $\mu^e,\ e\in\EE$,
describing the nilpotent operators  $\sx-e$
on generalized $e$-eigenspace $D_e$ of $\sx$ in $D$.
In the loop Grassmannian $\sb$ deforms
the partial resolution
$\TGG^a_\mu
\to \barr{\GG_\mu}$.

\ud{\em Various notation.}
\label{BorelMoore}
By $G_m=\C^{\ast}$ we denote the multiplicative group.
For an algebraic variety $X$ whose connected components
$X_c$ are of pure dimension we denote
by $\HH(X)$ its top-dimensional Borel-Moore homology $\pl_c\
H^{\rm BM}_{2\dim(X_c)}(X_c)$.

{\bf Acknowledgments}
We are grateful to
A. Braverman, I. Frenkel, M. Finkelberg, D. Gaitsgory, V. Ginzburg, H. Nakajima, G. Lusztig, A. Maffei, A. Malkin, O. Schiffmann, and W. Wang
for useful discussions,
and to MSRI, IH\' ES and IAS for their
hospitality and support. The third author is especially grateful to his scientific advisor Michael Finkelberg for posing the problem and for many helpful suggestions.
The research of I.M. was
supported by NSF.
The research of M.V. was
supported by NSF Postdoctoral
Research Fellowship in 2001-2003.
The research of V.K. was partially supported by the HSE University Basic Research Program, Russian Academic Excellence Project '5-100' by the
August M\"obius contest (2016) and
the Dobrushin stipend.
We are very grateful to two referees for Advances who noticed many typos and suggested useful improvements in
presentation that influenced the paper.

\se{
Quiver varieties of type A
}
\label{quiversAsection}

\sus{
Representations of quivers
}
We start with the the Dynkin quiver $(I,\Om)$ of type $A_{n-1}$ with
vertices $I=\{1,\dots,n-1\}$ and the  arrows  $\Om$
given by $1\CD @>p_1>>\cddd @>p_{n-2}>>\endCD\ n-1$.

Passing to the Nakajima (framed) version
of the quiver involves first
doubling the arrows to $H=\Om\sqcup\overline\Om$
where
$\Om\con\overline {\Om},\ \om\mm\barr\om$, is the reversal of
orientation.
For an arrow $h\in H$ we denote
by $h'\in I$ its initial vertex and by $h''\in I$ its terminal vertex.

The data for framed quiver varieties are two
$I$-graded vector spaces
$V=\oplus_{i\in I}V_i$ and $D=\oplus_{i\in I}D_i$.
Their dimension vectors
$v,d\in\N^I$
define a vector  space
\[
M(v,d)\ =\ \bigoplus_{h\in H}\Hom(V_{h'},V_{h''})\oplus
\bigoplus_{i\in I}\Hom(D_{i},V_{i})\oplus
\bigoplus_{i\in I}\Hom(V_{i},D_{i}).
\]
Following Lusztig and Maffei
we will consider an element in $M(v,d)$ as a quadruple
$(x,\bx,p,q)$ with

\begin{equation}
\begin{split}
&
x
=
(x_h)_{h\in \Om} \in\bigoplus_{h\in \Om}\Hom(V_{h'},V_{h''}),
\ \ \
\barr x
=
(x_h)_{h\in \barr\Om} \in\bigoplus_{h\in \barr\Om}\Hom(V_{h'},V_{h''}),
\\
&
p =
(p_i)_{i\in I} \in\bigoplus_{i\in I}\Hom(D_{i},V_{i}),
\ \ \  \ \ \ \ \
q =
(q_i)_{i\in I} \in\bigoplus_{i\in I}\Hom(V_{i},D_{i}).
\end{split}
\end{equation}



The total of the $x$-data is denoted $\tx=(x,\bx)=(\tx_h)_{h\in H}$.
For any $j>i$ we define the following polynomial operators
in $(\tx,p,q)$ :
$$
p_{j \to i}
\dff
\bx_{i} \dots \bx_{j-1} p_j:D_j\to V_i
\aand
q_{i \ra j} \dff
q_j x_{j-1} \dots x_{i}:V_i\to D_j
.$$

\sss{}
\lab{The group G(V)=prod GL(Vi) acts on M(v,d)}
The group $G(V)=\prod_{i\in I}GL(V_i)$ acts on $M(v,d)$
so that for
$g=(g_i)_{i\in I}$
\begin{equation}
g(x,\bx,p,q)
\ \dff\
(g_{i+1}x_ig_i^{-1},
g_i\bx_ig_{i+1}^{-1},
g_ip_i,
q_ig_i^{-1})_{i\in I}.
\end{equation}

We are interested in fibers of the
corresponding moment map
\begin{equation*}\mmm\colon M(v,d)\to \fg(V)\dff Lie[G(V)]\cong \pl \gl(V_i),
\end{equation*}
at elements
$c=(c_1,\dots,c_{n-1})\in
Z[\gl(V)]$.
The fiber
$\Mc(v,d)\dff
\mmm^{-1}(c)$
consists
of all $(x,\bx,p,q)$ such that
\begin{equation}
\label{deformedrelations}
\begin{split}
c_i+\bx_ix_i =x_{i-1}\bx_{i-1}+p_iq_i
\ \ \
\foor
2\leq i\leq n-2;
\aand
\\
c_1+\bx_1x_1 =p_1q_1,\ \ \ \ \ \
c_{n-1} = x_{n-2}\bx_{n-2}+p_{n-1}q_{n-1}
.
\end{split}
\end{equation}

\sus{
Generators of $G(V)$-invariant polynomials on $ \Mc(v,d)$
}
The following is version of a result of Lusztig
using some ideas of Maffei.

\theo{}
\label{invtheorem}
The algebra of invariant functions
$\OO(\Mcvd))^{G(V)}$
is generated by polynomials
of the form $\chi(q_{l\ra j}p_{i\ra l})$, where
$1\le l\le i,j< n$,
and $\chi$ is a linear form on
$\Hom(D_i, D_j)$.

\begin{proof}
Lusztig~\cite[1.2]{L98} considered two kinds
of invariant polynomials in $\OO(\Mvd)^{G(V)}$.

(a) Any cycle  $\ga=(h_1,h_2\dots,h_r)$
in $(I,H)$ defines a polynomial
$
\Tr(V_{h_1'},\tx_{h_r}\tx_{h_{r-1}}\dots \tx_{h_1})$.

(b) A path  $
(h_1,h_2\dots,h_r)$ in $(I,H)$
and a linear form $\chi$ on
$\Hom(D_{h_1'}, D_{h_r''})$ define a polynomial
$
\chi(q_{h_r''}\tx_{h_r}\tx_{h_{r-1}}\dots \tx_{h_1}p_{h_1'})
$.

He proved
\label{LusztigTheorem}
\cite[Theorem 1.3, 5.8]{L98}
that the algebra $\OO(\Mc(v,d))^{G(V)}$
is generated by restrictions of invariant polynomials of types (a) and (b)
above.

It just remains to switch from Lusztig's
generators to Maffei's generators \cite{M}
using the following lemma.
\end{proof}

\slem{}
Let $\PP$ be the space of operators on $D$ which are
polynomials  in products
$q_{l\ra j}p_{i\ra l}$ \ for\  $1\le l\le i,j\le n-1$,

(a)
For  any cycle
$h_1,\dots,h_r$ in $(I,H)$
the restriction of
the polynomial
$
\Tr(\tx_{h_r}\tx_{h_{r-1}}\dots \tx_{h_1})
$
to $\Mc(v,d)\sub M(v,d)$\
is of the form $\Tr(P)$
for some  $P\in\PP$.

(b) For any path $h_1,\dots,h_r$ in $(I,H)$
and  a linear form $\chi$ on
$\Hom(D_{h_1'}, D_{h_r''})$
the restriction of
the polynomial
$
\chi(q_{h_r''}x_{h_r}\cddd x_{h_1}p_{h_1'})
$
to $\Mc(v,d)$
is of the form
$\chi(P)$
for some  $P\in\PP$.

\begin{proof}
Easily follows from equations
(\ref{deformedrelations})
for
 $\Mc(v,d)\sub M(v,d)$.
\end{proof}

\sus{
Nakajima's quiver varieties \cite[3.12]{N98}
}
\label{quiverdefs}
The affine invariant theory quotient of
$\Mc(v,d)$ by $G(V)$ is denoted
\begin{equation}
\fMc_\zero\vd
\ \dff\
\Mc(v,d)//G(V)
\ =\ \Spec[ \OO(\Mc(v,d))^{G(V)}].
\end{equation}

\sss{The stable part
$\Mc_s(v,d)\sub \Mc(v,d)$
}
Following Nakajima \cite{N98} and Lusztig \cite[2.11]{L98}
we say that a quadruple
$(x,\bx,p,q)$ is
{\em stable} if for any $I$-graded subspace $V'$ of $V$
preserved by
$x$ and $\bx$
and containing
$\Im (p)$, we have $V'=V$. The subset of all stable quadruples in
$\Mc(v,d)$ is denoted by
$\Mc_s(v,d)$.

\slem{}\label{stablelemma}
\cite[Lemma 14]{M}
An element  $(x,\bx,p,q)$ of $\Mc(v,d)$ is stable if
and only if
for all $1\leq i\leq n-1$
\begin{equation}
\Im(x_{i-1})+\sum_{j=i}^{n-1} \Im(p_{j\ra i})\ =\ V_i.
\end{equation}

\sss{
The quiver variety $\fMc(v,d)$
}
This is the geometric quotient of the stable part $\Mc_s(v,d)$ by $G(V)$.
In particular the set of $\C$-points
is the quotient set $\Mc_s(v,d)/G(V)$.
Below we only consider such $(v,d)$ that $\fMc(v,d)$ is nonempty
(for explicit conditions on $(v,d)$ see \cite[10]{N98}, \cite[Lemma 7]{M}).

\sss{
$\fMc_\one$ and the Lagrangian $\fL(v,d)$
}
There is a canonical  map $\bbp:\fMc(v,d)\ra \fMc_\zero(v,d)$.
Following Maffei we denote its image
$$
\fMc_\one (v,d)\ \dff\ \Im(\bbp)\ \subset\fMc_\zero(v,d).
$$
Finally, we consider the  central fiber of $\bbp$ which is a Lagrangian
$\fL^c(v,d)\dff \bbp^{-1}(0)$ in $\fMc(v,d)$,  and
its  top-dimensional Borel-Moore homology
$
\HH(\fL(v,d))
$ (see \ref{BorelMoore}).

\sus{
Nakajima's construction of
$SL(n)$-modules
}
\label{slnmodules}
In this subsection $c=0$.

\stheo{}
\label{theonakajima}
\cite[\emph{10.ii}]{N98}
For any $d\in\N^{n-1}$ the sum  $\pl_v \HH(\fL(v,d))$ has the structure of a
simple $SL(n)$-module $V^{SL(n)}_d$ with the highest weight $d$
(meaning really $\sum_1^{n-1}\ d_i\om^i$ for the  fundamental weights
$\om^i$).
The summand
$\HH(\fL(v,d))$
is the weight space for
the weight
$d-Cv$, where $C$ is the Cartan matrix of type $A_{n-1}$.

In particular, the  weight space
$V^{SL(n)}_d(d-Cv)$ has a basis $\Irr\fL(v,d)$
of
irreducible components of $\bbp^{-1}(0)$.
Lusztig \cite{L00} calls this basis \emph{semi-canonical}.

\sss{
From $SL(n)$ to $GL_n$
}
\label{fromSLntoGLn}
\label{sltogl}
We may consider $\pl_v \HH(\fL(v,d))$ as a representation $V^{GL(n)}_{\cla}$
of $GL_n$
with highest weight $\cla$,
where $\cla=\cla(d)=(\cla_1, \cla_2,\dots,\cla_n)$
is a partition of
$N=\sum_{j=1}^{n-1}jd_j$ defined as follows:
$\cla_i=\sum_{j=i}^{n}d_j$ (here $d_n=0$). Then
$\HH(\fL(v,d))$ is the weight space $V^{GL(n)}_{\cla}(a)$,
where $a_i=v_{n-1}+\sum_{j=i}^{n}(d-Cv)_j$ (here $(d-Cv)_n=0$),
cf. \cite[8.3]{N94}.

\se{
Conjugacy classes of matrices
}
\label{orbitSection}
In this section we fix a vector space $\sN$ of dimension $\sn$,
let
$G = GL(\sN)$ and $\fg=\gl(\sN)=\End(\sN)$.
In \ref{definitionBases} we recall degenerations of
conjugacy classes in  $\gl(\sN)$
and the corresponding resolutions.
In
\ref{ontransverseslices}
we
consider normal slices to nilpotent orbits.

\sus{
Resolutions and degenerations of
closures of conjugacy classes
}
\label{definitionBases}
For an operator $\sx\in\End(\sN)$ we list
invariants $\EE_\sx,\tmu_\sx,\mu_x,M_\sx$
of the conjugacy class
$\OO_\sx\dff\ ^{GL(\sN)}\sx$. We use $(\EE_\sx,\tmu_\sx)$
to describe $\OO_\sx$.
An ordering $(a,\sb)$ of $M_\sx$
gives a resolution
$\tgab$ of $\barr{\OO_\sx}$ and a degeneration of $\barr{\OO_\sx}$
to $\barr{\OO_{\mu_\sx}}$ for the nilpotent orbit
attached to the partition $\mu_\sx$.
We also recall two related constructions of representations of $\gln$.

\sss{
The variety $\FF_n$ of $n$-step flags
}
Let $\NN$ be the nilpotent cone in $\fg=\End(\sN)$.
The variety $\FF_n=\FF_n(\sN)$ of {\em $n$-step flags} in $\sN$
consists of systems of subspaces
$F=(0=F_0\subseteq F_1\subseteq \dots\subseteq F_n=\sN)$
(it was studied by Ginzburg \cite{CG}).
The  $n$-term decompositions
$a=(a_1,\dots,a_n)\in\N^n$ of
$\sn=\sum_{i=1}^{n}a_i$),
parameterize the
connected components
$\FF^a_n=\{F\in\FF_n;\ \dim(F_i/F_{i-1})=a_i\}$\
of $\FF_n$.

\sss{
Fibers of moment maps
}
\label{ginzburgfiber}
The
cotangent bundle
$\TN_n\dff T^*\FF_n$
consists of all
$(x,F)\in\fg\times\FF_n$ such that $xF_i
\subseteq F_{i-1}
$. It lies inside
the ``extended cotangent bundle''
$\tfg_n\dff\tii T^*\FF_n$ given by all $(x,F)$
with $xF_i\sub F_i$.
The moment map for the $G$-action
on $\tfg_n$ is the
projection
$\tmm\colon \tfg_n\to \fg\dff\fg$.
Denote its restriction by
$\mmm\colon \TN_n\to \fg=\End(\sN)$.
Its
image
$\NN_n\dff \mmm(\TN)$
is given by the equation $x^n=0$.
We also denote
$\TN^a_n= T^*\FF^a_n$
and in the extended version
$\tfg_n\dff\ \tii T^*\FF_n$
with
$\tfg_n^a\dff \tii T^*\FF^a_n$.
The restrictions of $\mmm,\tmm$ to
$\TN_n^a,\tfg_n^a$
are denoted
$\mmm^a,\tmm^a$.


\slem
$\TN_n^a$
is a resolution of the closure of the nilpotent orbit $\OO_\mu$
where the dual partition $\cmu$ is given by
reordering
$a=(a_1,...,a_n)$.


\begin{proof}
First we notice  equality of dimensions
$$\dim(\TN^a)
=
2\dim(T^*\FF^a)
=
\sum_{i\ne j} a_ia_j
=
{\sn}^2-
\sum_{i=1}^n  a_i^2
=\
\sn^2-\sum_{i=1}^n |\cmu_i|^2
=
\dim(\OO_\mu)
.
$$
Since $\TN^a$ is irreducible, it suffices now to show that
$\TN^a$ covers $\OO_\mu$ isomorphically.

Let $y\in\OO_\mu$.
Any Jordan basis $\JJ$ of $\sN$ for $y$, corresponds to  squares in the Young diagram 
of $\mu$  so that
the  span of the $i\thh$ column is the $i\thh$ Jordan block.
Let $V_p$ be the span of the $p\thh$ row $\JJ_p\sub\JJ$
of $\mu$ and let $I_p=y^{p-1}V_p$. Then $y^{p-1}:V_p\con I_p\sub V_1=\Ker(y)$.

Let $a_1$ be the length $\cmu_p$ of the $p\thh$ row of $\mu$.
Then there is precisely one subspace $F_1$ of $\Ker(y)$ such that the Young diagram of
$y$ on $\sN/F_1$ is the diagram $\mu'$ obtained from $\mu$ by erasing the $p\thh$-row.
This is clearly satisfied for the subspace $F_1=I_p$.
For uniqueness,  observe that  for the highest row (on the height $\mu_1$),
$F_1$ must contain $I_{\mu_1}$ since otherwise one of Jordan blocks of $\sN/I_p $
would have length $\mu_1$ which is too large for $\mu'$. This  argument
can be iterated to show that $F_1$ contains all $I_j$ for $\mu_1\ge j \ge p$.

For any  
reordering $a=(a_1,...,a_n)$ 
of $\cmu$
we can  repeat this procedure for $y$ on $\sN/F_1$ etc,
to
construct 
a unique filtration $F_1\sub\cddd\sub F_n$  such that $(y,F)\in \TN^a$.
(For example, if
$a=\cmu$ then
simply $F_i=\Ker(x^i)$.)
\end{proof}

\sss{
Data $(\EE,\tmu)$
for conjugacy classes in
$\End(\sN)$
}
\label{conjclass}
To an operator
$\sx$ on $\sN$
we assign the set of its eigenvalues $\EE=\EE_\sx\sub\ \A^1$
and a collection
 $\tmu_\sx=\tmu=(\mu^e)_{e\in \EE}$ of partitions $\mu^e$ corresponding
to the nilpotent operator $\sx-e$
on  the generalized $e$-eigenspace
$\sN_e$ of $\sx$.
We will call such $\tmu$ an
{\em $\EE$-bipartition} of $\sn$ since $\sum_{e\in \EE}\ |\mu^e|=\sn$.

The data $(\EE,\tmu)$ determines the conjugacy class
$\OO_x$ through $\sx$ so we also  denote it $\OO_{\EE,\tmu}$.
This parameterizes all conjugacy classes
as $\EE$ goes through subsets of $\Ao$ of $\le n$ elements
and $\tmu$ through $\EE$-bipartitions
of $\sn$.


\sss{
Invariants $n_\sx,\mu_\sx,M_\sx$ of
an operator $\sx$
}

A bipartition $\tmu=(\mu^e)_{e\in\EE}$
defines the dual  bipartition
$\tmu^\sv
\dff\
((\mu^e)^\sv)_{e\in\EE}$.
It also defines a single partition
$\mu$ of $\sn$ such that its dual $\ch\mu$
is the   partition
obtained by ordering the multiset
of terms in the bipartition
$\tmu^\sv
$.\fttt{
A {\em multiset structure} on a set $A$ is
a multiplicity function
$A\to\N$\ie
an isomorphism type of
surjections $A_1\sur A$ with finite fibers.
An {\em ordering} on the multiset $A_1\sur A$
is an ordering on $A_1$.
}

For
any eigenvalue $e$ of $\sx$  denote by $n_e$ the maximal size
of  Jordan blocks for $e$\ie
the number of steps in the filtration of $\sN_e$ by $\Ker[(\sx-e)^i]$.
This
is also the number of terms in  the dual partition
$(\mu^e)^\sv=(\cmu_1^e\geq \cddd \geq \cmu^e_{n(e)}>0)$.
The sum
$n_\sx
=\
\sum_{e\in \EE} n_e
$ is the total number of terms in
the bipartition $(\tmu_\sx)^\sv$
and it is also
the degree of the minimal polynomial of $\sx$.

Let us say that a filtration
 $F=(0=F_0\sub F_1\sub\cddd\sub F_n=\sN)$
is a {\em Jordan filtration} for $\sx$
if there is an ordering $\sb=(b_1,...,b_{n_x})$ on the multiset
$\tii\EE_\sx$ of
roots of the minimal polynomial of $\sx$,
such that $F_i=\Ker[\prod_{1\le j\le i}\ (\sx-b_j)]$.
(The length of any Jordan filtration is $n_\sx$.)
The graded piece $F_i/F_{i-1}$ is described by
its dimension $a_F(i)=\dim(F_i/F_{i-1})$
and by the scalar $b_F(i)=b_i$ (the action of $\sx$).


%


We denote by $M_\sx$.
 the multiset
of pairs $(a_F(i),b_F(i)), 1\le i\le n_\sx$,
for any Jordan filtration $F$ (it
is independent of $F$).

The multiset
$A_\sx$ is given by all  first components of elements of $M_\sx$\ie
by dimensions of graded pieces of any Jordan filtration.
Notice that if $\sx$ lies in the orbit $\OO_{\EE,\tmu}$ then $M_\sx$  can be described as
the multiset of   terms in
the dual bipartition $\tmu^\sv\dff\
((\mu^e)^\sv)_{e\in\EE}$ of $\sn$.

\sss{
Resolutions  and
degenerations
}
\label{Deformation to a nilpotent orbit}
For an operator $\sx$ on $\sN$
define $\EE,\tmu$ and $M_\sx$ as above so that $\OO_\sx=\
\OO_{\EE,\tmu}$.
An {\em ordering}  $(a_1,b_1;...;a_n,b_n)$
of the multiset $M_\sx$
will be viewed as  a pair of vectors
$(a,\sb)=(a_1,...,a_n;b_1,...,b_n)$
(hence $b_i\in\EE$ and $a_i\in\N$).
The component  $a$ is an ordering
of
the multiset $A_\sx$ of dimensions of graded pieces in
a Jordan filtration
while
$\sb=(b_1,...,b_n)$ remembers the eigenvalues of $\sx$ on
the graded pieces.
Notice that any Jordan  filtration $F$ gives
an ordering $(a_F,\sb_F)$ and all orderings arise in this way.

The first component $a$
defines
$\TN^a=T^*\FF^a
\sub \tfg^a=\tii T^*\FF^a$
and then $\sb$ defines
\begin{equation}
\label{gnabth}
\tgab\
\dff\
\{(\sx,F)\in \tfg^a;\
\sx \text{ acts on } F_i/F_{i-1}
\text{ as }
b_i
\}.
\end{equation}

\slem{}
(a)
The variety $\tgab$ is smooth and  connected.
Its moment map image is
$\barr{\OO_{\EE,\tmu}}$
and the moment map $\tgab\to\barr{\OO_{\EE,\tmu}}$
is a resolution.

(b) For $\sb\in\A^n$ the maps
$\tgab\to\barr{\OO_{\EE,\tmu}}$ form a deformation
of $\TN^a\to \barr{\OO_\mu}$\ which is the case $\sb=0$.\fttt{
Here $\OO_{\EE,\tmu}$ degenerates to $\OO_\mu$ so that in each Jordan block
eigenvalue degenerates to zero but extensions are created between
blocks corresponding to different eigenvalues.
An example is a regular semisimple orbit degenerating to
the regular nilpotent orbit.
}

(c)
For the Lusztig stratum $\LL_\sx$ that contains
$\OO_\sx=\OO_{\EE,\tmu}$,\ one has
$\barr{\LL_\sx}\cap\NN=
\barr{\OO_{\mu_\sx}}$.
\fttt{
The Lusztig stratification of $\fg$ is the stratification
induced by the $G$-action.
The stratum through $\sx$ with a Jordan decomposition
$\sx=s+n$
is $\LL_x\dff\ ^G(Z_r(\fl)+n)$
where $\fl\dff Z_\fg(s)$
and $Z_r(\fl)=\ \{y\in\fl;\ Z_\fg(y)=\fl\}$.
(cf. \cite[5.5]{Mir} and references therein.)
}

%
%

\begin{proof}
(a) 
Variety $\tgab$
is a torsor for the vector bundle $T^*\FF^a$ over $\FF^a$,
so it is smooth and connected.
The claim that the fiber at $\sx\in\OO_{\EE,\tmu}$ is a point
reduces to the case when $\sx$ has one eigenvalue
since $\sx$-invariant filtrations $F$ decompose
into an $\EE$-sum of filtrations on generalized eigenspaces,
Then $\sx$  may as well be nilpotent and this case is the lemma
\ref{ginzburgfiber}.

The remaining claims follow from equality of dimensions of
$\tgab$ and $\OO_{\EE,\tmu}$. Since
we know that $\dim(\tgab)=\dim(\TN^a)=\dim(\OO_\mu)$,
we just need
$$
\dim(\OO_{\EE,\tmu})=\dim[\gl(\sN)-\dim[Z_{\gl(\sN)}(\sx)]
=\
\sn^2-\sum_{e\in\EE}\ \dim[Z_{\gl(\sN_e)}(\sx-e)]
$$
$$
=
\sn^2-\sum_{e\in\EE}\ \sum_i\ |\cmu^e_i|^2
=\dim(\OO_\mu)
.$$

(b)
The family of $ T^*\FF^a$ torsors $\tgab $ is clearly flat.
Here, the deformation space is
$$
\tfg^{a,\A^n}
\dff\
\{(\sx,F)\in \tfg^a;\
\sx \text{ acts on } F_i/F_{i-1}
\text{ as a scalar}
\}=\ \cup_{\sb'\in\An}\ \tfg^{a,\sb'}
.$$
At the same time the moment map image
$
\barr{ \OO_{\EE,\tmu}}
=\tmm(\tgab)
$
moves to $\mmm(\TN^a)=\barr{\OO_\mu}$
and
$\dim(\OO_\mu)=\dim(\TN^a)
=\dim(\tgab)
=\dim( \barr{\OO_{\EE,\tmu}} )$.

(c)
Notice that $\LL_\sx$ is the union of all
$\barr{\OO_{\io(\EE),\io_*\tmu}}$ where $\io:\EE\inj\Ao$
is injective.
The space 
$
\tfg^{a,\A^n}
$
is closed in $\tfg^a$ and its   moment map image
$\tmm^a(\tfg^{a,\A^n})$
is by (a) the closure
of $\barr{\LL_\sx}$.
Now,
$\tmm(\tfg^{a,\A^n})\cap \NN=
\tmm(\tfg^{a,0})
=\tmm(\TN^a)
=
\barr{ \OO_{\mu_\sx}}$.
\end{proof}

\sss{
Two geometric bases of
representations
of $\gln$
}
At $x\in \NN_n$ we denote the fibers of moment map
by $\FF_{n,x}^a\dff \mmm_a^{-1}x\sub\FF_n^a$
and
$\tii\FF_{n,x}^a\dff \tmm_a^{-1}x\sub\tfg_n^a$.
Let also $\HH(Z)$ denote the  top Borel-Moore homology of $Z$.

\stheo{}\label{theoginzburg}
Let  $\la$ be a partition of $\sn$ with at most $n$ terms
and let $x,y$ be nilpotent operators on $\sN$ of types
$\la$ and $\cla$.
Then the irreducible
$\gl(n)$-module $V_{\la}$ with the highest weight $\la$ has realizations as
$\HH(\FF_{n,y})$
(see \cite[4.2]{CG})
and as $\HH(\tii\FF_{n,x})$
(see \cite{BG}).

\sus{
Normal slices to  nilpotent orbits
}
\label{ontransverseslices}

We will say that a normal
(transverse)
slice in $\fg$ to
a nilpotent orbit $\al$ at a point $x\in\al$,
is a submanifold $\SS$ of $\fg$ such that
\bi
({\bf IN})
({\em Infinitesimal normality})\
The tangent space
$\sT_x(\fg)$ equals
$\sT_x(\al)\ \pl\  \sT_x(\SS)$.
\item
({\bf C})
({\em Contraction.})
There is an action of $G_m$ on $\SS$ which
contracts it to $x$ and preserves intersections with
the Lusztig strata
in $\fg$.\fttt{
The Lusztig stratification of $\fg$ is the stratification
induced by the $G$-action.
The stratum through $\sx$ with a Jordan decomposition
$\sx=s+n$
is $\LL_x\dff\ ^G(Z_r(\fl)+n)$
where $\fl\dff Z_\fg(s)$
and $Z_r(\fl)=\ \{y\in\fl;\ Z_\fg(y)=\fl\}$.
(cf. \cite[5.5]{Mir} and references therein.)
}
\ei

\slem
For a normal slice $\SS$ at an element
$x$ of a nilpotent orbit $\al$ :
\begin{enumerate}
\item $\SS\cap \al=\{x\}$;
and $^G\SS$ is open in $\fg$;
\item $\SS$ meets Lusztig stratum $\be$ iff $\al\sub \barr\be$;
\item $\SS$ meets conjugacy classes transversally;
\i
for any partial flag variety $\PP=G/P$ the base changes
$\tii T^*\PP\tim\fg\ \SS$
and $T^*\PP\tim\fg\ \SS$
of the slice $\SS$ are smooth and connected.
\end{enumerate}

\begin{proof}
(1)
Since nilpotent orbits are Lusztig strata, $G_m$ contracts
$\al\cap \SS$ to $x$. So if $\al\cap \SS$ were more then a point then
it would have positive dimension contradicting property ({\bf IN}).
Next,
$^G\SS$ is open at the point $x$ by
$({\bf IN})$ and therefore at each  $y\in \SS$ by
({\bf C}).

(2)
 If $\SS$ meets $\be$ then the contraction of $\barr{\be\cap \SS}$
to $x$ shows that $\barr\be $ contains $x$ and $\al$.
On the other hand, if $\al\sub\barr\be$ then
the neighborhood $^G(\SS\sub\fg)$ of $\al$ meets $\be$.

(3) Our  transversality claim for a subspace and a foliation 
means that for  $y\in \SS$, the sum
$\Si_y\dff \sT_{^y}(^Gy)+\sT_{y}(\SS)$  is all of $\fg$.
For this notice that
as $u\in \Gm$ approaches $0$, we have $z=\ ^uy\to x$, hence
$\dim(\Si_y)=\dim(\Si_{^uy})\ge \dim(\Si_x)=\dim(\fg)$.

(4)
Near $y\in \SS$, $\fg$ is a product of
$\SS$ and the orbit $G\cd y$, so smoothness is inherited from
that of $T^*\PP$ and $\tii T^*\PP$.
The $G_m$-action lifts to
$\tii T^*\PP\tim\fg\ \SS\ \subb\
T^*\PP\tim\fg\ \SS$ and contracts the spaces to the fiber at $x$.
Connectedness then follows from connectedness of (generalized)
Springer fibers cf. \cite{Sp}.
\end{proof}

\slem
\label{sufficientnormaldata}
Consider  a  vector subspace $C\sub \fg$ complementary to
$\sT_x(\al)=[\fg,x]$.
For the submanifold $x+C$
to be  a normal slice to the orbit $ ^Gx$ at $x$
it suffices that there exists $h\in\fg$ which
is\ (i) semisimple,\ (ii)
integral (i.e., eigenvalues of $\ad\ h$ are integral),\
(iii)  $[h,x]=2x$,\
(iv) $h$ preserves
$C$
and\
(v) the eigenvalues of $h$ in $C$ are
${\le 1}$.

\begin{proof}
Such $h$ lifts to a homomorphism $\io: G_m\to G$.
Now the action of $G_m$  on $\fg$ by
$
s\ast \sy=\ s^{2}\ \cd\ ^{\io(s)\inv}\sy, \ s\in G_m,\ \sy\in\fg
$;\ \
preserves $x+C$ and contracts it to $x$.
\end{proof}

\sss{
Normal slice $T_x$
}
\lab{Normal slice Tx}
\label{basisfreeslice}
For  a nilpotent
$x$  on $\sN$ we will use an  $sl(2)$-triple $x, h, y$
to construct a pair $(h,C)$ as in the lemma
\ref{sufficientnormaldata}.\fttt{
In practice we start with a  choice of a minimal subspace $\CC$ generating  $\sN$ under $x$,
and a decomposition of $\CC$ into lines. This gives a canonical $sl_2-$triple.
}
The subspaces $C$  will be of a special form:\
(i) $C$ is  $y$-invariant,\
(ii) the $h$-eigenvalues in $C$
are $\leq 0$.
For instance the  Slodowy slice corresponds to the case
$C = Z_{\fg} (y)$.

We will use the  $sl(2)$-decomposition of the vector space $\sN$
as a sum\
$\sN = \bigoplus_i M_i\otimes L_i$
where $L_i=L_i^{sl(2)}$ is a simple $sl(2)$-module of highest weight $i$, and $
M_i=\Hom_{sl(2)}(L_i,\sN)$.
Then
$\fg=\End(\sN)$ is $\bigoplus_{i, j} \Hom(M_j, M_i)\otimes L_j^{\ast} \otimes L_i,
$ and we
choose
\begin{equation}
C = \bigoplus_{i, j}\ \Hom(M_j, M_i)\otimes
\Ker(y,L_i)
\otimes
\Ker(y^{i+1},L_j^{\ast})
\subseteq \End(\sN).	
\end{equation}
Notice that
$C$ generates $\End(\sN)$ as a $\C[x]$-module
and it is a minimal subspace with this property
(it is easy to see that this holds for the subspace
$
\Ker(y,L_i)
\otimes
\Ker(y^{i+1},L_j^*)
$ of $L_i
\otimes
L_j^*
$).
This property implies  that $C$ is complementary to $[\fg,x]$.


It is elementary to see that $h, C$ satisfy the conditions of Lemma \ref{sufficientnormaldata}, and thus
$T_x = x + C$ is a normal slice.

\sss{
Normal slice to a pronilpotent}
\lab{pronilpotent}
Let $z$ be regular nilpotent operator on the vector space $\OO$ with a cyclic line  $\LL$.
On the multiple  $\sN=U\ten \OO$ of $\OO$ we consider the operator
$x=1\ten z$.
Then the above formula for the subspace $C$ is
$$
C=\ \End(U)\ten \LL\ten \OO^*\ \cong\ \Hom(U\ten\OO,U\ten \LL)
.$$
Notice that this gives a normal slice
$T_x=x+C$ even if $e$ is pronilpotent, say the $z$-multiplication
on $\OO=\C[[z]]$.

\sss{
``As close to regular nilpotent as possible"
}
\label{maximalpartition}
\lab{regular slice}
If $x$ is a nilpotent of type
$\la
$
then
$\sN \cong  \bigoplus_{i=1}^m  L_{\la_i-1}
$ hence
$\End(\sN) = \bigoplus_{i, j = 1}^m
L_{\la_i-1}
\otimes
L_{\la_j-1}^{\ast}
$.
Then $C = \bigoplus_{i, j = 1}^m\
\Ker(y^{\la_i},L_{\la_j-1}^{\ast}) \otimes \Ker(y,L_{\la_i-1})
$.

As a $y$-module
$\Ker(y^{\la_i},L_{\la_j-1}^{\ast}) \otimes \Ker(y,L_{\la_i}-1)
\ \cong\
\Ker(y^{\la_i},L_{\la_j-1}^{\ast})
\ \cong\ L_{\min(\la_i,\la_j)-1}$ is a single Jordan block of size
$\min(\la_i,\la_j)$.
So,  if $\la=\ 1\hp{m_1}\cddd s\hp{m_s}$
then  the
$y$-action  on $C$
contains
the Jordan block of size $i$
with multiplicity
$2(m_i+\cddd+m_s)-1$.
Then the partition
$\la_C$ corresponding to $y$ on $C$
is the largest possible partition
for the action of $y$ on any $y$-invariant subspace of $\fg^h_{\leq 0}$
complementary to $\sT_x(\al)=[\fg,x]$.
(For instance,  if $x$ is regular, then $y$ restricted to our $C$ will
be regular.)
This is opposite from  Slodowy's situation where $y$ acts
on $C = \Ker(y,\fg)$ as $0$.
For this reason we will call $T_x$ the ``{\em regular}'' normal slice.


\sss{
Slice $T_x$ in a Jordan basis
and $C\sub\Hom(\sN,\U)$
}
\label{sliceJordanBasis}
\label{nilpnormalslice}
Again,   here $x$ is
a nilpotent operator of type $\la=(\la_1\geq\dots\geq \la_m)$
on $\sN$.
Let us choose a Jordan basis
$\es_{k,i}$ for $1\leq k\leq \la_i$ of $\sN$\ie
such that  $x \es_{k,i}= \es_{k-1,i}$.
Let $\U\sub \sN$ be the span of the
generators $\es_{\la_i,i},\ 1\le i\le m$ of Jordan blocks.
Then  $C$ is a subspace of $\Hom(\sN,\U)$
consisting of all $\ga\in\Hom(\sN,\U)$
such that
for the matrix elements
$\ga_{l,j}^{i} $
corresponding to
the component
$\C\es_{l,j}\to\C \es_{\la_i,i}$,\
one has
$\ga_{l,j}^{i} =0$ for $ l>\la_i$.
For example,
if $\la=(\la_1\geq\la_2) = (3, 2)$ the operators in
$T_x$ will have the form
(in the
basis $\es_{k,i}$)
\begin{equation}\label{tsliceexample}
\left(
\begin{matrix}
0 & 1 & 0 & | & 0 & 0 \\
0 & 0 & 1 & | & 0 & 0 \\
\ga_{1, 1}^{1}
&
\ga_{2, 1}^{1}
&
\ga_{3, 1}^{1}
& |&
\ga_{1, 2}^{1}
&
\ga_{2, 2}^{1}
\\
\hline
0 & 0 & 0 & | & 0 & 1 \\
\ga_{1, 1}^{2}
&
\ga_{2, 1}^{2}
& 0 & | &
\ga_{1, 2}^{2}
&
\ga_{2, 2}^{2} \\
\end{matrix}
\right ).
\end{equation}

\sss{
Intersection of the normal slice $T_x$ with
the conjugacy class $\OO_{\sy}$
}
\label{tildeta}
Since $x$ is of type $\la$ the slice $T_x$ will also  be denoted $T_\la$.
Consider $\barr{\OO_\mu}\subb\OO_\la$
and let $a$ be any permutation of terms in $\cmu$.
The restriction of the moment map $\mmm^a:\TN^a\to\NN$ to the slice
$T_x$ is the space
$
\tT^a_x\dff\
T_x
\
\tim_{\NN}
\
\TN^{a}
$.


\slem
Consider an element $x$ of a nilpotent orbit
$\OO_\la$.

(a)
Let
$\OO_\la\sub\barr{\OO_\mu}$ and
let $a=(a_1,\dots,a_n)$ be
related to $\mu$ as above.
Then the variety
$\tT_{x}^{a}$ is smooth and connected
with the moment map image $T_x\cap\barr{\OO_\mu}$.
The map
$\mmm_a: \tT_{x}^{a}\to T_{x}\cap \barr{\OO_\mu}$
is projective and generically finite.

(b)
$T_x$ meets a   conjugacy class 
$\OO_\sy$ iff
it meets the Lusztig stratum
$\LL_\sy$ through $\sy$\ie iff
$x\in\barr{\LL_\sy}$. If $\OO_\sy=\OO_{\EE,\tmu}$ as in
\ref{Deformation to a nilpotent orbit},
this is equivalent to the relation
$\OO_x\sub\barr{\OO_{\mu_\sy}}$ of two nilpotent orbits.
 Then the variety
$\tT^{a,\sb}_x\dff
T_x \tim_{\barr{\OO_{\EE,\tmu}}}\ \tgab$ is smooth
and connected, its moment map image
is $T_x\cap \barr{\OO_{\EE,\tmu}}$
and it is finite over $T_x\cap\OO_{\EE,\tmu}$.

\begin{proof}
(a)
Smoothness and connectedness follow from lemma
\ref{ontransverseslices}.
The claims about the map $\mmm_a: \tT_{x}^{a}\to T_{x}\cap \barr{\OO_\mu}$
follow from the corresponding claims for the map
$\TN^a\to\barr{\OO_\mu}$.



(b)
If $T_x$ meets $\OO_\sy$ then it certainly meets $\LL_\sy$.
If $T_x$ contains a  point $\sy\in\LL_\sy$
we present its conjugacy class as
$\OO_{ \io(\EE),\io_*\tmu}$
for some injection $\io:\EE\inj\Ao$.
Since the $u$-action is
given by conjugation (which preserves $\OO_{ \io(\EE),\io_*\tmu}$)
and by
multiplication with the scalar $u^2$, we see that
$u\ast \sy$ is in
$\OO_{u^2\cd\io(\EE),(u^2\io)_*\tmu}$.
Therefore, as $u\to 0$ this orbit degenerates to
$\barr{\OO_{\mu_\sy}}$, hence $x=\lim_{u\to 0}\ u\ast\sy$
lies in  $\barr{\OO_{\mu_\sy}}$.

In the opposite direction, if $x\in\barr{\OO_{\mu_\sy}}$
then $T_x$ meets $\OO_{\mu_\sy}$ by lemma
\ref{ontransverseslices}.2.
Since $^GT_x\sub\fg$ is open, the deformation of $\OO_{\mu_\sy}$ through
all $\OO_{c\EE,(c\io)_*\tmu}$ for $c\in\Gm$ implies that
$T_x$ meets $\OO_{c\EE,(c\io)_*\tmu}$ for small $c$,
and then also for $c=1$ by using the $\Gm$-action on $T_x$.
So, $T_x$ meets $\OO_{\EE,\tmu}=\OO_\sy$.
The remaining claims now follow as in similar proofs above.
\end{proof}

\np
\se{
Beilinson-Drinfeld Grassmannians
}
\label{grassmannSection}

In this section $U$ is a vector
space of dimension $m$ and $G$ denotes  the group
$GL(U)$.
We will fix a basis
$\bfe_1,\dots,\bfe_m
$ of $U$, hence
$G=GL_m$.

In \ref{Grassmannians GG Cn}
we recall the loop Grassmannians  $\GG$
of the group $GL(U)$\ie the moduli of $\C[[z]]$-forms of $U((z))$,
and its various versions.
In
\ref{Isomorphism of a pronilpotent slice to z and a neighborhood of sL}
we identify a neighborhood $\UU$ of a point
$\sL\in \GG$, with
the pronilpotent part of the  ``regular''normal slice $T_z$ 
to the  conjugacy class in $\gl(\sL)$ of the operator
$z$ on  $\sL$. 
This restricts to isomorphisms of normal slices in $\GG$ and in matrices
(\ref{A filtration on the slice TT si sub GG0 is given by nilpotent slices in matrices}).
These isomorphisms  lift to resolutions of slices in
\ref{Comparison of resolutions of slices}.
In \ref{Embedding the regular slice into the rational BD-Grassmannian}
we remove the  pronilpotency restriction in $T_z$ by passing to the Beilinson-Drinfeld
deformation of $\GG$.
Finally, the representation theoretic aspects
are
in \ref{Perverse sheaves on loop Grassmannians}.


\sus{
Grassmannians $\GG_\Cn$
}
\lab{Grassmannians GG Cn}
\lab{globalpicture}
Let $C$ be a  smooth curve over $\C$.
The loop Grassmannian construction
for $G=GL(U)$
yields
over each symmetric power
$C\syp{n}$
a reduced ind-scheme
$
\GG_\Cn\to \Cn
$ (this is the reduced part of the space
$GRAS_G$ 
from  4.3.14 of \cite{BD}, which is there called  the ``big Grassmannian'', see also \cite{TZ}). 
The $\C$-points in the fiber at
$b\in C\syp{n}$
are
vector bundles $\LL$ over $C$
which are identified  with $U\ten\OO_C$ off $b$.

Here, ``off $b$'' means ``over
$C-\ub$'',
for the support subset $\ub\sub C$ of $b$
(we will denote $C-\ub$ just by $C-b$).
The condition implies that
the vector bundle
$\LL$ is
a subsheaf of
the quasicoherent sheaf
$j_*j^*(U\ten \OO_C)$ for
$
j
:
C-\ub
\inj
C
$.
We will call pairs $(b,\LL)$ above ``$C$-lattices'' for $U$.
We omit $b$ when it seems obvious. 

The  {\em rational
Grassmannian} $\GG^{rat}_C$ is the moduli of vector bundles
$\LL$ on $C$ equipped with a rational trivialization $\ze:U\ten {\C(C)}\con \LL|_{\C(C)} $
(see  \cite{TZ}). 
This extends to  an identification over $C-b$ for some $b$, so $\GG^{rat}_C$ is 
 a certain colimit 
of spaces $\GG_\Cn$.



For an open cover $C=V_1\cup V_2$ we can glue $V_i$-lattices $(b_i,\AA_i)$, 
with $b_i$ disjoint from $V_1\cap V_2$,
to a $C$-lattice $(b_1\cup b_2,\AA_1\AA_2)$ where $\AA_1\AA_2$ is $\AA_i$
on $V_i$. 
Then
$\Ga(C,\AA_1\AA_2)=\cap_i\  \Ga(V_i,\AA_i)$.
Over any curve $C$ we have the sheaf $U_C\dff U\ten\OO_C$ 
and if $C$ is affine we identify it with $U\ten\OO(C)$.
Then for $C'$ open in $C$ we can extend any $C'$-lattice $(b,\LL)$
to a $C$-lattice $\LL U_{C-b}$.

\sss{
Local loop Grassmannians
}
\lab{Local loop Grassmannians}
\label{Loop Grassmannian and lattices}

For  a finite subscheme $b\in C\syp{n}$
denote by
$\hatt b$ 
the formal neighborhood  of $b$
and by  $\tii b\dff \hatt b-b$ its punctured version.
One can extend the preceding section to local curves such as $\hb$, so the loop Grassmannian $\GG_{\hatt b}$ is defined as  the moduli of
$\hatt b$-lattices (vector bundles 
over $\hatt b$, trivialized to $U_{\tii b}$ over $\tii b$).
Then the fiber
$\GG_b$ of $
\GG_{C\syp{n}}$ at $b\in C\syp{n}$ is identified with
$\GG_{\hatt b}$ 
by restricting vector bundles from $C$ to $\hb$.
The inverse extends a $\hb$-lattice 
$L$ to a $C$-lattice $(b,L U_{C-b})$.

Denote by
$\sO\dff \C[[z]]\sub
\sK\dff \C((z))$  the
formal power series and Laurent series.
Then the group theoretic Loop Grassmannian is
$\GG\dff \GK/\GO$.

A choice of a 
formal parameter
$z$ at a point
$c\in C$ 
identifies a
$\hatt c$-lattice $L$
with  the $\sO$-submodule $\Ga(\hatt c,L)$
 of $U_\sK=U\ten_\C\sK$,
such that the map $L\otimes_{\sO} \sK\to U_\sK$
is an isomorphism.
For instance
to a  coweight $\la\in\Z^m$ of $G$, one attaches
the $\hatt c$-lattice $\sL_\la=\ \pl_1^m\ \C[[z]]\cd z^{-\la_i}\es_i\in\GG_c$.

The group $\GK$ acts on $\GG_\hc$, the stabilizer of $\sL_0=U_\sO$ is $\GO$ and this identifies
$\GG_c=\GG_{\hatt c}$
with  
$\GG\dff \GK/\GO$.
Any coweight $\la$ generates an orbit
$\GG_\la\dff G(\sO)\cd \sL_\la$
and this parameterizes 
$G(\sO)$-orbits
in $ \GG$
by 
dominant coweights   (i.e., $\la_1\ge\cddd\ge \la_m$).

The connected component
$_s\GG$ 
of $\GG$ corresponding to $ s\in \Z$, consists of all $\hc$-lattices $L\sub U_\sK$ such that 
 $\dim(L)\dff \dim(L/L\cap \sL_0)-\dim(L_0/L\cap \sL_0)$ equals $s$.
We will also use $\GG^+\sub\GG=\GG_\hc$ consisting of
all $\hatt c$-lattices $L$ that contain $\sL_0=U\ten_\C\sO$.

\srem
The {\em factorization}
property of loop Grassmannians says that at any $b\in\Cn$ the fiber
factors into a power of the local loop Grassmannian
\begin{equation}\label{bdproduct}
\BDG_b\cong\ \prod_{c\in \ub} \GG_c .
\end{equation}

\sss{
Convolution spaces
$\TGG^a_{\mu}$
}
\label{convolutionlocal}

The $G(\bOO)$-orbits $\GG_\mu$ in
the intersection 
$ _s\GG\cap\GG^+$
correspond to
partitions
$\mu$
of $s$ into at most $m$ parts (there are none for $s<0$).
For such $\mu$ the  orbit $\GG_\mu$ consists of 
$\hc$-lattices $L$  such that $z$ acts on
$L/\sL_0$ as a nilpotent of type $\mu$.

Any choice of a reordering
$a=(a_1,\dots, a_n)$ of
the dual partition $\cmu$
defines a variety of
$n$-step flags of $\hz$-lattices in $U_\sK$:
$$
\TGG^a_\mu
\dff\
\{\sL_0= L_0\sub L_1\sub\cddd\sub L_n\in \GG_{\mu};\
\dim L_i/L_{i-1}=a_i \aand  z(L_i)\sub L_{i-1} \}
$$
Since each  $a_i$ is $\le m$ it corresponds to a
fundamental coweight $\om_{a_i}$ of $GL_m$ and 
$\GG_{\om_p}\cong Gr_p(m)$, Then
$\TGG^a_\mu$
can be described as the
convolution scheme
$\GG_{\om_{a_1}}\ast\cddd\ast \GG_{\om_{a_n}}$\fttt{
In terms of $\varpi:\GK\to G(\sK)/G(\sO)=\GG$ for $A\sub \GG$
let us  denote
$\tii A\dff \varpi\inv A$.
Then the convolution of $G(\sO)$-invariant $A_i\sub\GG$
is
$A_1\ast\cddd\ast A_q\dff\tii A_1\tim_{G(\sO)}\cddd\tim_{G(\sO)}\tii A_q\tim_\GO pt
$.
}.

\srems (0)
The convolution map
$\pi^a_{\mu}:\ \TGG^a_{\mu}\ra\barr{\GG_\mu},\
\pi^a_\mu(L_0,...,L_n)=L_n$, is
a partial resolution of singularities \cite{MV1}.

(1)
The fiber $(\pi^a_{\mu})^{-1}(L)$
at a point $L$ in
$\barr{\GG_\mu}$
is by definition isomorphic to  the Springer-Ginzburg fiber
$
\FF^a_x
$
defined in
\ref{ginzburgfiber}
for the operator $x$ on $D=L/\sL_0$, given by the $z$-action on
$L/\sL_0$.

\sus{
Map $\Phi$ from BD-Grassmannians $\GG_\An$ to linear operators
}

\lab{From loop Grassmannians to linear operators}
\label{ normal slices in loop Grassmannians and nilpotent cones}

This transition
is global, it uses $\P
\dff\Po$ with coordinates $z^\pmo$ over $\A=\P-\yy$ and $\Am=\P-0$.
\lab{Normal slices}

In this section we are always in a setting of one of two  ``formally open covers'' $(V_0,V_\yy)$ of $\P$,
either $(\hz,\Am)$ or $(\A,\hy)$. In each case we will study an open part
of the moduli of $V_0$-lattices defined by ``complementarity'' with a fixed $V_\yy$-lattice.
The results in two cases are parallel.

\sss{
Maps 
\ud{$\Phi$}$^{\II,K}:
\UU^{I,\KK}_\An
\to 
\gl(H^1(\P,I\KK))
$
on
spaces  of complementary lattices 
}

\lab{Maps Phi b}

A  choice of an $\hy$-lattice $I\sub U_\ty= U((z\inv))$
defines an open submoduli $\UU_\An^I\sub \GG_\An$
of $\A$-lattices
$\LL$
that are  {\em complementary} to $I$  in the sense that
the composition of a  restriction and a quotient
$\Ga(\A,\LL)\to\Ga(\ty,I)=U_\ty\sur U_{\tii\yy}/I$ is an isomorphism\ie 
$H^1(\P,\LL I)=0$.

Also, for any $\A$-lattice $\KK\sub U_\A$ let $\GG_\An^\KK$ be the $\A$-lattices
$\LL\in\GG_\An$ that contain $\KK$.
Denote $\UU^{I,\KK}_\An\dff\ \UU^I_\An\cap \GG_\An^\KK$.
We also glue $I$ and $\KK$ over $\ty$ to get a $\P$-lattice 
$
I\KK$.

\slem
For any $\A$-lattice
$\LL$ in $\UU^{I,\KK}_\An$ we have
$
\Ga(\A,\LL/\KK)
\con
H^1(\P,I\KK)
$.

\pf
The isomorphism $\Ga(\A,\LL)\con U_{\tii\yy}/I$ 
induces $\Ga(\A,\LL/\KK)\con U_{\tii\yy}/[I+\Ga(\A,\KK)]$
which is $H^1(\P,I\KK)$ since $I\KK$ is equal to:\
$U_\ty$ on $\ty$, \  $\KK$ on $\A$, and\  $I$ on $\hatt\yy$.
\epf

Now we can   define the map that sends an $\A$-lattice $\LL$ to the $z$ action on $\Ga(\LL/\KK)\cong H^1(\P,I\KK)$
$$
\uPhi^{I,\KK}:
\UU^{I,\KK}_\An\to \gl(H^1(\P,I\KK))
.$$
As $\AA$-lattices $\KK$ become small
these maps glue to 
$\uPhi^{I}:
\UU^{I}_\An\to \gl({U_\ty/I})$.



\srems (0) Let $\UU^{rat,I}_\A$ be the image of all $\UU^I_\An$
in the rational Grassmannian $\GG^{rat}_\A$.
The map $\uPhi^I$ factors through $\uPhi^{rat,I}:
\UU^{rat,I}_\A\to\gl(\uLL)$ since it does not use the 
estimate $b$ on singularity of trivialization.

(1) There are also the analogous maps
$\uPhi^{\II,K}:\UU^{\II,K}_\hz\to H^1(\P,\II K)$
where 
an   $\Am$-lattice $\II$ and a $\hz$-lattice $K$
define subspaces of the ordinary loop Grassmannian 
$\GG_\hz\subb \UU^\II_\hz\subb\UU^{\II,K}_\hz$.


\sss{
A reference point 
\underline{$\LL$} in 
$\UU^I_\An$
} 

\lab{A reference point uLL in UU I An}

For a  choice of $\uLL\in \UU^I_\An$
we can now move $\uPhi^I:\UU^I_\An\to \gl({U_\ty/I})$
to $\vPhi^I:\UU^I_\An\to \gl(\uLL)$.\fttt{
Similarly, an   $\Am$-lattice $\II$ and a complementary $\hz$-lattice $\sL$ 
define 
$\vPhi^{\II}:\UU^{\II}_\hz\to \gl(\sL)$.
}

For $\LL\in\UU^I_\An$ denote the complementarity isomorphism  by $\io_\LL:\LL\con  U_\ty/I$.
Then $\vPhi^I(\LL)=\ ^\io(z|_\LL)$,
for $\io=\io_\uLL\inv \io_\LL:\LL\con\uLL$.

Furthermore, any $\LL\in \UU^I_\An$
is the graph $\LL=(1+f)\uLL\sub U_\ty$ of a unique linear map $f:\uLL\to I$.
Then $\io=1+f$ since $(\uLL\con \LL\sub U_\ty)$ equals  $(1+f)\ci(\uLL\sub U_\ty)$. 
So,  
$\vPhi^I(\LL)=\ ^{1+f}(z|_\LL)$.

According to the remark \ref{Maps Phi b}.0, maps 
$\vPhi^I:\UU^I_\An\to \gl(\uLL)$ factor to 
$\vPhi^I:\UUratI_\An\to \gl(\uLL)$.

\sus{
Map $\Psi_0$ from nilpotent operators to
the loop Grassmannian $\GG_\hz$ 
}
\lab{Isomorphism Psi}
\label{From operators to lattices}

As above, we consider an $\Am$-lattice
$\II$ and a $\hz$-lattice $\sL$ in the 
open submoduli  $\UU^\II_\hz\sub\GG_\hz$ of $\hz$-lattices complementary to $\II$.
We will actually start with a $\C$-form $\sU$ of $\sL$ that also defines $\II$,
since $\sU$  gives a  ``regular'' normal slice $T_z=z+\Hom(\sL,\sU)$ to the conjugacy class of $z$ in 
$\gl(\sL)$. For $g=z+\vphi\in T_z$, we will find two incarnations $f$ and $\tvphi$ 
of $\vphi$ which will give two formulas $\Psiz(g)=(1+f)\sL=(1-z\inv\tvphi)\inv\sL$ 
for the $\hz$-lattice associated to $g$.

\sss{
$\C$-forms and neighborhoods of lattices
}
\lab{A C-form of a lattice in GG_0 defines its neighborhood}

Any $\C$-form
$\sU$ of
$U_\Gm=U[z^\pmo]$, gives
an $\A$-lattice $\uLL\dff\sU[z]$ and an  $\Am$-lattice $\II\dff z\inv\sU[z\inv]$ so that 
over $\Gm$ both are identified with $\sU[z^\pmo]=U_\Gm$.
We restrict these to $\hz$ and $\hy$
to get $\sL=\sU[[z]]\sub \sU((z))=U((z))$ and
$I=z\inv\sU[[z\inv]]\sub \sU((z\inv))=U((z\inv))$.

We have defined $\UU^I_\An$ as all $\LL\in \GG_\An$ that satisfy the condition:\ 
(i) $\LL\pl I=U(((z\inv))$; and
$\UU^\II_\hz$ as all  $L\in \GG_\hz$ such that: \  
(ii) $L\pl\II=U((z^\pmo))$.
Then
$
\UU^{I}_\An \sub \GG_\An
$  and
$
\UU^{\II}_\An \sub \GG_\An
$  are neighborhoods of
$
\uLL$ and $\sL$. Say, 
$\uLL =\sU[z]$  
is complementary to
$I=z\inv\sU[[z\inv]]$ in $\sU((z\inv))=U((z\inv))$.

\srems (a)
The extension operation $L\mm \LL\dff LU_\Gm$ 
(from \ref{Grassmannians GG Cn}),
identifies  $\GG_\hz$ with 
the moduli $ H^1(\A,\Gm,GL(U))$ of 
vector bundles $\LL$ over 
$\A$ trivialized to $U_\Gm$ over $\Gm$. 
An example is $\sL\mm \uLL$. Similarly for  $I\subb\II$ which are
lattices over $\hy\sub\Am$.

(b) For  $L$ and $\LL$ as in (a),
the following complementarity claims are equivalent:
(i)$\eq$(ii)$\eq$(iii),\ where the last one is 
 $L\pl I=U[[z^\pmo]]\dff\prod_\Z \C z^p$.
For instance  (iii)$\imp$(ii)
by intersecting with $U((z))$ and the converse is by completion.



\sss{
Isomorphism $\Psi_0$ of  a
pronilpotent normal
 slice at  $z\in \gl(\sL)$ 
and a neighborhood of $\sL
\in\GG_0$ 
}
\lab{Isomorphism of a pronilpotent slice to z and a neighborhood of sL}

As in \ref{A C-form of a lattice in GG_0 defines its neighborhood} we use the notation
$\sU,\sL,\II,\UU^{\II}_\hz$.
Let  $\Hom(\sL,\sL)=\gl(\sL)$ denote the continuous linear operators.
Say,  its subspace $\Hom(\sL,\sU)$ consists of operators
$\sL\to \sU$ that vanish on some sublattice
$K\sub \sL$.

Also, our $\C$-form $\sU$ of $\sL$ gives  $T_z\dff z+\Hom(\sL,\sU)$ which is the ``regular'' normal slice
in $\gl(\sL)$ to the pronilpotent operator
$z$ (see \ref{pronilpotent}).
Let $\NN_z\sub\gl(\sL)$ consist of all operators
$A$ for which there is a  sublattice $K\sub \sL$
such  that $A=z$ on $K$ and $A$ is nilpotent on $\sL/K$.

\spro
(a)
The following map
is an isomorphism
of the pronilpotent part of the normal slice
to the coadjoint orbit at $z\in \gl(\sL)$
to a neighborhood of $\sL$ in $\GG_\hz$:
\fttt{
Here,
$\vphi\in \Hom(\sL,\sU)\sub\gl(\sU)$ and we use the composition
$
\sL
\ \aa{(z+\vphi)^{k-1}}
\lra
\
\sL
\
\aa{\vphi}
\ra
\
\sU
\
\aa{z^{-k}}
\ra
\
z^{-k}\sU
$.
}
$$
\Psi_0:\
T_z\cap\NN_z
\ \con\ \UU^{\II}_\hz
,\ \ \ \
\text{by}
\ \ \
\Psiz(z+\vphi)\
\dff\
[1+\sum_{k=1}^{\infty}z^{-k}\vphi(z+\vphi)^{k-1}]\sL
.$$

(b) 
Let $g=z+\vphi$ and  $f=\sum_{k=1}^{\infty}z^{-k}\vphi(z+\vphi)^{k-1}$, then 
$1+f: \sL\con \Psiz(g)$ and this identifies the action of $z$ on $\Psiz(g)$
and $g$ on $\sL$. Moreover, $\vPhi^I\ci \Psi_0=id$.


(c) $\Psi_0(z+\vphi)\cap \sL$ is the largest $\hz$-lattice inside $\Ker(\vphi)$.


\begin{proof}
(a$_1$)
Since we know one lattice
$\sL$ in $\UU^\II_\hz$,
all lattices $L\in \UU^\II_\hz$
are graphs   $(1+f)\sL$
of continuous $\C$-linear maps $f: \sL\to \II$.
Here,
$
\II
=
\pl_{k>0 }\ z^{-k}\sU$
decomposes
$
f$ into $\sum_{k>0}z^{-k}f_k
$
with  $f_k\in\Hom(\sL, \sU)$.
(The sum is finite since $f$ vanishes on some sublattice of $\sL$.)

For which $f$'s is $(1+f)\sL$ a  $\hz$-lattice? Invariance of
$(1+f)\sL$ under $z$ means that
for each $v\in \sL=\prod_{i\ge 0} z^{i}\sU $ one has
 $z(1+f)v
=(1+f)w$ for some $w\in \sL$\ie
$zv+\sum_{l\ge 0}z^{-l}f_{l+1}v=
w+\sum_{k>0}z^{-k}f_kw$.
The summand in $\sL$ is $zv+f_1v=w$
and for $p>0$ the summand in $z^{-p}\sU $
is
$z^{-p}f_{p+1}v=z^{-p}f_pw$.
Then $g=z+f_1$ lies in $\gl(\sL)$ and
$z$-invariance means that:\ (i)
the only choice for $w$ is $gv$ and\
(ii) $f_1$ freely determines $f$ by $f_p=f_1g^{p-1}$
for $p>1$.

So, for $f$'s such that $(1+f)\sL$ is $z$-invariant 
we have found that 
$z(1+f)=(1+f)g$ on $\sL$\ie $^{(1+f)\inv}(z|_\sL)$ equals
$g=z+f_1$. 


(b) We just checked the first claim. The second follows since 
$\vPhi^I(\LL)=\ ^{1+f}(z|_\LL)$ by \ref{A reference point uLL in UU I An}.

(a$_2$)
It remains  to see that the $z$-invariant space $(1+f)\sL=(1+\sum_{k>0}z^{-k}f_k)\sL$
is a $\hz$-lattice. This is equivalent to asking that
$f_k=f_1(z+f_1)^{k-1}$ is zero  on $\sL$ for some $k$.

Since $f_1$ is continuous it vanishes on sufficiently small sublattices $K\sub\sL$.
Then on $K$ we have
$g=z$ and $(1+f)K=K$. So, the condition
on $f_k$ reduces to vanishing  on $W$ for some complementary subspace $K\pl W=\sL$.
If $1+f$ gives a $\hz$-lattice then
$g$
is nilpotent on
$\sL/K$
since it behaves the same as
$z$ on $(1+f)\sL/K$.
Conversely if $g$ is nilpotent on $\sL/K$
then
$f_1g^p$ vanishes on $W\cong\sL/K$ for large $p$.



(c) By definition,
$\Psiz(z+\vphi)$ is $(1+f)\sL$ with
$f_1=\vphi$ and $f=\sum_1^\yy\ z^{-k}\vphi(z+\vphi)^{k-1}$.
Since $f\sL\sub I$ we have $(1+f)\sL\cap\sL=\Ker(f)$.
Then $\Ker(f)$ is  the intersection of kernels of
maps
$
z^{-k}
\vphi
(z+\vphi)^{k-1}:
\sL\to z^{-k}\sU
$.
Case $k=1$ gives  $\Ker(f)\sub \Ker(\vphi)$. Also, on any
sublattice $K$ of
$ \Ker(\vphi)$ we have
$(z+\vphi)K=zK\sub K$ hence
$\vphi(z+\vphi)^{k-1}K=0$.
\end{proof}

\srem
The isomorphism in the proposition
translates a moving lattice
$L\in\UU^\II_\hz$ with the ``constant'' operator $z$,
into a moving operator
$g$ on a fixed  vector space $\sL$.


\sss{
 Example: neighborhood $\UU^\si_\hz\sub\GG_\hz$ of $\sL_\si
$}

\lab{Example: neighborhood UUsi of sLsi}
\lab{UU si}


We will now specialize  
 \ref{A C-form of a lattice in GG_0 defines its neighborhood} to the case  
when the 
$\C$-form $\sU$
comes from a coweight $\si\in\Z^m$ as  
 $\sU_\si\dff\pl_i\ z^{-\si_i}\C\bfe_i$.
This gives $\II_\si\dff z\inv\sU_\si[z\inv]$
and  the neighborhood
$\UU^{\si}_\hz\dff \UU^{\II_\si}_\hz$ of $\sL_\si=\sU_\si[[z]]$.

\sss{
Slices $\TT_\si\sub\GG_\hz$  as ``regular'' pronilpotent slices in matrices
}

\lab{A filtration on the slice TT si sub GG0 is given by nilpotent slices in matrices}
\lab{The nilpotent part of the slice}

The negative congruence subgroup $L^{<0}G$
is the kernel of the evaluation $G[z\inv]\to G$ at $z=\yy$.
The standard normal slice at $\sL_\la$ to $\GG_\la=\GO\sL_\la$ is the orbit
 $\TT_\la\dff L^{<0}G\cd\sL_\la$.

For a coweight $\si$ and a choice of  
$S\in\Z $ such that $\sL_\si\subb {z^S\sL_0}$,
define the data $\sN,N,\bz,
\la$ so that 
 $\bz$ is the action of $z$ on
the vector space $\sN\dff \sL_\si/z^S\sL_0$ of dimension $\sn$. Also, when we  reorder $\si-S
\dff (\si_1-S,...,\si_m-S)$ into a partition $\la$ 
of $\sn$
then $\bz$ lies in the nilpotent orbit $\OO_\la$
 in the nilpotent cone $\NN\sub\gl(\sN)$.

We will now restrict the description of the neighborhood $\UU^\si_\hz$ from 
\ref{Isomorphism of a pronilpotent slice to z and a neighborhood of sL}
to the slice $\TT_\si$, to $\GG_\si=\GO\sL_\si$ by moving from $z$ on  $\sL_\si/z^S\sL_\si$
to $z$ on  $\sL_\si/z^S\sL_0$\ie a from a    box to a general Young diagram.

Let $T_\bz=T_\la$ be
the ``regular'' normal
slice (at $\bz$) to $\OO_\la\sub\gl(\sN)$ 
from
\ref{regular slice}
(where  we choose the $z$-generating subspace $\CC\sub\sN$ as $\sU_\si$).
We will embed
$T_\bz$ 
into the slice 
 $ T_z$ from
\ref{Isomorphism of a pronilpotent slice to z and a neighborhood of sL}
as  operators  on $\sL$
which on $z^S \sL_ 0$
agree with the multiplication by $ z$.
Then 
its nilpotent part 
$T_ \bz \cap\NN$  
lies in $T_z\cap\NN_z$
(\ref{Isomorphism of a pronilpotent slice to z and a neighborhood of sL}).



\stheo
(a)
The isomorphism $ \Psiz:
T_z\cap\NN_z
\con \UU^\si_\hz$ from
\ref{Isomorphism of a pronilpotent slice to z and a neighborhood of sL},
restricts to an isomorphism of the
nilpotent part of  the ``regular'' normal
slice $T_\bz$
and the $S\thh$ piece 
of the slice $\TT_\si\sub \GG_\hz$:
$$
\Psiz^S: T_\bz\cap\NN
\con\
\TT_\si\cap \GG_\hz^{z^S\sL_0}
.$$


(b)
If
$\OO_\la\sub\barr{\OO_\mu}$
for a nilpotent orbit $\OO_\mu$ in $\gl(\sN)$
then $\Psiz^S$ restricts to
an isomorphism of intersections
$\Psi^S_{0,\mu}:
T_\bz\cap
\barr{\OO_\mu}
\con
\TT_\si\cap \barr{\GG_{\mu+S}}$.

\srem This identifies 
the slice $\TT_\si=L^{<0}G\cd\sL_\si$ with the ``regular'' normal slice to the orbit of
the pronilpotent operator $z$ on $\sL_\si$
(\ref{pronilpotent}).

\bep
(a) According to  the proposition
\ref{Isomorphism of a pronilpotent slice to z and a neighborhood of sL}.c,
for $g=z+\vphi$, $\Psiz(g)$ contains $z^S\sL_0$ iff $\vphi=0$ on
$z^S\sL_0$\ie
$\vphi:\sN\to \sU_\si$ for $\sN=\sL_\si/z^S\sL_0$.
Then we can consider $g$ as an operator 
$\sN\to\sU$.


Observe that
$\TT_\si
=L^{<0}G\cd\sL_\si$ is a torsor for
$L^{<0}G\cap\ ^\si (L^{<0}G)$, so it  lies in the neighborhood
$^\si (L^{<0}G)\cd\sL_\si$ of $\sL_\si$ which is exactly $\UU^\si_\hz$.
It remains to show that
for
 $\vphi\in \Hom(\sN,\sU_\si)$ such that
$g=\bz+\vphi$ is nilpotent,
the corresponding $\hz$-lattice
$L = \Psiz(g)$ in $\UU^\si_\hz\cap \GG_\hz^{z^S\sL_0}$
lies in $ \TT_\si$
iff $\vphi$ lies  in the subspace
$C\dff T_\bz-\bz$ of $\Hom(\sN,\sU_\si)$.

Group $\RR$ of loop rotations
is the group $G_m$ acting on
$U((z))$ by $s\bu z^k\bfe_i\dff (sz)^k\bfe_i$, hence also  on
the space $\GG_\hz$ of lattices.
One knows
(\cite{MV1}) that a $\hz$-lattice $L\in \GG$ is in  $L^{<0}G\cd
\sL_\si\dff \TT_\si$
 if and only if\
$\lim_{s\to\infty} s\bu L=\sL_\si$.

Let us
write $\Psiz(g)$ as $(1+f)\sL_\si$ for $f=\sum_1^\yy\ z^{-k}\vphi(z+\vphi)$
as in the proof of proposition
\ref{Isomorphism of a pronilpotent slice to z and a neighborhood of sL}.
Since $\sL_\si$ is fixed by  $\RR$,
we have
$
s\bu (1+f)\sL_\si
=(1+s\bu f) \sL_\si
$.
Therefore,
$\lim_{s\to\infty} s\bu \Psiz(g)=\sL_\si$ is equivalent to
$\lim_{s\to\infty}s\bu f=0$.
Since $f=\sum_1^\yy\ z^{-k}\vphi(z+\vphi)$,
$$
s\bu f=\sum_{k=1}^{\infty}(sz)^{-k} (s\bu \vphi )(sz + (s\bu \vphi ))^{k-1}
=\   \sum_{k=1}^{\infty} \
z^{-k}\ s\inv (s\bu \vphi )\ (z+s^{-1}(s\bu \vphi ))^{k-1}
.$$
So,
$\lim_{s\to\infty}s\bu f=0$ is equivalent to
$\lim_{s\to\infty}s^{-1}(s\bu \vphi )=0$.

Here,
$\vphi :\sL_\si\to \sU_\si$
has the $(i,j)$-component maps
in
$
\Hom(
z^{-\si_j}\OO ,z^{-\si_i}\C )
\cong\
\pl_{p\ge 0}
z^{\si_i-\si_j-p}\C$.
On the  $(i,j,p)$ component $\bu$-action of $s$ is by $s^{\si_i-\si_j-p}$ so
this component has to vanish 
unless
$\si_i-\si_j-p\le 0$.
The remaining components are for 
$p<\si_i-\si_j=
\la_j-\la_i$\ie
$\la_i<\la_j-p$.
However, this condition on $\vphi$ is  the description of
$C\sub\Hom (\sN,\sU_\si)$ from \ref{nilpnormalslice}.


(b) follows  from (a) since
$\Psiz(T_x\cap\OO_\mu)
\sub
\TT_\si\cap \GG_{\mu+S}$
(according to the proposition
\ref{Isomorphism of a pronilpotent slice to z and a neighborhood of sL}.b
the nilpotent operators $g$ on $\sL_\si/z^S\sL_0$
and $z$ on $\Psiz(g)/z^S\sL_0$ have the same  type).
\enp

\sss{
Transform $\vphi \mm \tvphi$ and the congruence subsemigroup $C_\yy$
}
\lab{Congruence subsemigroup}
We return to the setting of $\sU,\sL,\II$ from 
\ref{Isomorphism of a pronilpotent slice to z and a neighborhood of sL}.
If we write  $v\in \sL$  as  $\sum_0^\yy\ z^av_a$ with $v_a\in \sU$,
then
any  continuous linear operator $\vphi:\sL\to\sU$
can be written  as
$
\vphi(\sum_0^\yy\ z^av_a)=\sum_0^\yy z^{-a}\vphi_av_a$
where  $\vphi_a\in \End(\sU)$ is $\sU \aa{z^a}\to \sL\aa{\vphi}\to \sU$.
Then 
$$\tvphi\dff\ \sum_{i\ge 0} z^{-i}\vphi_i\ \in
\End(\sU)[z\inv]
.
$$

We will call 
$C_\yy\dff 1+z\inv\End(\sU)[z\inv]\sub\ G[z\inv]$ the
{\em congruence subsemigroup} (for $\sU$).

\spro
There is a canonical bijection $T_z\con C_\yy$ by $z+\vphi\mm\ 1-z\inv\tvphi$.

\pf
We have already constructed this map.
Its inverse ``restricts'' $\tvphi$ to $
\vphi\dff \ (\sL\sub \sU_\sK\aa{\tvphi}\to\sU_\sK\aa{pr_{\sU}}\to 
\sU)$.
\epf

\sss{
A rational formula for  $\Psi_0$
}

\lab{A rational formula for isomorphism Psi0}
Here we  use the transform $\vphi \mm \tvphi$ from \ref{Congruence subsemigroup}.

\spro
For $g=z+\vphi\in T_z\cap\NN_z$
we have
$$
\Psi_0(g)
\ =\  (1-z\inv\tii\vphi)\inv \sL
.$$

\pf
We have $\Psi_0(g)=(1+f)\sL$ for $f=\sum_{l\ge 0}z^{-l-1}\vphi(z+\vphi)^l$, and we need 
$(1+f)\sL=(1-z\inv\tii\vphi)\inv \sL$.
We will notice that it suffices to  check that
$1+ f$ equals $(1+z\inv\tii\vphi)\inv$
on $\sU\sub \sL$.
The reason is that
as $\C[[z]]$-modules,
$(1-z\inv\tii\vphi)\inv  \sL$
and
$(1+f)\sL$
are generated by
$(1-z\inv\tii\vphi)\inv \sU$
and
$(1+f)\sU$.
The first claim
is obvious.
By proposition \ref{Isomorphism of a pronilpotent slice to z and a neighborhood of sL}.b
the second claim is really that
$\sU$ generates $\sL$ as a $\C[[g]]$-module.
This in turn holds because
on $\sL$ we have $g=z$ modulo $\sU$.

Now we will calculate
$
1+f
=\ 1+\sum_{l\ge 0}z^{-l-1}\vphi(z+\vphi)^{l}
$
on $\sU$.
Notice that
each $\vphi(z+\vphi)^{l}$ is a sum of products
$\vphi z^{a_1}\cddd \vphi z^{a_s}$  with $a_1+\cddd+a_s=(l+1)-s$.
Each factor $\vphi z^a$ preserves $\sU$
and acts on it as the linear operator $\vphi_a$.

So,
on $\sU$ operator $z^{-l-1}\vphi(z+\vphi)^{l}$
is the sum 
of
$
z^{-s-\sum_1^s\ a_j}
\prod_1^s \vphi_{a_j}
$ over
all $a_*$ as above.
Since $\vphi_a$ is $\sK$-linear
this can be rewritten as
$
\prod_1^s z^{-1-a_j}\vphi_{a_j}
$.
So, on $\sU$,
$$
1+f=\ \sum_{s\ge 0}\ \sum_{a_\bu\in\N^s} \prod_1^s\ z\inv(z^{-a_j}\vphi_{a_j})
\ =\ \sum_0^\yy \ (z\inv\tii\vphi)^s=(1-z\inv\tii\vphi)\inv
.$$
\epf

\sus{ 
Embedding
the regular slice $T_z$  into
the rational BD-Grassmannian 
}

\lab{Embedding the regular slice into the rational BD-Grassmannian}
\lab{Map Phi from linear operators to the BD-Grassmannian GGAn }
\label{Isomorphism of slices}
\lab{Isomorphism of slices in loop Grassmannians and in matrices}
\lab{globalpsi}
As in
\ref{Isomorphism of a pronilpotent slice to z and a neighborhood of sL},
for a $\C$-form
$\sU$ of $U[z^\pmo]$ we consider  $\sL=\sU[[z]]\subb \uLL=\sU[z]$ with topology where
the $\hz$ and $\A$ sublattices form  bases of neighborhoods of  $0$.
Then the  ``regular'' normal slice $T_z=z+\Hom(\sL,\sU)\sub \gl(\sL)$ 
to the conjugacy class
of $z$ in $
\gl(\sL)$, embeds into  its $\A$-version $z+\Hom(\uLL,\sU)\sub \gl(\uLL)$.
This will allow us to extend the embedding of the pronilpotent part of  $T_z$
into $\GG_\hz$
(\ref{Example: neighborhood UUsi of sLsi}), to an embedding
of all of $T_z$ into the 
 Beilinson-Drinfeld deformation of $\GG_\hz$. Actually, our formulation will
instead of the system of spaces $\GG_\An$ 
use the  {\em rational Grassmannian}
$\GGr_\A$ from
\ref{Grassmannians GG Cn} (then  
one does not have to choose a bound $b\in\An$
on singularity of the  trivialization). 


\sss{}
\lab{Slice Tz and UU ratI AGm}
For   $I=z\inv \sU[[z\inv]]\in \GG_\hy$ let  $\Phi^{rat,I}:
\UU^{rat,I}_\A\to\gl(\uLL)$ be the map from
\ref{A reference point uLL in UU I An}. 
We will use the transform $\vphi\mm\tvphi$ from \ref{Congruence subsemigroup}.
Finally, let $\cUUratI$ consist of all $\A$-lattices $\LL\in\UUratI_\A$ such that $\LL|_\Gm\subb U_\Gm$.

\stheo
(a) 
The ``regular'' normal slice $T_z$ at the operator $z$ in $\gl(\sL)$
is isomorphic to the moduli of $\A$-lattices 
$
\LL\in \GGr_\A$
that are complementary to $I$ and such that $\LL|_\Gm\subb U_\Gm$.
This moduli is also a (free)  orbit of the 
congruence subsemigroup $C_\yy$ through $\uLL$ in $\GGr_\A$. The isomorphism is
$$
\Psi: T_z\to \cUUratI,\ \ \
\Psi(z+\vphi)\dff\ (1-z\inv\tii\vphi)\inv\uLL
=\ 
(1+\sum_{k=1}^{\infty}z^{-k}\vphi(z+\vphi)^{k-1})
\uLL
.$$

(b) For $g=z+\vphi\in T_z$ 
the actions of $z$ on $\Psi(g)$
and $g$ on $\uLL$
are the same under 
the canonical identification
$\uLL\con \Psi(g)$. 

(c) The inverse of $\Psi$  is a  restriction 
of the map $\Phi^{rat,I}:
\UUratI\to\gl(\uLL)$\ie
$\Phi^{rat,I}\ci \Psi=id_{T_z}$.

\pf
(a$_1$) {\em 
$(1-z\inv\tii\vphi)\inv\uLL$ is an $\A$-lattice.}
Since $A=1-z\inv\tvphi\in\End(U)[z\inv]$ has value $1$ at $z=\yy$,
it is an invertible rational endomorphism of the vector bundle $U_\Po$, so
$A^\pmo$ act on $\GGrat_\A$.

(a$_2$)
 {\em The two formulas for $\Psi$ coincide.}
In \ref{A rational formula for isomorphism Psi0}
we checked that $A\inv$ and $ 
(1+\sum_{k=1}^{\infty}z^{-k}\vphi(z+\vphi)^{k-1}
$ have the same effect  on $
\sL=\prod_{a\in\N} z^a\sU
$
which we can think of as a completion of 
$\uLL=\pl_{a\in\N} z^a\sU$. The claim for $\uLL=gr(\sL)$ follows.\fttt{
Strictly speaking in the setting of \ref{A rational formula for isomorphism Psi0}
the sum in $k$'s is finite but this is not used in the proof.}

(b)
In
proposition
\ref{Isomorphism of a pronilpotent slice to z and a neighborhood of sL}.b
it has been checked  that 
for $g=z+\vphi\in T_z\cap\NN_z$ and  $f=\sum_{k=1}^{\infty}z^{-k}\vphi(z+\vphi)^{k-1}$, 
the isomorphism
$1+f: \sL\con \Psiz(g)$  identifies the action of $z$ on $\Psiz(g)$
and $g$ on $\sL$. The proof works the same for  $g\in T_z$ and  $1+f:\uLL\con\Psi(g)$.

(a$_3$) {\em $\Psi$ maps $T_z$ to $\cUUratI$.}
$\Psi(g)=(1+f)\uLL$ is complementary to $I$ since $\uLL$ is.
Also, over $\Gm$ the operator 
$A=1-z\inv\tvphi\in\End(U)[z\inv]$ is regular hence
$A(\uLL|_\Gm)=A U_\Gm\sub U_\Gm$, hence $\Psi(g)|_\Gm=A\inv U_\Gm \subb U_\Gm$.

(a$_4$)
{\em $\Psi$ is surjective.} For $\LL\in\UUratI_\A$,  complementarity with $I$ shows that 
$\LL=(1+F)\uLL$ for a unique linear map $F:\uLL\to I$.
Since $(1+F)\uLL=\LL$  is $z$-invariant,
we have $f=\sum_{k>0} z^{-k}P(z+P)^{k-1}$ for unique $P:\uLL\to \sU$.
[This was noticed in the proof of 
\ref{Isomorphism of a pronilpotent slice to z and a neighborhood of sL}.a, except that 
there one considers maps from $\sL$ to $\II$ rather than from $\uLL$ to $I$, but the proof is the same.]
\
If $\LL|_\Gm$ also contains $U_\Gm$ (which also lies in $\uLL$),
then $P$ kills $z^S\uLL$ for $S>>0$. Therefore, $P:\uLL \to \sU$ extends to  some $\vphi:\sL\to \sU$.
Then $\LL=\Psi(z+\vphi)$.

(c) follows from (b) and the formula 
  $\Phi^I(\LL)=\ ^{1+f}(z|_\LL)$ from \ref{A reference point uLL in UU I An}.

(a$_5$) Now we know that $\Psi:T_z\con\cUUratI$ is an isomorphism.
Then the  isomorphism $T_z\con C_\yy$ from proposition
\ref{Congruence subsemigroup} interprets $\cUUratI$ as a free right orbit of $C_\yy$.
\epf

\srems
(1)
$\Psi$ extends the isomorphism $\Psi_0:T_z\cap\NN_z
\con \UU^I_\hz$ from
\ref{Isomorphism of a pronilpotent slice to z and a neighborhood of sL},
under the correspondence from the remark 
\ref{A C-form of a lattice in GG_0 defines its neighborhood}.b.
(For this compare the second formula for 
$\Psi$ and the proposition \ref{A rational formula for isomorphism Psi0}.)

(2) Claim (b) allows one to compare interesting subschemes of  $T_z$ and  $\cUUratI$ as in
\ref{Isomorphism of a pronilpotent slice to z and a neighborhood of sL}.a,
and
\ref{A filtration on the slice TT si sub GG0 is given by nilpotent slices in matrices}.

(3)
 For $g\in T_z$, the support of [$\Psi(g)|_\Gm]/ U_\Gm$ is the  spectrum of $g$ with $0$ removed.
(This follows from 
(b) since for $S>>0$,
 $z$ acts on $\Psi(g)/z^S\uLL$ the same as $g$ on $\uLL/z^S\uLL$
and  the restriction to $\Gm$  removes $0$.)


\sex
The map $\Psi$ can be viewed as a (weak) extension from line bundles to vector bundles of the
(normalized)
 Abel-Jacobi map.
A special case is the
embedding of $\gl(m)$ into $\GG_{\A\syp{m}}$
by $x\mm\ \fra{z}{z-x}
\cd \ox^{\pl m}$
that lies above
$gl_m//GL_m$.
Its restriction to
the nilpotent cone $\NN_m$ in $\End(U)$
recovers Lusztig's embedding of $\NN_m$ into the local loop Grassmannian $\GG$ at $0\in\A$
\cite{Lu}, which is  the origin of this paper.


\sus{
Comparison of resolutions of slices
}
\lab{Comparison of resolutions of slices}
\label{normal slices in loop Grassmannians and nilpotent cones}
Let $x=\bz\in\OO_\la\sub\barr{\OO_\mu}$
and let $a$ be any reordering of terms in the dual partition
$\cmu$.
Recall the
partial resolution
 $\pi^a_{\mu}:  \TGG^a_{\mu}\ra\barr{\GG_\mu}$
from the theorem \ref{convolutionlocal}.c
and
$
\TN^{a}
$
from
\ref{tildeta}.
We can lift the map $\psi_\mu\dff \Psi^N_{0,\mu}:
T_{\bz}\cap\barr{\OO_\mu}
\con
\TT_\si\cap\barr{\GG_\mu}
$ from 
\ref{A filtration on the slice TT si sub GG0 is given by nilpotent slices in matrices},
to the map
$$
\tpsi_\mu: \tT_{\bz}^{a}\dff T_\bz
\
\tim_{\NN}
\
\TN^{a}
\ \ra  \
(\TT_\si\cap\barr{\GG_\mu})\ \tim_{\barr{\GG_\mu}}\ \TGG^a_{\mu}
,$$
since any $(\bz+\vphi)$-invariant $n$-step flag in $\sN$ will give rise
to an $n$-step flag of $\hz$-lattices in $L=\Psi_0(\bz+\vphi)$.

\scor
Maps $(\tpsi_\mu,\psi_\mu)$ form an isomorphism
of resolutions
\begin{equation}
\begin{CD}
\tT_{\bz}^{a}
@>{\tpsi_\mu}>\cong>
\TGG^a|_{\TT_\si\cap\barr{\GG_\mu}}
 \\
@V{\mmm_a}VV
@V{\pi^a_\mu}VV
\\
T_{\bz}\cap\barr{\OO_\mu}
@>{\psi_\mu}>\cong>
\TT_\si\cap\barr{\GG_\mu}
.\end{CD}
\end{equation}

\sus{
Perverse sheaves on loop Grassmannians
}
\lab{Perverse sheaves on loop Grassmannians}
\label{satake}
In this subsection $G$ denotes $GL_m$ or $PGL_m$.

\let\Perv\PP

Let $\Perv_{G(\sO)}(\GG_{G})$
be the category of
$G(\sO)$-equivariant
perverse sheaves on $\GG_{G}$.
We will denote by $\IC_{\mu}=\IC(\barr{\GG_\mu})$
the intersection cohomology complex on the closure of the
orbit $\GG_{\mu}$.
There is a tensor product (convolution) construction $\ast$ \cite{MV1, MV2}
which makes
the category $\Perv_{G(\sO)}(\GG_{G})$
into a tensor category.
Let
$\Rep(\cG)$
be the category of rational representations of the
Langlands dual group $\cG$ over $\C$.

\sss{
Geometric Satake correspondence \cite[7.1]{MV2}
}
The global cohomology  functor
is an equivalence of  tensor categories
$H^*(\GG,-):\ (\Perv_{G}(\GG),\ast)
\con\ (\Rep(\cG),\ten)$.
Under this equivalence the
intersection cohomology sheaf $\IC_{\mu}=IC(\barr{\GG_\mu})$ corresponds to the
irreducible representation $V^\cG_\mu$ of $\cG$
with the highest weight $\mu$.


\sss{
Multiplicity spaces
$
\Hom_{GL_{m}}(V^\cG_{a_1}\otimes\dots\otimes V^\cG_{a_n},V^\cG_\la)
$
}
By Satake equivalence this space is realized  as
$
\Hom_{\Perv_\GO(\GG)}(\IC_{a_1}\ast\dots\ast \IC_{a_n},\IC_\la)
$.
The convolution
$\IC_{a_1}\ast\dots\ast \IC_{a_n}$ is the direct image
of the  IC-sheaf $\IC(\GG_{a_1}\ast\dots\ast \GG_{a_n})$
of the smooth simply connected space $\GG_{a_1}\ast\dots\ast \GG_{a_n}$
under the map
$\pi^a_\mu:\GG_{a_1}\ast\dots\ast \GG_{a_n}\to
\barr{\GG_\mu}
$ from \ref{convolutionlocal}\
(see \cite{MV2}).
The map
$\pi^a_\mu$ is
semismall
 and therefore
the multiplicity  space is the highest Borel-Moore homology
$
\HH[(\pi^a_\mu)^{-1}(\sL_{\la})]
$ of the fiber at $\sL_\la$ (see~\cite{MV2}).

\se{
Formulation of the Isomorphism Theorem
}
\label{mainresults}
Here we identify in type A the
quiver varieties
with certain varieties that arise for conjugacy
classes and for loop Grassmannians (theorem \ref{Isomorphism theorem}).

\sus{
Three settings: combinatorial data  and  spaces
}
\label{data}
The data
for $\gl_n$ quiver varieties
will be \
 $c,v,d$.
The data related to $\gl(N)$ conjugacy classes 
and to subvarieties of $GL(m)$ loop Grassmannians
 will be the same:\ $\la,\EE,\tmu,a,\sb$.

In \ref{qdatatogl}--\ref{gldatatoclass}
we pass from quiver data $d,v,c$ to the data $\la,\EE,\tmu,a,\sb$.
The geometric meaning of $\la,\EE,\tmu,a,\sb$  is explained 
for 
conjugacy classes at the end of 
\ref{gldatatoclass} and for loop Grassmannians 
in \ref{Rational convolution Grassmannians}.

\sss{
Quiver data
$(v,d)$ give  
weights $\la,\mu$ for $GL_m$ and $a$ for $GL_n$
}
\label{qdatatogl}
Let $d,v  \in\N^{n-1}$ be dimension vectors for representations
of Nakajima's $A_{n-1}$-quiver. Then let
\begin{enumerate}
\item
$N=\sum_{j=1}^{n-1}jd_j$ and  $m=\sum_{j=1}^{n-1}d_j$;
\item
 $\la,\cla\in\N^n$ be dual partitions
of $N$ where $
\cla_i=\sum_{j=i}^{n-1}d_j$ (hence $\cla_n=0$);
\item
$a\in\N^n$ be defined
using the Cartan matrix $C$ of type $A_{n-1}$,
cf. \cite[8.3]{N94}, by
\begin{equation}\label{a}
a_i=v_{n-1}+\sum_{j=i}^{n-1}(d-Cv)_j
;\end{equation}


\item
Let $\mu,\cmu$ be the dual partitions such that 
$\cmu$ is a reordering of $a$.
\end{enumerate}
Here,  $\cla$ can be viewed as a dominant  weight of $GL_n$ and $a$
is a weight of $V^{GL(n)}_{\cla}$,
cf. \ref{fromSLntoGLn}.

\sss{
Quiver datum $c$ gives  spectral datum $\sb\in\A^n$
}
\label{c to sb}
We will encode an element
$c=(c_1,\dots,c_{n-1})$ of the center of $\fg(V)$
(where
$\fg(V)=\prod_{i=1}^{n-1}\gl(V_i)$ and $\dim V_i=v_i$),
as an element
$\sb=
(b_1,\dots, b_n)
$ of $\A^n$ where
\begin{equation}
 b_i=c_1+\dots +c_{i-1}\ 
 (\text{hence}\ b_1=0).
\end{equation}


\sss{
Pair $ (a,\sb)$ is a resolution datum
for a conjugacy class
$\OO_{\EE,\tmu}\sub\gl(N)$
}
\label{gldatatoclass}


First $\sb$ gives $\EE=\{b_1,...,b_n\}$ and then  for  $e\in\EE$\
let $I_e\dff  \{i\in[1,n];\ b_i=e\}$ be the positions of $e$ in $\sb$.
Then let
$a^e\dff (a_i)_{i\in I_e}$ be the values of $a$ at these positions.
By ordering  $a^e$ we get the   partition $(\mu^e)^\sv$
and its dual is called $\mu^e$.
Then the bipartition  $\tmu$ is the family $(\mu^e)_{e\in\EE}$.

These  data $(\EE,\tmu)$ define a conjugacy class 
$\OO_{\EE,\tmu}$ (see \ref{conjclass}), where $\sx\in \OO_{\EE,\tmu}$
has eigenvalues
$\EE$ and at each  $e\in \EE$ the partition $\mu^e$ are sizes of Jordan blocks with eigenvalue $e$.
Also, there is a resolution    $\tgab$ of  $\barr{\OO_{\EE,\tmu}}$ 
and the 
the nilpotent conjugacy class $\OO_\la$
lies in the
closure of the Lusztig stratum 
of ${\OO_{\EE,\tmu}}$  (\ref{Deformation to a nilpotent orbit}.c).

\sss{
Rational convolution Grassmannians 
}
\lab{Rational convolution Grassmannians}
In the setting of loop Grassmannians
the same data $\la,\EE,\tmu,a,\sb$ as for the nilpotent slices
  arises
in the following way.

{\bf (i) }
We start in the rational Grassmannian $\GGrat_\A$.
For a dominant $\la\in\N^m$
we consider the $\A$-lattice
$\uLL_\la$ with the $\C$-form $\sU_\la=\pl_1^m z^{-\la_i}\C \bfe_i$.\
The space $\sN=\uLL_\la/\uLL_0=\sL_\la/\sL_0$ (for $\sL_\la=\uLL_\la|_\hz$) has dimension
$N=|\la|$. The 
 factorization $\bz$  of  $z$ to $\sN$ lies in the nilpotent conjugacy class $\OO_\la$
with a ``regular'' normal slice
 $T_\bz=
\bz+\Hom(\sL_\la/\sL_0,\sU_\la)$.
The
$\hy$-lattice $I_\la=z\inv \sU[z\inv]]$ is complementary to $\uLL_\la$. It
defines
the moduli
$\UU(\la)\sub\GG^{rat}_\A$
of $\A$-lattices $\LL$ that are complementary to $I_\la$,
contain the $\A$-lattice
$\uLL_0=U[z]$ and equal $U_\Gm$ on $\Gm$.

Consider $\EE\sub\C$ and an $\EE$-bipartition $\tmu$ of $N$
(see \ref{conjclass}).
Define $\UU(\la,\EE,\tmu)$ to consist of all $\LL\in\UU(\la)$
such that the action of  $z$ on $\LL/\uLL_0$
lies in the closure of the conjugacy class $\OO_{\EE,\tmu}\sub
\gl(\LL/\uLL_0)$
(see \ref{conjclass}).

\label{convolutionglobal}

{\bf (ii) }
Global convolution Grassmannians $\tGGra_\A$.
For $a\in\N^n$ with $\sum_{i=1}^n a_i=N$,
this is the space 
whose $\C$-points are
$n$-step flags of $\A$-lattices:
$$
\tGGra_\A
=\{ \LL_0\sub \LL_1\sub \dots \sub \LL_n
~|~ \LL_i\in\GGrat_\A \aand
\dim\LL_i/\LL_{i-1}=a_i
\text{ for }
1\leq i\leq n  \} ,
$$
where $\LL_0=U[z]$.
The map $\pi: \tGGra_\A\to \GGrat_\A$ is the projection to
$\LL_n$.
For 
$b\in \A^{n}$ let 
$
\tGGrab
\sub\tGGra_\A$ consist of filtrations such that $ z=b_i$ on 
$\LL_i/\LL_{i-1}$.


\sus{
Isomorphism theorem
}
\lab{Isomorphism theorem}
Now we can formulate our main result that relates quiver varieties, nilpotent slices
and loop Grassmannians. The translation of relevant data is in  \ref{data}.

\stheo{}
\label{main}
Consider  the above 3 kinds of data and  transitions between these data as in \ref{data}.
There are natural algebraic isomorphisms
$\phi, \tphi, \psi, \tpsi$
such that the following diagram commutes:

\begin{equation}
\begin{CD}
\fMc(v,d)
@>{\tphi}>{\simeq}>
\tgab|_{T_{\la}\cap{\barr{\OO_{\EE,\tmu}}}}\ 
\ @>{\tpsi}>{\cong}>\
\
\tGGrab_{\A}|_{\UU(\la,\EE,\tmu)}
\\
@V{p^c}VV
@V{\tmm}VV
@V{\pi}VV
\\
\fMc_\one (v,d)
\  @>{\phi}>{\simeq}>\
T_{\la}\cap \barr{\OO_{\EE,\tmu}}
@>{\psi}>{\cong}>
\UU(\la,\EE,\tmu)
\sub
\GGrat_\A
\end{CD}
\end{equation}

\sss{
Proof for isomorphisms $\psi,\tpsi$
}
Here we prove the claims for the right square.
The left square will be postponed to
section
\ref{proofIsomorphismTheorem}.

The middle column of the theorem concerns a  vector space space $\sN$ 
with a nilpotent operator $x$
and  a ``regular''normal  slice $T_\la$
to a nilpotent orbit $\OO_\la\sub \gl(\sN)$ at some $x \in \OO_\la$.
In the setting of the last  column (see \ref{Rational convolution Grassmannians}),
we will realize these data so that 
$\sN=\uLL_\la/\uLL_0$,\
$x$ is the factorization $\bz$  of $z$
and $T_\la$ is $T_\bz$.

Recall the isomorphism $\Psi:T_z\con
\cUUratI$ (theorem \ref{Slice Tz and UU ratI AGm})
where $I=I_\la$ from \ref{Rational convolution Grassmannians}).
Then the map $\psi$ is its restriction
to the subspace $T_\bz$ which consists of all operators $g\in T_z$
that equal $z$ on $\uLL_0$. This translates to the condition on $\LL=\Psi(g)\in\cUUratI$
that $\LL\subb\uLL_0$, so 
we have $\psi:T_\bz\con\UU(\la)$
in notation from
\ref{Rational convolution Grassmannians}).
Now, according to the theorem \ref{Slice Tz and UU ratI AGm}.b, the intersection
of $T_\bz=T_\la$ with the closure of the conjugacy class $\OO_{\EE,\tmu}$ in $\gl(\sN)$
corresponds  to 
$\UU(\la,\EE,\tmu)\sub \UU(\la)$.

Finally, the isomorphism $\psi$ lifts to the isomorphism $\tpsi$
since any $g$-invariant $n$-step flag in $\sN$ correspond by $\psi$ 
to an $n$-step flag of $\A$-lattices.
\qed

\sus{
Consequences
and Remarks
}
We comment on the case $c=0$,
consider two applications of our Isomorphism Theorem
and check an equality of dimensions.

\sss{
Case $c=0$
}
Here, we can describe the images of the maps
$\psi$ and $\tpsi$ and obtain a more precise result
stated in the introduction and \cite{MVyb}. In particular,
$(\psi\circ\phi)(0)=\sL_\la\in
\GG=\GG_\hz$,
and
$\tpsi\circ\tphi$
restricts to an isomorphism
\begin{equation}\label{fiberiso}
\tpsi\circ\tphi: \fL(v,d)\simeq\pi^{-1}(\sL_\la).
\end{equation}
We believe that one should be able to generalize these statements for arbitrary $c$.

\sss{
Compactifications of quiver varieties
}
The closure of $\UU(\la,\EE,\tmu)$
in $\GGrat_\A$
is a compactification of $\fMc_\one (v,d)$.
Analogously, the closure
of the image of $\fMc(v,d)$  in $\tGGra$
gives  a compactification of the quiver variety
$\fMc(v,d)$.

\sss{
A decomposition of the loop Grassmannian
}
\label{decompositionAffGrass}

The following is a corollary of the main theorem. Here $c=0$.

\scor{} We can decompose $\barr{\GG_\mu}$ into the following disjoint union:
\begin{equation}
\barr{\GG_\mu}=\bigsqcup_
{\substack{y\in G\cdot \sL_\la \\ \la\leq\mu  }}
\fMc_\zero(v,d)_y ,
\end{equation}
where $\la$ varies over the set of dominant coweights of $G$,
$G\cdot \sL_\la$ is the $G$-orbit of $\sL_\la$ in $\GG$, and
$\fMc_\zero(v,d)_y$ is a copy of quiver variety $\fMc_\zero(v,d)$
for every point $y\in G\cdot\sL_\la$,
with $v,d$ obtained from $\la,\mu$
by reversing the procedures of \ref{qdatatogl}.

\begin{proof}
The dominant cocharacters
$\la\in X^{+}_{*}(T)$
parameterize the
$G(\C[z\inv])$-orbits
$G(\C[z\inv])\cd\sL_\la$
in $\GG$
and the maps
$ L^{<0}G\tim\ G\cd  \sL_\la
\ \ra\
G(\C[z\inv])\cd\sL_\la$
are isomorphisms.
Then:
\begin{equation}
\barr{\GG_\mu}=\bigsqcup_{\substack{\la\in X^{+}_{*}(T)\\
y\in G\cdot \sL_\la}}
(L^{<0}G\cdot y)\cap\barr{\GG_\mu}
\cong
\bigsqcup_{\substack{\la\in X^{+}_{*}(T)\\
y\in G\cdot\sL_\la}} \fMc_\zero(v,d)_y
\end{equation}
since every $(L^{<0}G\cdot y)\cap\barr{\GG_\mu}$, for $y\in G\cdot\sL_\la$
is isomorphic to a copy of $\fMc_\zero(v,d)$.
\end{proof}

\sss{Remarks}
(1)
 The paper~\cite{BF1} extends
our construction to the  affine type A.

(2)
We would also like to mention another example of a decomposition
of an infinite Grassmannian into a disjoint union of quiver varieties.
Generalizing a result of G. ~Wilson \cite{W},
V. ~Baranovsky, V.~ Ginzburg,  and A.~ Kuznetsov \cite{BGK}
constructed a decomposition of (a part of) adelic Grassmannian
into a disjoint union of \emph{deformed} versions of quiver varieties
$\fMc(v,d)$ associated to affine quivers of type A.

\sss{
Dimensions for  $c=0$
}
We check that the varieties
$\fMc(v,d)$ and $\tT_{x}^{a}$ have the same dimension.
According to Nakajima \cite[Corollary 3.12]{N98}
$\fMc(v,d)$, if nonempty, is a smooth variety of dimension
${}^{t}v(2d-Cv)$
where $C$ is the Cartan matrix of type $A_{n-1}$.\fttt{
If $a _i > m $ for some $i$ then the smooth quiver variety is empty
by 
\ref{quiverdefs}.1.
}
If $\cla$ and $\cmu$ are defined by $v,d$ as in \ref{qdatatogl},
then we have
\begin{equation}
\begin{split}
\dim\fMc(v,d) & ={}^{t}v(2d-Cv)=
2\sum_{i=1}^{n-1}v_id_i-2\sum_{i=1}^{n-1}v_i^2
+2\sum_{i=1}^{n-2}v_iv_{i+1}  \\
& =\sum_{i=1}^{n-1}[(\cla_i)^2-(\cmu_i)^2]=
\dim\tT_{x,\mu}.
\end{split}
\end{equation}

\nc{\MMM}{\text{M}}
\se{
Quiver varieties and conjugacy classes of matrices
}
\label{quiverConjClassSection}

Here we   consider a particular case
of the Isomorphism Theorem. We choose $d=(\MMM, 0,\dots, 0)$
and $v=(v_1,\dots, v_{n-1})$ so that
$\MMM\geq v_1\geq v_2\geq\dots\geq v_{n-1}\ge 0$.
In 
\ref{Reduction} we will  explained 
how the  general case of the theorem reduces to this
special setting  by
the ``Main lemma'' from section~\ref{proofMainLemma}.

Define algebraic morphisms $\tphi: \fMc(v,d)\ra\tgab$ and
$\phi: \fMc_\one (v,d)\ra\barr{\OO_{\EE,\tmu}}$
by:
\begin{equation}
\begin{aligned}
\tphi: (x,\bx,p,q)& \mm\ \b(q_1p_1, \{0\}\sub\Ker (p_1)\sub\Ker( x_1p_1)
\sub\Ker (x_{n-1}\cddd x_1p_1)\b) , \\
\phi:(x,\bx,p,q)& \mm\ q_1p_1 .
\end{aligned}
\end{equation}

The following theorem is a common generalization of
results of \cite{KP} and \cite{N94},
cf. \cite{CB}.


\Theo{}
\label{gkp} The maps $\phi,\tphi$ defined above are
isomorphisms of algebraic varieties and the following diagram commutes
\begin{equation}
\begin{CD}
\fMc(v,d)
@>{\tphi}>>
\tgab
\\
@V{p}VV @V{\tmm}VV
\\
\fMc_\one (v,d)
@>{\phi}>>
\barr{\OO_{\EE,\tmu}}
\end{CD}
\end{equation}

\begin{proof}  Following the logic of \cite{N98, M},
it is not hard to check that $\tphi$ is a bijective
morphism between two smooth varieties of the same dimension
and thus an isomorphism. One checks directly that the map $\phi$ is a closed immersion
and it is surjective since both $p$ and $\tmm$ are surjective.
\end{proof}

\sss{}
In the case  when all the numbers
$0, c_1, c_1+c_2, \dots, c_1+c_2+\dots+ c_{n-1}$
are pairwise distinct, then the quiver variety $\fMc(v,d)$
is isomorphic to the conjugacy class of  semisimple elements
with eigenvalues $b_i=\sum_{1\le k<i}\ c_k,\ 1\le i\le n$, where
$b_i$ appears with multiplicity $a_i=v_{i-1}-v_i$
(we use $v_0=\MMM$ and $v_n=0$).

Quiver variety
$\fMc_\zero(v,d)$ is also isomorphic to a conjugacy class which
is generally different from the conjugacy class considered above.
The two classes coincide when the $SL(n)$ weight $d-Cv$ is dominant,
i.e. when $a_1\geq a_2\geq \dots \geq a_n$.

\se{
Main Lemma
}
\label{proofMainLemma}
This lemma provides an isomorphism
$\Phi: \Mc(v,d)\con S$ with a subspace $S$
of ``transversal'' elements
of
$\Mc(\tv,\td)$ for certain parameters $(\tv,\td)$.
This map restricts to an isomorphism
$\Phi^s: \Mc_s(v,d)\ra S^s$
of subspaces of stable elements.
The point is that the new parameters $\tv,\td$ are of the special kind considered in section
\ref{quiverConjClassSection}.

\sus{D'apr\` es Maffei}\label{dapresMaffei}
\label{mafNotation} We borrow Maffei's \cite{M} notations and conventions.
Let $v=(v_1,\dots,v_{n-1})$ and
$d=(d_1,\dots,d_{n-1})$ be two $(n-1)$-tuples of integers
and let us define $(n-1)$-tuples $\tv$ and $\td$ as follows:
\begin{equation}
\begin{split}
\td_1 &\dff \sum_{j=1}^{n-1}jd_j \aand \td_i \dff 0, \text{ for } i>1, \\
\tv_i &\dff v_i+\sum_{j=i+1}^{n-1}(j-i)d_j. \\
\end{split}
\end{equation}

Our goal is to construct a map from $\Mc(v,d)$ to $\Mc(\tv,\td)$,
that is we have to send a quadruple $(x,\bx,p,q)\in \Mc(v,d)$
to a quadruple $(\TA,\TB,\tga,\tde)\in \Mc(\tv,\td)$.
First of all,
the $I$-graded vector spaces $\TV_i$ and $\TD_i$ such that
$\dim\TV_i=\tv_i$ and $\TD_i=\td_i$
are constructed as follows. Let $D_j^{(k)}$ be a copy of $D_j$.
\begin{equation}
\begin{split}
\TD_1 &=\bigoplus_{1\leq k\leq j\leq n-1} D_j^{(k)}
\aand
\TD_i =0, \text{ for } i>1, \\
\TV_i &=V_i\oplus\bigoplus_{1\leq k\leq j-i\leq n-i-1} D_j^{(k)}.\\
\end{split}
\end{equation}
We consider the group $GL(V)$ as the subgroup of $GL(\tii V)$ acting as identity
on $D_j^{(k)} $. 
We need the following subspaces of $\TV_i$.

\begin{equation}
D _i ^{\prime} =
\bigoplus_{\substack{  i+1 \leq j \leq n-1 \\
1 \leq k \leq   j-i}}  D_j^{(k)},  \qquad
D _i ^{+} = \bigoplus_{\substack{i+2 \leq j \leq n-1 \\
2 \leq k \leq   j-i}} D_j^{(k)},   \qquad
D _i ^{-} = \bigoplus_{\substack{ i+2 \leq j \leq n-1 \\
1 \leq k \leq   j-i-1}}  D_j^{(k)}.
\end{equation}

In order to make the notation more homogeneous we set
$\TV_0\dff \TD_1$, $\TA_0=\tga_1$, $\TB_0=\tde_1$.

We will name the blocks of the maps $\TA_i$ and $\TB_i$ as follows
\begin{equation}\label{blocks}
\begin{aligned}
\pi_{D^{(h)}_j} \TA_i|_{D^{(h')}_{j'}}
&= {}^{i}\!t^{j',h'}_{j,h} &
\qquad   \pi_{D^{(h)}_j} \TB_i|_{D^{(h')}_{j'}}
&= {}^{i}\! s^{j',h'}_{j,h}  \\
\pi_{D^{(h)}_j} \TA_i|_{V_i } &= {}^{i}\! t^{V}_{j,h} &
\qquad \pi_{D^{(h)}_j} \TB_i|_{V_{i+1} } &= {}^{i}\!  s^{V}_{j,h} \\
\pi_{V_{i+1}}\TA_i|_{D^{(h')}_{j'}} &= {}^{i}\!  t^{j',h' }_{V} &
\qquad
\pi_{V_{i}}\TB_i|_{D^{(h')}_{j'}} &= {}^{i}\!s^{j',h'}_{V}
\end{aligned}
\end{equation}

We define also the following operator $z_i$ on $D_i^{\prime}$
\begin{equation}
z_i|_{D_j^{(1)}} = 0
\aand,                         \\
z_i|_{D_j^{(h)}} = Id_{D_j}: D_j^{(h)} \ra D_j^{(h-1)}
\end{equation}

\sss{}
Following Maffei let us introduce the following degrees:
\begin{equation}
\begin{aligned}
 \deg({}^{i}\!t^{j',h'}_{j,h} )  &=  \min(h-h'+1, h-h'+1+j'-j), \\
 \deg({}^{i}\!s^{j',h'}_{j,h} )  &=  \min(h-h', h-h'+ j'-j).
\end{aligned}
\end{equation}

\sss{} A quadruple $(\TA,\TB,\tga,\tde)\in \Mc(\tv,\td)$
is called \emph{transversal} if it satisfies the following two groups
of relations for $0\leq i\leq n-2$
\begin{enumerate}
\item first group (Maffei)
\begin{equation} \label{transone}
\begin{aligned}
   {}^{i}\!t^{j',h'}_{j,h}  &= 0
& &\mbox{ if } \deg(t^{j',h'}_{j,h}) < 0  \\
   {}^{i}\!t^{j',h'}_{j,h}  &= 0
& &\mbox{ if } \deg(t^{j',h'}_{j,h}) = 0
\mbox{ and } (j',h') \neq (j, h+1) \\
   {}^{i}\!t^{j',h'}_{j,h}  &= Id_{D_j}
& &\mbox{ if } \deg(t^{j',h'}_{j,h}) = 0
\mbox{ and } (j',h') = (j, h+1) \\
   {}^{i}\!t^{V}_{i,j,h}        &= 0         & &   \\
   {}^{i}\!t^{j',h'}_{V}    &= 0
& &\mbox{ if } h' \neq 1 \\
   {}^{i}\!s^{j',h'}_{j,h}  &= 0
& &\mbox{ if } \deg(s^{j',h'}_{j,h}) < 0  \\
   {}^{i}\!s^{j',h'}_{j,h}  &= 0
& &\mbox{ if } \deg(s^{j',h'}_{j,h}) = 0
\mbox{ and } (j',h') \neq (j, h) \\
   {}^{i}\!s^{j',h'}_{j,h}  &= Id_{D_j}
& &\mbox{ if } \deg(s^{j',h'}_{j,h}) = 0
\mbox{ and } (j',h') = (j, h) \\
   {}^{i}\!s^{V}_{j,h}        &= 0
& &\mbox{ if } h \neq j-i  \\
   {}^{i}\!s^{j',h'}_{V}    &= 0         & &
\end{aligned}
\end{equation}
\item second group
\begin{equation}\label{transtwo}
\pi_{D^{(h)}_j} \TB_i \TA_i|_{D^{(h')}_{j'}} - x_i=0
\qquad \text{ unless } h=j-i
\end{equation}
\end{enumerate}

Let us denote the set of all transversal elements in
$\Mc(\tv,\td)$ by $S$.
The set of all stable transversal elements is denoted by
$S^s=S\cap
\Mc_s(\tv,\td)$.

\sss{} We will need more notation. First of all denote
\begin{equation}
b^{i}_{j}=c_{i+2}+\dots + c_{j}\qquad \text{ for } -1\leq i\leq n-3,
\text{ and } i+2\leq j\leq n-1.
\end{equation}

Now we introduce some invariant polynomials of
$q_{i\ra j}p_{j\ra i}$ as follows.
First,
\begin{equation}
P(i,1,j)=q_{i+2\ra j}p_{j\ra i+2}
\end{equation}
and for $2\leq h'\leq j-i-1$
\begin{equation}
\begin{split}
P(i,h',j) & =q_{i+h'+1\ra j}p_{j\ra i+h'+1}
+(-1)^{j-i-h'-1}
\si_{j-i-h'}(b^{i}_{i+2},\dots,b^{i}_{j-1})
\\
 & +\sum_{k=1}^{j-i-h'-1} (-1)^k
\si_k(b^{i}_{i+2},\dots,b^{i}_{i+h'-1+k})
q_{i+h'+1+k\ra j}p_{j\ra i+h'+1+k}.
\end{split}
\end{equation}
where $\si_k$ is the $k$-th elementary symmetric function.


\sus{Main Lemma} We can now formulate our main lemma

%

We will always think of the maps $\tii A_i, \tii B_i$ as a block-matrix  with respect to
the given decomposition of $\tii V $, $\tii D$ and when we use a projection on
one of our subspaces, it will be a projection with respect  to the given 
decompositions. 
To simplify the notation we give a name to these blocks:

\slem\begin{enumerate}
\item[(i)] There exists a unique $G(V)$-equivariant map
$\Phi: \Mc(v,d)\ra S$
such that
\begin{equation}\label{basic}
\begin{aligned}
\pi_{V_{i+1}} \TA_i|_{V_i} &= x_i         &
\pi_{V_{i}} \TB_i|_{V_{i+1}}  &= \bx_i   \\
{}^{i}\!t^{i+1,1}_{V}      &=  p_{i+1} & \ \ \ \ \
{}^{i}\!s^{V}_{i+1,1}      &=  q_{i+1}
\end{aligned}
\end{equation}
\item[(ii)] The blocks of $\TA_i, \TB_i$ not defined in the equations
(\ref{transone}) and (\ref{basic})
are described as follows:
\begin{equation}
\begin{aligned}
{}^{i}\!t^{j',1}_{V}      &=  p_{j'\ra i+1} &\ \ \ \ \
{}^{i}\!s^{V}_{j,j-i}      &=  q_{i+1\ra j} \label{V} \\
\end{aligned}
\end{equation}
When $j'\neq j$ we have
\begin{equation}\label{tj}
\begin{aligned}
{}^{i}\!t^{j',h'}_{j,h} & =0  &
\text{ if } & (j',h') \neq (j, h+1) \\
{}^{i}\!s^{j',h'}_{j,h}  &= 0  &  \text{ if } &
(j',h') \neq (j, h) \text{ and } h\neq j-i
\end{aligned}
\end{equation}
and
\begin{equation}\label{sj}
{}^{i}\!s^{j',h'}_{j,j-i}=
q_{i+h'+1\ra j}p_{j'\ra i+h'+1}
\end{equation}
When $j=j'$ we have
\begin{equation}
{}^{i}\!t^{j,h'}_{j,h}  =
\begin{cases}
0, & \text{ if } h'=1 \\
(-1)^{h-h'+1}
\left(\begin{matrix}
h-1 \\ h'-2
\end{matrix}\right) c_{i+1}^{h-h'+1},
& \text{ if } 2\leq h'\leq h+1
\end{cases}
\end{equation}
And finally,
\begin{equation}
{}^{i}\!s^{j,h'}_{j,h}=
\begin{cases}
\left(\begin{matrix}
h-1\\ h'-1
\end{matrix}\right)
c_{i+1}^{h-h'}, &
\text{ if } h\neq j-i \\
P(i,h',j)+
\left(\begin{matrix}
h-1\\ h'-1
\end{matrix}\right)
c_{i+1}^{h-h'},
& \text{ if } 1\leq h'\leq h,\text{ and }  h=j-i
\end{cases}
\end{equation}
\item[(iii)] For $x\in\Mc(v,d)$ we have $\Phi(x)\in S^s$ if and only if
$x\in\Mc_s(v,d)$. Thus the restriction of $\Phi$ to the stable points
provides the $G(V)$-equivariant map $\Phi^s: \Mc_s(v,d)\ra S^s$
\item[(iv)] The maps $\Phi$ and $\Phi^s$ are isomorphisms of algebraic varieties.
\end{enumerate}

\begin{proof}
Following Maffei, we prove the lemma by decreasing induction on $i$.
If $i=n-2$ the maps $\TA_{n-2}$ and $\TB_{n-2}$ are completely defined by
the relations (\ref{basic}) and (\ref{transone}) and it is easy to see that
$\TA_{n-2}\TB_{n-2}=c_{n-1}$.

Assume that $\TA_k, \TB_k$ are defined for $k>i$
by the formulas in the lemma.

We have the following equations for
$\TA_i$ and $\TB_i$:
\begin{equation}\label{m=l+c}
\TA_i\TB_i =\TB_{i+1}\TA_{i+1} +c_{i+1}
\end{equation}
\begin{equation}\label{n=0}
\pi_{D^{(h)}_j} \TB_i \TA_i|_{D^{(h')}_{j'}} - z_i=0\
\text{ unless } h=j-i.
\end{equation}
Observe that
\begin{equation}\nonumber
\pi_{V_{i+1}} \TA_i \TB_i |_{V_{i+1}}= A_iB_i+p_{i+1}q_{i+1}
= B_{i+1}A_{i+1} +c_{i+1}= \pi_{V_{i+1}} \TB_{i+1} \TA_{i+1} |_{V_{i+1}}
+c_{i+1}.
\end{equation}
Then, in agreement with formulas (\ref{V})
\begin{equation}\nonumber
\begin{aligned}
\pi_{V_{i+1}} \TA_i \TB_i |_{D_{j}^{(h)}} =
\pi_{V_{i+1}} \TB_{i+1} \TA_{i+1} |_{D_j^{(h)}}=
\de_{h,1} B_{i+1} p_{j\ra i+2}=
\de_{h,1} p_{j\ra i+1}, \notag \\
\pi_{D_{j}^{(h)}} \TA_i \TB_i |_{V_{i+1}} =
\pi_{D_{j}^{(h)}} \TB_{i+1} \TA_{i+1} |_{V_{i+1}}=
\de_{h,j-i-1} q_{i+2\ra j} A_{i+1}
=\de_{h,j-i-1}q_{i+1 \ra j}.\notag
\end{aligned}
\end{equation}
where
$
\de_{p,q}
$
is the Kronecker symbol.

Now, in order to simplify the notation a bit we set
$t^{j',h'}_{j,h}\dff \ {}^{i}\!t^{j',h'}_{j,h}$
and $s^{j',h'}_{j,h}\dff \ {}^{i}\!s^{j',h'}_{j,h}$

Case I: $j\neq j'$. In this case
the equation (\ref{m=l+c}) and translates into the following
equations for $t^{j',h'}_{j,h}$ and $s^{j',h'}_{j,h}$:
\begin{equation}\label{m=lts}
s^{j',h'}_{j,h+1}+\sum_{\substack{h'<h''<h+1 \\ h'-j'<h''-j''<h+1-j}}
t^{j'',h''}_{j,h}s^{j',h'}_{j'',h''}+t^{j',h'}_{j,h}=
\begin{cases}
0, & \mbox { if } h\neq j-i-1 \\
q_{i+h'+1\ra j}p_{j'\ra i+h'+1}, & \mbox { if } h= j-i-1. \\
\end{cases}
\end{equation}
while the equation (\ref{n=0}) translates into the following
equations for $t^{j',h'}_{j,h}$ and $s^{j',h'}_{j,h}$, $h\neq j-i$:
\begin{equation}\label{n=0ts}
t^{j',h'}_{j,h}+\sum_{\substack{h'-1<h''<h \\
h'-1-j'<h''-j''<h-j}}s^{j'',h''}_{j,h}t^{j',h'}_{j'',h''}
+s^{j',h'-1}_{j,h}=0
\end{equation}

We claim that the system of equations (\ref{m=lts}) and (\ref{n=0ts})
has a unique solution indicated in the statement of the lemma.
We will prove this claim by
induction on $h$ and $h'$.

First of all, observe that from the
equation (\ref{n=0ts}) we have $t^{j',1}_{j,1}=0$.

We make two induction assumptions ($k\geq 1$):
\begin{enumerate}
\item  $t^{j',h'}_{j,h}=0$ for all $(h',h)$
such that $h'\leq h\leq k$
for all $j\neq j'$ at the same time.
\item  $s^{j',h'}_{j,h+1}=0$ for all $(h',h)$ such that
$h'< h\leq k+1\leq j-i$
for all $j\neq j'$ at the same time.
\end{enumerate}

Induction Step 1. Consider the equation (\ref{n=0ts}) for $h=k+1$.
By assumption (2) we have $s^{j',h'-1}_{j,k+1}=0$ and
$s^{j'',h''}_{j,k+1}=0$ for $j''\neq j$. If $j''=j$, then $j''\neq j'$
and by assumption (1) $t^{j',h'}_{j'',h''}=0$ for $h''\leq k$.
Now from equation (\ref{n=0ts}) we see that $t^{j',h'}_{j,k+1}=0$
for $h'\leq k+1$.

Induction Step 2. Consider the equation (\ref{m=lts}) for $h=k+1$.
By induction step (1) $t^{j',h'}_{j,k+1}=0$ and
$t^{j'',h''}_{j,k+1}=0$ for $j''\neq j$. If  $j''=j$, then $j''\neq j'$
and by assumption (2) $s^{j',h'}_{j'',h''}=0$.
Now from equation (\ref{m=lts}) we see that $s^{j',h'}_{j,k+2}=0$
for $h'< k+2$.

Finally, if $h+1=j-i$, then the equations (\ref{m=lts})
and the induction steps 1 and 2 yield:
\begin{equation}
s^{j',h'}_{j,j-i}=
q_{i+h'+1\ra j}p_{j'\ra i+h'+1}.
\end{equation}

Case II: $j=j'$. In this case we fix $j$ and
simplify the notation further a bit,
by setting $t^{h'}_{h}\dff t^{j,h'}_{j,h}$ and
$s^{h'}_{h}\dff s^{j,h'}_{j,h}$.
Now, taking into account Case I,
the equation (\ref{m=l+c}) and translates into the following
equations for $t^{h'}_{h}$ and $s^{h'}_{h}$:
\begin{equation}\label{m=ltsj}
s^{h'}_{h+1}+\sum_{h'<h''<h+1}
t^{h''}_{j,h}s^{h'}_{h''}+t^{h'}_{h}=
\begin{cases}
0, & \text { if } h\neq j-i-1 \text { and } h\neq h' \\
c_{i+1} & \text { if } h\neq j-i-1 \text { and } h= h' \\
P(i,h',j), & \mbox { if } h= j-i-1 \text { and } h\neq h'\\
P(i,h',j)+c_{i+1}, & \mbox { if } h= j-i-1 \text { and } h=h'\\
\end{cases}
\end{equation}
(In order to compute the right hand side, we need to use the
following combinatorial formula
$$
\si_a(c, c+b_1, \dots, c+b_p)
=\sum_{l=0}^a c^l
\left(\begin{matrix}
p-a+l+1 \\ l
\end{matrix}\right)
\si_{a-l}(b_1,\dots, b_p)
$$
for $a, p\in \Z$, $1\leq a\leq p$. We assume here that $\si_0(b_1,\dots, b_p) = 1$.)

The equation (\ref{n=0}) translates into the following
equations for $t^{j',h'}_{j,h}$ and $s^{j',h'}_{j,h}$, $h< j-i$:
\begin{equation}\label{n=0tsj}
t^{h'}_{h}+\sum_{h'-1<h''<h}
s^{h''}_{h}t^{h'}_{h''}
+s^{h'-1}_{h}=0
\end{equation}

Again, we claim that the system of equations (\ref{m=ltsj}) and
(\ref{n=0tsj})
has a unique solution indicated in the statement of the lemma.
Again, we will prove this claim by
induction on $h$ and $h'$.

First of all, observe that from the
equation (\ref{n=0tsj}) we have $t^{1}_{1}=0$.

We make two induction assumptions ($k\geq 1$):
\begin{enumerate}
\item $t^{h'}_{h}$ is given
by equations (\ref{tj})
 for all $(h',h)$
such that $h'\leq h\leq k$.
\item $s^{h'}_{h}$
is given
by equations (\ref{sj})
for all $(h',h)$ such that
$h'< h\leq k+1\leq j-i$.
\end{enumerate}

Proceeding by induction as in Case I
and using the formula (for $b, l\in \Z$, $0\leq b\leq l-2$)
$$
\sum_{a=b}^l (-1)^{l-a}
\left(
\begin{matrix}
l \\ a
\end{matrix}\right)
\left(\begin{matrix}
a+1 \\ b+1
\end{matrix}\right)
=0
$$
it is easy to see that all $t^{h'}_{h}$ and $s^{h'}_{h+1}$ are
given by formulas (\ref{tj}) and (\ref{sj}) respectively.

We have proved the assertions (i) and (ii)
of the lemma.
The assertion (iii) follows from the construction
and Lemma \ref{stablelemma} exactly as in \cite[Lemma 19]{M}.
The assertion (iv) follows from the construction, cf.
\cite[Lemma 19]{M}.\end{proof}

\sss{} It is important for us to record the formula
for $\TB_0\TA_0=\tde_1\tga_1$.
To simplify notation, we set
$$
b_{l}\dff b^{-1}_l=c_1+\dots + c_l,
$$
and
$$
P'(h',j)\dff \sum_{k=1}^{j-h'-1} (-1)^k
\si_k(b_{1},\dots, b_{h'-2+k})
q_{h'+k\ra j}p_{j\ra h'+k}
+(-1)^{j-h'-1}
\si_{j-h'}(b_{1},\dots,b_{j-1}).
$$

Now we have
\begin{equation}\label{dega1c}
(\tde_1\tga_1)^{j',h'}_{j,h}=
\begin{cases}
\Id_{D_j}, &\text{ if } h'=h+1,\ j'=j, \\
q_{h'\ra j}p_{j'\ra h'}+ \de_{j,j'}P'(h',j),
& \text{ if } h=j, \\
0, & \text{ otherwise }.
\end{cases}
\end{equation}
where
$\si_k$ is the $k$-th elementary symmetric function,
and we assume that the value of $\si_k$ at the empty
collection of variables is zero.

Finally, let us record the specialization of the above
formula for the case $c=0$.
Clearly, in this case $P'(h',j)=0$ and we have
\begin{equation}\label{dega1}
(\tde_1\tga_1)^{j',h'}_{j,h}=
\begin{cases}
\Id_{D_j}, &\text{ if } h'=h+1,\ j'=j, \\
q_{h'\ra j}p_{j'\ra h'},  & \text{ if } h=j, \\
0, & \text{ otherwise. }
\end{cases}
\end{equation}

\se{
Proof of the Isomorphism Theorem
}
\label{proofIsomorphismTheorem}

In this section we complete the proof of the Isomorphism Theorem (Theorem \ref{main}.)
The immersions $\psi$ and $\tpsi$
 were constructed in section \ref{globalpsi}.

\sus{
The isomorphisms $\phi$ and $\tphi$
}
The argument in this subsection is for the case $c=0$.
The argument for a general $c$ is completely
analogous. In the proof we mostly follow the logic of
\cite{M}.

\slem{} Let $(\TA,\TB,\tga,\tde)\in S$ and let $\tg\in G(\TV)$
be such that $\tg(\TA,\TB,\tga,\tde)\in S$. Then $\tg_i(V_i)\sub V_i$
and if we denote $g_i=\tg_i|_{V_i}$ we have
\begin{equation}
\tg(\TA,\TB,\tga,\tde)=g(\TA,\TB,\tga,\tde) .
\end{equation}

\begin{proof} The proof is lifted verbatim from \cite[Lemma 22]{M}.
\end{proof}

\sss{
Reduction to the situation from section
\ref{quiverConjClassSection}
}
\label{Reduction}
Let $D = \TD_1$ as in \ref{mafNotation}. Then $\dim D = N = \td_1 \dff \sum_{j=1}^{n-1}jd_j$.
Observe that $(\tv,\td)$ as constructed in \ref{mafNotation} must satisfy the conditions of section
\ref{quiverConjClassSection}
in order for $\fMc_\one(\tv,\td)$ and $\fMc(\tv,\td)$ to be nonempty, cf. \cite[1.4]{M}
and therefore, if nonempty,  $\fMc_\one(\tv,\td) \simeq \barr{\OO_\mu}$ and
$\fMc(\tv,\td) \simeq T^*\FF^a$ (for $c = 0$),
where $\mu, a$ are defined as in \ref{qdatatogl}.
(For a general $c$ the nilpotent orbit $\OO_\mu$ deforms into a general conjugacy class, cf.
\ref{Deformation to a nilpotent orbit}., \ref{gldatatoclass}.)
Now recall the definition (cf. \ref{nilpnormalslice}) of the transverse slice $T_x$
to the orbit $\OO_\la$ where $\la$ is obtained from $(v, d)$ as in \ref{qdatatogl}.
Let $T_{x,\mu} = T_x\cap\barr{\OO_\mu}$ be as in \ref{nilpnormalslice} and let $\tT_x^a$
be as in \ref{tildeta}.

\sss{}
\label{constructphi}
Now we will construct the maps $\phi_0$ and $\tphi$
completing the following commutative diagrams.

\begin{equation}
\begin{CD}
\Mc(v,d)@>{\Phi}>> S \\
@VVV @VVV \\
\fMc_0(v,d)@>{\phi_0}>> \fMc_0(\tv,\td) \\
\end{CD}
\qquad
\begin{CD}
\Mc_s(v,d)@>{\Phi^s}>> S^s \\
@VVV @VVV \\
\fMc(v,d)@>{\tphi}>>  \fMc(\tv,\td)\\
\end{CD}
\end{equation}

We denote $\phi\dff \phi_0|_{\fMc_\one (v,d)}:\fMc_\one (v,d)\to \fMc_\zero(\tv,\td)$.
Since
$\fMc_\one (\tv,\td) \simeq \barr{\OO_\mu}$
an element of $\fMc_\one (v,d)$ will be sent
by $\phi$ to an operator $y + f\in \End(D)$,
where y is nilpotent of type $\la$ and $f$ is given by the
\emph{explicit formulas} (\ref{dega1})
(and (\ref{dega1c}) for arbitrary $c$). A simple inspection
shows that $\Im\ \phi\sub T_{x,\mu}$, and $\Im\ \tphi\sub\tT_x^a$.

\slem{} The map $\phi$ is a closed immersion.

\begin{proof} It is enough to prove that $\phi_0$ is closed immersion.
Recall that
\begin{equation}
\begin{split}
\fMc_\zero(v,d) & =\Mc(v,d)//G(V)=\Spec \OO(\Mc(v,d))^{G(V)}, \\
\fMc_\zero(\tv,\td) & =\Mc(\tv,\td)//G(\TV)=\Spec \OO(\Mc(\tv,\td))^{G(\TV)}.
\end{split}
\end{equation}
We will prove that the restriction map
$\phi^*: \OO(\Mc(\tv,\td))^{G(\TV)}\to\OO(\Mc(v,d))^{G(V)}$
is surjective.

By Theorem \ref{invtheorem} the algebra
$\OO(\Mc(\tv,\td))^{G(\TV)}$ is generated by
$\widetilde\chi(\tde_1\tga_1)$ where
$\widetilde\chi$ is a linear form
on $\Hom(\TD_1,\TD_1)$. If $\tde_1\tga_1$ is of the form
(\ref{dega1}) and
$$
\widetilde\chi=\chi\in
\Hom(D_{j'}^{(h')},D_{j}^{(j)})^{*}\sub
\Hom(\TD_1,\TD_1)^{*},
$$
then for $1\leq h'\leq\min(j,j')$ we have
$$
\widetilde\chi(\tde_1\tga_1)=
\chi(\pi_{D_{j}^{(j)}}(\tde_1\tga_1)|_{D_{j'}^{(h')}})=
\chi(q_{h'\ra j}p_{j'\ra h'}),
$$
which are all the generators of the algebra $\OO(\Mc(v,d))^{G(V)}$
according to the Theorem \ref{invtheorem}.
\end{proof}

\slem{} The map $\tphi: \fMc(v,d) \to\tT_{x}^{a}$ is proper
and injective.

\begin{proof} We have the following diagrams
\begin{equation}\label{sdiagram}
\begin{CD}
\fMc(v,d)@>{\tphi}>> \tT_x^a \\
@V{p}VV @V{\mmm_a}VV \\
\fMc_\zero(v,d)@>{\phi_0}>> \fMc_\zero(\tv,\td)
\end{CD}
\qquad
\begin{CD}
\fMc(v,d)@>{\tphi}>> \tT_x^a \\
@V{p}VV @V{\mmm_a}VV \\
\fMc_\one(v,d)@>{\phi}>>  T_{x,\mu}
\end{CD}
\end{equation}
Since $\phi$
is a closed immersion and the morphisms $p$ and $\mmm_a$ are
projective, we see that $\tphi$ is proper.
Since all orbits in $\Mc_s(v,d)$
and $\Mc_s(\tv,\td)$ are closed, $\tphi$ is injective.
\end{proof}

\slem{} The map $\tphi: \fMc(v,d) \to\tT_x^a$ is an
isomorphism of algebraic varieties.

\begin{proof} Since $\tphi$ is a proper
injective morphism between connected smooth varieties of the same
dimension, $\tphi$ is an analytic isomorphism and therefore
an algebraic isomorphism.
\end{proof}

\slem{} The map $\phi:\fMc_\one(v,d)\to T_{x,\mu}$ is an isomorphism
of algebraic varieties.

\begin{proof} Since $\mmm_a$ is surjective, from the diagram
(\ref{sdiagram}) we see that $\phi$ is surjective.
Since $\phi$ is a surjective closed immersion,
and both $\fMc_\one(v,d)$ and $T_{x,\mu}$ are reduced varieties
over $\C$, $\phi$ is an algebraic isomorphism.
\end{proof}

\se{
Application to representation theory: $(\glm,\gln)$-duality
}
\label{applicationsRepTheory}
Here we observe that the relationship between quiver varieties and loop Grassmannians
provides a natural framework for $(GL_m,GL_n)$ duality.
In this section we use notation
$V=\C^m$, $W=\C^n$
and denote by
$V_\la=V^\glm_\la$ and
$W_{\cla}=W_\cla^\gln$
the irreducible representation of $\glm$and $\gl_n$
with highest weights $\la$ and $\cla$
(with weight spaces $V_\la(a)$ etc).

\sus{
Skew $(GL_m,GL_n)$ duality
}
For the $\gl_m\times \gl_n$ bimodule $V\otimes W$
we have the following decomposition \cite[4.1.1]{Ho}:
\begin{equation}\label{howew}
\wedge^N(V\otimes W)=\bigoplus_{\la}V_\la\otimes W_{\cla},
\end{equation}
where $\la$ goes through all partitions of $N$ which fit into the $n\times m$
box.

\sss{} Considering
$V\otimes W$ as a $\gl_m$ module $V\otimes\C^n$, we have the following
decomposition:
\begin{equation}\label{wedge}
\wedge^N(V\otimes W)=\bigoplus_{a_1+\dots +a_n=N}
\wedge^{a_1}V\otimes\dots\otimes\wedge^{a_n}V.
\end{equation}

Considered as a representation of the torus
$(\C^{\times})^n\sub \gl_n$ the vector space
$\wedge^{a_1}V\otimes\dots\otimes\wedge^{a_n}V$ has weight
$a=(a_1,\dots,a_n)$. Thus decompositions (\ref{howew}) and (\ref{wedge})
imply that
\begin{equation}\label{weightmult1}
\Hom_{\gl_m}
(\wedge^{a_1}V\otimes\dots\otimes \wedge^{a_n}V,
V_{\la}) \simeq W_{\cla}(a).
\end{equation}

\sss{
Geometric skew duality
}
We construct
a based version of the isomorphism (\ref{weightmult1}),
i.e., a geometric skew $(GL_n,GL_m)$ duality.
More precisely, with $N,v,d,a,\la$ as in \ref{qdatatogl},
we identify
the right hand side
with
$\HH(\pi^{-1} (\sL_\la))$, where $\sL_{\la}$ is a $\hz$-lattice in
the loop Grassmannian $\GG$,
and the left hand side
with
$\HH(\fL(v,d))$ by Theorem \ref{theonakajima}.
The identification of irreducible components
$\Irr\pi^{-1} (\sL_\la)=\Irr\fL(v,d)$,
which follows from the isomorphism (\ref{fiberiso})
matches the
natural basis of the space
of intertwiners
$\Hom_{GL_m}
(\wedge^{a_1}V\otimes\dots\otimes \wedge^{a_n}V,
V_{\la})$
arising from the loop Grassmannian construction
(i.e., $\Irr\pi^{-1} (\sL_\la)$), and the
natural basis of the weight space $W_{\cla}(a)$ in the Nakajima
construction
(i.e., $\Irr\fL(v,d)$). Altogether:
$$
\Hom_{GL_m}
(\wedge^{a_1}V\otimes\dots\otimes \wedge^{a_n}V,V_{\la})
\simeq\HH(\pi^{-1} (\sL_\la))
\simeq\HH(\fL(v,d))
\simeq W_{\cla}(a).
$$
Dually, we have
\begin{equation}\label{weightmult2}
\Hom_{\gl_n}
(\wedge^{c_1}W \otimes\dots\otimes \wedge^{c_m}W, W_{\cla})=V_{\la}(c)
.\end{equation}

\sus{Symmetric $(GL_m,GL_m)$ duality}
Analogously, if we consider the $N$-th symmetric power
$\Sym^N(V\otimes V)$
of the $\gl_m\times \gl_m$ bimodule $V\otimes V$, we have the following
 decomposition (a particular case of \cite[2.1.2]{Ho}):
\begin{equation}\label{howes}
\Sym^N(V\otimes V)=\bigoplus_{\la}V_\la\otimes V_{\la},
\end{equation}
where the sum is over all partitions $\la$ of $N$ with
at most $m$ parts.

Considering $V\otimes V$ as a $\gl_m$ module $V\otimes\C^m$, we
have the following
decomposition:
\begin{equation}\label{sym}
\Sym^N(V\otimes V)=\bigoplus_{c_1+\dots +c_m=N}
\Sym^{c_1}V\otimes\dots\otimes\Sym^{c_m}V.
\end{equation}
Thus decompositions (\ref{howes}) and (\ref{sym})
imply the following formula
\begin{equation}\label{weightmultsym}
\Hom_{\gl_m}
(\Sym^{c_1}V\otimes\dots\otimes \Sym^{c_m}V,
V_{\la})=V_{\la}(c),
\end{equation}
where $V_{\la}(c)$ is the weight space corresponding to weight $c$
of the $\gl_m$ highest weight module $V_{\la}$.

Combining the equations (\ref{weightmult2}) and
(\ref{weightmultsym}) we get
\begin{equation}\label{weightmultmix}
\Hom_{\gl_n}
(\wedge^{c_1}W \otimes\dots\otimes \wedge^{c_m}W,
W_{\cla})=
\Hom_{\gl_m}
(\Sym^{c_1}V\otimes\dots\otimes \Sym^{c_m}V,
V_{\la}).
\end{equation}

\sss{
Geometric symmetric duality
}
Geometry allows us to find a \emph{based}
isomorphism of the left and right hand side of (\ref{weightmultmix}).
We let $N,v,d,a,\la$ be as in \ref{qdatatogl}
except that we use notation $c$ for $a$.
First of all, the  quiver tensor product
theorems of Malkin \cite{Mal} and Nakajima \cite{N01b}
say that a natural basis
of the multiplicity space
in the left hand side of (\ref{weightmultmix})
is given by  the relevant irreducible components
$\Irr_{rel}[\tfg^{c,x}]$ of the ``Spaltenstein fiber''
 at $x$, defined by
$$
\tfg^{c,x}\ \dff\
\
\{F\in \FF^{c};\
x(F_i)\sub F_i
\text{ and } x \text{ acts on } F_i/F_{i-1}
\text{ as a regular nilpotent }\}.
$$

Now consider another convolution space
\begin{equation}
\begin{split}
\TGG^c & =\barr{\GG_{c_1\om_1}}\ast\cdots\ast\barr{\GG_{c_m\om_1}} \\
& =
\{\ L_0\sub\ L_1\sub\cddd\sub L_m ;\
\dim L_i/L_{i-1}=c_i,\  z|_{L_i/L_{i-1}}
\text{ is a regular nilpotent }\},
\end{split}
\end{equation}
where $\om_1$ is the first fundamental weight of $GL_m$.
We have a map $\pi: \TGG^c\to \GG$ defined by
$\pi:
(L_0\sub L_1\sub \dots \sub L_n)\mapsto L=L_n$.
Consider $\pi^{-1}(\sL_\la)$ for $\sL_\la\in\GG$.
It follows from the Geometric Satake Correspondence that
the set of relevant irreducible components 
$\Irr\pi^{-1} (\sL_\la)
$
indexes a basis in the right hand side of (\ref{weightmultmix}).

It is clear that the varieties
$\pi^{-1}(\sL_\la)
\simeq \tfg^{c,x}
$. This gives us a bijection
$\Irr \tfg^{c}_x\ \con\ \Irr\pi^{-1} (\sL_\la)$.
Summarizing, we have a based isomorphism
\begin{equation}\nonumber
\begin{split}
\Hom_{GL_n}
(\wedge^{c_1}W\otimes\dots\otimes \wedge^{c_m}W,
W_{\cla}) & \simeq \HH(\tfg^{c}_x)  \simeq\HH(\pi^{-1} (\sL_\la)) \\
& \simeq
\Hom_{GL_m}
(\Sym^{c_1}V\otimes\dots\otimes \Sym^{c_m}V,
V_{\la}).
\end{split}
\end{equation}

\sss{Remark}
The second author has greatly benefited from a class taught by W.~ Wang
at Yale \cite{Wa1}.
The ``geometric symmetric duality'' above has a
lot in common with the construction described in \cite{Wa2} and
we believe that the ``geometric skew duality'' construction answers a
question posed by Weiqiang Wang.

\section{
{Appendix by Vasily Krylov.}
Explicit isomorphism of quivers and loop Grassmannians
}

\label{Appendix}




In this section we use another approach to construct
the isomorphism between quiver varieties
and slices in the affine Grassmannian. We will denote this isomorphism by $\ze$.
We  prove that the isomorphism $\ze$ 
coincides with
the above morphism $\psi\ci\phi$. Our approach allows us to construct $\psi\ci\phi$  as an explicit morphism from quiver varieties to slices in the affine Grassmannian (note that by the construction $\phi$ is not the explicit morphism from $\mathfrak{M}^{0}_{{\bf{0}}}(v,d)$ to $T_\la \cap \overline{\mathcal{O}}_\mu$ since to define $\phi$ we embed both $\mathfrak{M}^{0}_{{\bf{0}}}(v,d)$, $T_\la \cap \overline{\mathcal{O}}_\mu$ inside some bigger variety and prove that their images in this variety coincide).
This approach also gives the explicit formula for the morphism $\psi\ci\phi$ (see theorem ~\ref{102}). Let us briefly explain geometric construction of the isomorphism $\ze$.

Our formula will actually
define  $\ze$ only on the  ``regular" part (certain open dense subvariety) of a quiver variety.
We check that here it coincides with the isomorphism $\psi\ci\phi$ and therefore
$\ze$ extends uniquely to the quiver variety.
The construction of $\ze$ is in some sense 2-dimensional
since over the  regular part we actually identify the above two moduli
with a certain moduli of vector bundles on $\P^2$.

Let $\mathfrak{M}(v,d)$ be a quiver variety with dimension vectors $v,d$, recall that $v,d$ are collections $v_1,\ldots,v_{n-1}$, $d_1,\ldots,d_{n-1}$ of nonnegative integer numbers. 
Let $\lambda, \mu$ be as in subsection~\ref{data}. Let variety $\mathfrak{M}_{\zero}^{reg}(v,d)$ be the regular part of $\mathfrak{M}(v,d)$ (see subsection ~\ref{111} for the definition). We are assuming that $\mathfrak{M}_{\zero}^{reg}(v,d) \neq \varnothing$ that is equivalent to $v_{i+1}+v_{i-1}+d_i-2v_i \geqslant 0$ for every $i$ (see~\cite[proposition~10.5 and remark~10.9]{N98}).  It then follows that the morphism $p\colon \mathfrak{M}^0(v,d) \ra \mathfrak{M}^0_{\bf{0}}(v,d)$ is surjective. Indeed $p$ is proper, hence its image is closed, on the other hand the image of $p$ contains open (so dense) subset  $\mathfrak{M}_{\zero}^{reg}(v,d)$ (see~\cite[proposition~3.24]{N98}). It follows that  $\mathfrak{M}^0_{\bf{0}}(v,d)=\mathfrak{M}^0_{\bf{1}}(v,d)$. Let us denote $V:=\bigoplus_{i} V_{i}$, $D:=\bigoplus_{i}D_{i}$. Let us consider $\mathfrak{M}_{\zero}^{reg}(V,D)$ that is the regular part of Gieseker  variety (see subsection~\ref{112} for the definition). Let us consider the natural inclusion from $\mathfrak{M}_{\zero}^{reg}(v,d)$ to  $\mathfrak{M}_{\zero}^{reg}(V,D)$. It follows from Nakajima's results that there exists a $\mathbb{C}^{\ast}$-action on $\mathfrak{M}(V,D)$ such that $\mathfrak{M}_{\zero}^{reg}(v,d)$ is a connected component of the fixed point subvariety $(\mathfrak{M}_{\zero}^{reg}(V,D))^{\mathbb{C}^{\ast}}$.

From the $ADHM$ description (see subsection ~\ref{112}) it follows that
the variety $\mathfrak{M}_{\zero}^{reg}(V,D)$ is isomorphic to the moduli space of vector bundles $E$ of rank $\operatorname{dim}D$ on $\mathbb{P}^{2}$ with fixed trivialization $D \otimes \mathcal{O}_{\ell_\infty} \iso E|_{\ell_\infty}$ at the line at infinity $\ell_\infty \subset \mathbb{P}^2$ and fixed second Chern class.

Thus $\mathfrak{M}_{\zero}^{reg}(v,d)$ identifies with the moduli space (to be denoted by $\operatorname{Bun}_{GL_{m},-\lambda}^{-\mu}(\mathbb{A}^{2}/\mathbb{G}_{m})$) of certain $\mathbb{C}^{\ast}$-equivariant vector bundles on $\mathbb{P}^{2}$ satisfying certain conditions (see Subsection ~\ref{401} for the precise statement). From \cite{BF1} it follows that we have an isomorphism
\begin{equation*}
\operatorname{Bun}_{GL_{m},-\lambda}^{-\mu}(\mathbb{A}^{2}/\mathbb{G}_{m}) \iso (L^{<0}GL_{m} \cdot \mathsf{L}_{\lambda} \cap L^{\geq 0}GL_{m} \cdot \mathsf{L}_{\mu}).
\end{equation*}
Thus we obtain the desired isomorphism $\mathfrak{M}_{\zero}^{reg}(v,d) \iso (L^{< 0}GL_{m} \cdot \mathsf{L}_\lambda \cap L^{\geq 0}GL_{m} \cdot \mathsf{L}_\mu)
$ as the composition
\begin{equation*}
\mathfrak{M}_{\zero}^{reg}(v,d) \iso \operatorname{Bun}_{GL_{m},\lambda}^{\mu}(\mathbb{A}^{2}/\mathbb{G}_{m}) \iso (L^{< 0}GL_{m} \cdot \mathsf{L}_\lambda \cap L^{\geq 0}GL_{m} \cdot \mathsf{L}_\mu).
\end{equation*}

The appendix is organized as follows.

In subsection~\ref{reg_part} we construct an isomorphism between regular parts of quiver varieties and intersections $L^{<0}GL_m \cdot \mathsf{L}_\la \cap L^{\geq 0}GL_m \cdot \mathsf{L}_\mu$ and compute it explicitly.

In subsection~\ref{121} we prove that the isomorphism we have constructed is given by the formula~\ref{102} below.

In subsection~\ref{301} we prove theorem~\ref{102} below by showing that the isomorphism we have constructed coincides with the morphism $\psi \circ \phi$ (restricted to the regular part).

\sus{}
Recall that $\mathcal{G} \simeq GL_{m}(\sK)/GL_{m}(\sO)$. For a cocharacter $\lambda$ of $GL_{m}$ we denote by $z^{\lambda}$ the class of $\operatorname{diag}(z^{\la_1},\ldots,z^{\la_m})$ in $T(\sK)/T(\sO) = \mathcal{G}_{T} \subset \mathcal{G}_{GL_{m}}$ where $T$ is the diagonal torus in $GL_{m}$.

\srem Note that the isomorphism $\mathcal{G} \simeq GL_{m}(\sK)/GL_{m}(\sO)$ identifies $\mathsf{L}_\la$ with $z^{-\la}$.

Recall the isomorphisms $\phi\colon\mathfrak{M}^{0}_{\zero}(v,d) \iso T_{\lambda} \cap \overline{\mathcal{O}}_{\mu}$ (see section 8,~\cite[subsection 3.3]{MVyb}) and $\psi=\psi_\mu\colon T_{\lambda} \cap \overline{\mathcal{O}}_{\mu} \iso L^{< 0}GL_{m} \cdot \mathsf{L}_\lambda \cap\overline{L^{\geq 0}GL_{m} \cdot \mathsf{L}_\mu}$
(see subsections~\ref{Theorem_iso!},~\ref{From operators to lattices},~\ref{normal slices in loop Grassmannians and nilpotent cones}) where $\lambda, \mu$ are as in Subsection~\ref{data}. 
It follows from (\ref{dega1}) that the isomorphism $\phi$ can be described as follows.
Recall a vector space $\widetilde{D}_1:=\bigoplus\limits_{1 \leq k \leq j \leq n-1} D^{(k)}_j$ of section~\ref{dapresMaffei}, here $D^{(k)}_j$ is a copy of $D_j$. Recall also the nilpotent operator $\mathsf{x}\colon \widetilde{D}_1 \ra \widetilde{D}_1$ that sends $D_{j}^{(k)}$ identically to $D_{j}^{(k-1)}$ for $k \geqslant 2$ and sends $D^{(1)}_j$ to zero.
For any $f_1 \in \operatorname{End}(\widetilde{D}_1)$ we consider its blocks $f^{j',h'}_{j,h}\colon D^{(h')}_{j'} \ra D^{(h)}_{j}$. Then we have
$\phi(x_i,\bar{x}_i,p_i,q_i)=\mathsf{x}+f_1 \in \mathcal{N}$, where
\begin{equation}\label{phi}
f^{j',h'}_{j,h} =  \begin{cases}
q_jx_{j-1}\ldots x_{h'+1}x_{h'}\bar{x}_{h'}\bar{x}_{h'+1}\ldots \bar{x}_{j'-1}p_{j'}~\text{if}~h=j,\\
0,~\text{otherwise}.
\end{cases}
\end{equation}

The isomorphism $\psi$ can be described as follows. Let us fix a basis $\{{\bf{e}}_1,\ldots,{\bf{e}}_m\}$ of $D$ such that $D_k=\operatorname{span}_{\mathbb{C}}({\bf{e}}_{1+d_{k+1}+\ldots+d_{n-1}},\ldots,{\bf{e}}_{d_{k}+\ldots+d_{n-1}})$.
Recall also the space $\widetilde{D}_1$ and set $D^{(k)}:=\bigoplus_{k \leq j \leq n-1}D_j^{(k)}$. Recall the basis $\{{\bf{e}}_1,\ldots,{\bf{e}}_m\}$ of $D$ and  denote by $\{{\bf{e}}_{k,1},\ldots,{\bf{e}}_{k,d_{n-1}+\ldots+d_k}\}$ the corresponding basis of $D^{(k)}$. We have the natural identification $\widetilde{D}_1 \iso \mathsf{L}_\la/\mathsf{L}_0$ that sends ${\bf{e}}_{k,j} \in D^{(k)}$ 
to the class of $z^{-k}{\bf{e}}_j$ in $L_\la/L_0$. Then $\psi$ is given by the following formula (compare with section~\ref{Slice Tz and UU ratI AGm}): 
\begin{equation*}
 \mathsf{x}+f_1 \mapsto \Big(1+\sum_{k=1}^\infty z^{-k}f_1(z+f_1)^{k-1}\Big)\mathsf{L}_{\la}.
\end{equation*}

\sss{\bf Theorem}
\label{102}
The isomorphism $\psi \circ \phi$ is given by the formula:
\begin{equation}
\label{72}
(x_{i},\bar{x}_{i},p_{i},q_{i}) \mapsto (1+z^{-1}\sum_{r,l=0}^\infty z^{-l}q \bar{x}^{r}x^{l}p)\mathsf{L}_\la.	
\end{equation}
where $(x,\bar{x},p,q):=(\oplus x_{i},\oplus \bar{x}_{i},\oplus p_{i},\oplus q_{i})$.

The proof will be given in subsection~\ref{12}.

\srem
Note that $1+z^{-1}\sum_{r,l=0}^\infty z^{-l}q \bar{x}^{r}x^{l}p $ can be compactly written as 
\begin{equation*}
1+q(\bar{x}-1)^{-1}(x-z)^{-1}p.
\end{equation*}

\sus{Geometric interpretation of the formula~(\ref{72})}\label{reg_part}
In this Section we give a geometric interpretation of the morphism from $\mathfrak{M}^{0}_{\zero}(v,d)$ to $L^{< 0}GL_{m} \cdot \mathsf{L}_\lambda \cap\overline{L^{\geq 0}GL_{m} \cdot \mathsf{L}_\mu}$ given by the formula~(\ref{72}). We prove that the variety $\mathfrak{M}^{reg}_{\zero}(v,d)$ can be identified with a certain moduli space $\operatorname{Bun}_{GL_{m},-\lambda}^{-\mu, \sum {v_{i}}}$
of $\mathbb{C}^{\ast}$-equivariant bundles on $\mathbb{P}^1\times \mathbb{P}^1$ with a fixed trivialization on $\mathbb{P}^{1}\times \infty \cup \infty \times \mathbb{P}^{1}$ (see section~\ref{119}). Using this identification we see that the morphism given by the formula~(\ref{72}) sends a vector bundle $E \in \operatorname{Bun}_{GL_{m},-\lambda}^{-\mu, \sum {v_{i}}}$ to a point $z^{-\lambda} \cdot E|_{\{1\} \times \mathbb{P}^1}$ of the affine Grassmannian $\GG$ (here one should note that the $\mathbb{C}^\ast$-equivariance of $E$ allows us to uniquely extend the trivialization of $E$ on $\infty \times (\mathbb{P}^{1} \setminus \{0\})$ to the trivialization of $E$ on $\mathbb{P}^1 \times (\mathbb{P}^1 \setminus \{0\}) $, hence, $E|_{\{1\} \times \mathbb{P}^1}$ defines a point of the affine Grassmannian $\GG$).

\sss{A variety $\mathfrak{M}_{\zero}^{reg}(v,d)$ and cocharacters $\lambda, \mu$}
\label{111}
A quadruple $(x,\bar{x},p,q) \in M^{0}(v,d)$ is called costable if for any $I$-graded subspace $V'$ of $V$ contained in $\operatorname{Ker}(q)$ and preserved by $x$ and $\bar{x}$ we have $V'=0$. Denote by $M^{0}_{reg}(v,d)$ the set of stable and costable quadruples in $M^{0}(v,d)$. Let $\mathfrak{M}_{\zero}^{reg}(v,d)$ \cite[Section 3]{N02} be the quotient of $M^{0}_{reg}(v,d)$ by the free action of $G(V)$. The morphism $\mathsf{p}:\mathfrak{M}^{0}(v,d) \rightarrow \mathfrak{M}_{\zero}^{0}(v,d)$ realizes $\mathfrak{M}_{\zero}^{reg}(v,d)$ as an open subvariety of $\mathfrak{M}_{\zero}^{0}(v,d)$.

Recall that $m=d_1+\ldots+d_{n-1}$, $D=\bigoplus_i D_i$. Recall the basis ${\bf{e}}_1,\ldots,{\bf{e}}_m$ of $D$ such that $D_k=\operatorname{span}_{\mathbb{C}}({\bf{e}}_{1+d_{k+1}+\ldots+d_{n-1}},\ldots,{\bf{e}}_{d_{k}+\ldots+d_{n-1}})$.
For the dimension vectors $(v,d)$ we define the following dominant cocharacters $\lambda, \mu$ of $GL(D)$: the cocharacter $\lambda$ acts with eigenvalue $t^{i}$ on $D_{i}$  and the cocharacter  $\mu$ acts with eigenvalue $t^{i}$ on the subspace of dimension $v_{i+1}+v_{i-1}+d_{i}-2v_{i}$, more detailed  in the basis ${\bf{e}_1},\ldots,{\bf{e}}_m$ the cocharacters $\mu,\, \la$ are given by 
\begin{equation*}
\mu:=(\underbrace{t^n,\ldots,t^n}_{v_{n-1}},\ldots,\underbrace{t^i,\ldots,t^i}_{v_{i+1}+v_{i-1}+d_i-2v_i},\ldots,\underbrace{1,\ldots,1}_{v_1}),
\end{equation*} 
\begin{equation*}
\la:=(\underbrace{t^{n-1},\ldots,t^{n-1}}_{d_{n-1}},\ldots,\underbrace{t^i,\ldots,t^i}_{d_{i}},\ldots, \underbrace{t,\ldots,t}_{d_1}).
\end{equation*}
We will sometimes consider $\la,\mu$ as elements of $\mathbb{Z}^m$: 
\begin{equation*}
\mu=(\mu_1,\ldots,\mu_m)=(\underbrace{n,\ldots,n}_{v_{n-1}},\ldots,\underbrace{i,\ldots,i}_{v_{i+1}+v_{i-1}+d_i-2v_i},\ldots,\underbrace{0,\ldots,0}_{v_1}),	
\end{equation*}
\begin{equation*}
\la=(\la_1,\ldots,\la_m)=(\underbrace{n-1,\ldots,n-1}_{d_{n-1}},\ldots,\underbrace{i,\ldots,i}_{d_{i}},\ldots, \underbrace{1,\ldots,1}_{d_1}).
\end{equation*}

Note that these definitions are compatible with the definitions of $\la,\mu$ given in subsection~\ref{data}. 


\sss{ADHM description}
\label{112}
Pick $u \in \mathbb{Z}_{\geqslant 0}$ and  $m \in \mathbb{Z}_{\geqslant 0}$. Let $V,\,D$ be vector spaces of dimensions $u$, $m$ respectively.

Let $\operatorname{Bun}_{GL_{m}}^{u}(\mathbb{A}^{2})$ denote the moduli space of rank $m$ vector bundles on $\mathbb{P}^{2}$  of second Chern class $u$ with a trivialization at the line at infinity $\ell_{\infty}$.

We set $\mathfrak{M}_{\zero}^{reg}(V,D)= \{(x,\bar{x},p,q) \in \mmm^{-1}(0)$ $|$ stable and costable $\}/GL(V)$,
where $(x,\bar{x},p,q)$ are Jordan quiver quadruples:

 
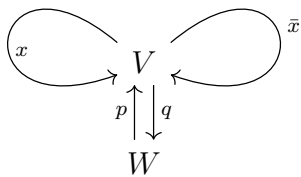
\begin{tikzcd}
V \ar[
    ,loop
    ,out=140
    ,in=200
    ,distance=5.5em
    ]{}{x} \ar[
    ,loop
    ,out=40
    ,in=-20
    ,distance=5.5em
    ]{}{\bar{x}} \arrow[xshift=0.7ex]{d}{q}  \\ \arrow[xshift=-0.7ex]{u}{p}
D
\end{tikzcd}

The ADHM description \cite[Theorem 2.1]{N99} identifies $\operatorname{Bun}_{GL_{m}}^{u}(\mathbb{A}^{2})$ with $\mathfrak{M}_{\zero}^{reg}(V,D)$.

The vector bundle $E_{(x,\bar{x},p,q)}$ corresponding to a quadruple $(x,\bar{x},p,q)$ can be obtained as the middle cohomology of the following monad:

\quad


\begin{tikzpicture}[node distance=0.4cm, auto]
\node (P) {$V \otimes \mathcal{O}_{\mathbb{P}^{2}}(-1)$};
\node(A)[right of=P] {$$};
\node(K)[right of=A] {$$};
\node(J)[right of=K] {$$};
\node(H)[right of=J] {$$};
\node(S)[right of=H] {$$};
\node(X)[right of=S] {$$};
\node(C)[right of=X] {$$};
\node(V)[right of=C] {$$};
\node(B)[right of=V] {$$};
\node(N)[right of=B] {$$};
\node(M)[right of=N] {$$};
\node(1)[right of=M] {$$};
\node(D)[right of=1] {$V \otimes \mathcal{O}_{\mathbb{P}^{2}}$};
\node(L)[above of=D] {$\oplus$};
\node(F)[above of=L] {$V \otimes \mathcal{O}_{\mathbb{P}^{2}}$};
\node(Z)[below of=D] {$\oplus$};
\node(G)[below of=Z] {$D \otimes \mathcal{O}_{\mathbb{P}^{2}}$};
\node(2)[right of=D] {$$};
\node(3)[right of=2] {$$};
\node(4)[right of=3] {$$};
\node(0)[right of=4] {$$};
\node(5)[right of=0] {$$};
\node(6)[right of=5] {$$};
\node(7)[right of=6] {$$};
\node(8)[right of=7] {$$};
\node(9)[right of=8] {$$};
\node(10)[right of=9] {$$};
\node(11)[right of=10] {$$};
\node(12)[right of=11] {$$};
\node(14)[right of=12] {$$};
\node(13)[right of=14] {$$};
\node(15)[right of=13] {$$};
\node(16)[right of=15] {$$};
\node(17)[right of=16] {$$};
\node(18)[right of=17] {$$};
\node(19)[right of=18] {$$};
\node(20)[right of=19] {$$};
\node(21)[right of=20] {$V \otimes \mathcal{O}_{\mathbb{P}^{2}}(1)$};
\draw[->] (D) --(21) node[below,midway] {${\bf{b}}$=$\begin{bmatrix}
    -(z_{0}\bar{x}-z_{2})  & z_{0}x-z_{1} & z_{0}p   \\

\end{bmatrix}$};
\draw[->] (P) --(D) node[below,midway] {${\bf{a}}$=$\begin{bmatrix}
    z_{0}x-z_{1}      \\
    z_{0}\bar{x}-z_{2}       \\
    z_{0}q
\end{bmatrix}$};
\end{tikzpicture}
the trivialization of $E_{(x,\overline{x},p,q)}$ at $\ell_{\infty}$ is induced by the morphism $\iota\colon D \otimes \mathcal{O}_{\ell_\infty} \ra (D \oplus V \oplus V) \otimes \mathcal{O}_{\ell_\infty}$ given by $s \mapsto (s,0,0)$ (note that the image of $\iota$ lies in the kernel of ${\bf{b}}|_{\ell_\infty}$ and is transversal to the image of ${\bf{a}}|_{\ell_\infty}$ so $\iota$ indeed induces the isomorphism $D \otimes \mathcal{O}_{\ell_\infty} \iso ({E_{(x,\overline{x},p,q)}})|_{\ell_\infty}$).

\sss{Torus fixed points} \label{401} Recall the notations of \cite[Subsection 4.4]{BF1}. Consider the action of $\mathbb{C}^{\ast}$ on $\mathbb{A}^{2}$ which sends $(x,y)$ to $(tx,t^{-1}y)$. Note that $GL_{m}\times \mathbb{C}^{\ast}$ acts on $\operatorname{Bun}_{GL_{m}}^{u}(\mathbb{A}^{2})$: the first factor acts by changing a trivialization at $l_{\infty}$ and the second factor via its action on $\mathbb{A}^{2}$. Every cocharacter $\eta \colon \mathbb{C}^{\ast} \rightarrow GL_{m}$ determines the diagonal action of $\mathbb{C}^{\ast}$ on $\operatorname{Bun}_{GL_{m}}^{u}(\mathbb{A}^{2})$. Let $\operatorname{Bun}_{GL_{m},\eta}^{u}(\mathbb{A}^{2}/\mathbb{G}_{m})$ denote the fixed point set of this action. The point $(0,0)\in\mathbb{A}^{2}$ is fixed under the $\mathbb{C}^{\ast}$-action, moreover note that the elements of $\operatorname{Bun}_{GL_{m}}^{u}(\mathbb{A}^{2})$ have no nontrivial automorphisms. So for every $E\in \operatorname{Bun}_{GL_{m},\eta}^{u}(\mathbb{A}^{2}/\mathbb{G}_{m})$ the group $\mathbb{C}^{\ast}$ acts on the fiber $E_{(0,0)}$ of $E$ at the point $(0,0)\in\mathbb{A}^{2}$ i.e. we have a homomorphism $\mathbb{C}^{\ast} \rightarrow GL(E_{(0,0)})$. For a cocharacter $\xi\colon \mathbb{C}^{\ast} \ra GL_m$ we denote by $\operatorname{Bun}_{GL_{m},\eta}^{\xi,u}(\mathbb{A}^{2}/\mathbb{G}_{m})$ the subvariety of $\operatorname{Bun}_{GL_{m},\eta}^{u}(\mathbb{A}^{2}/\mathbb{G}_{m})$ formed by all $E\in \operatorname{Bun}_{GL_{m},\eta}^{u}(\mathbb{A}^{2}/\mathbb{G}_{m})$ such that the conjugacy class of the homomorphism $\mathbb{C}^{\ast} \ra GL(E_{(0,0)})$ above coincides with the one for $\xi\colon \mathbb{C}^{\ast} \ra GL_m$.
It follows from~\cite[theorem~5.1(1)]{BF1} that the variety $\operatorname{Bun}_{GL_{m},\eta}^{\xi,u}(\mathbb{A}^{2}/\mathbb{G}_{m})$ is nonempty for only one value of $u$ so we will simply denote by $\operatorname{Bun}_{GL_{m},\eta}^\xi(\mathbb{A}^{2}/\mathbb{G}_{m})$ the corresponding variety.


\sss{Induced torus action on $\mathfrak{M}_{\zero}^{reg}(V,D)$}
 The diagonal $\mathbb{C}^{\ast}$-action on $\operatorname{Bun}_{GL_{m}}^{u}(\mathbb{A}^{2})$ corresponding to a cocharacter $\eta$ defines a $\mathbb{C}^\ast$-action on $\mathfrak{M}_{\zero}^{reg}(V,D)$ via the $ADHM$ isomorphism~\ref{112}.

\sss{\bf Lemma}
\label{115}
The action above can be described as follows (compare with~\cite[equation~(4.3)]{Hen}):
 $$(x,\bar{x},p,q) \mapsto (tx,t^{-1}\bar{x},p\eta(t)^{-1},\eta(t)q).$$

\begin{proof} Take  $t\in\mathbb{C}^{\ast}$. Consider the vector bundle $tE_{(x,\bar{x},p,q)}$ that is obtained from $E_{(x,\bar{x},p,q)}$ by the action of
$t$. It can be described as the middle cohomology of the following monad:
\quad

\begin{tikzpicture}[node distance=0.4cm, auto]
\node (P) {$V \otimes \mathcal{O}_{\mathbb{P}^{2}}(-1)$};
\node(A)[right of=P] {$$};
\node(K)[right of=A] {$$};
\node(J)[right of=K] {$$};
\node(H)[right of=J] {$$};
\node(S)[right of=H] {$$};
\node(X)[right of=S] {$$};
\node(C)[right of=X] {$$};
\node(V)[right of=C] {$$};
\node(B)[right of=V] {$$};
\node(N)[right of=B] {$$};
\node(M)[right of=N] {$$};
\node(6)[right of=M] {$$};
\node(D)[right of=6] {$V \otimes \mathcal{O}_{\mathbb{P}^{2}}$};
\node(L)[above of=D] {$\oplus$};
\node(F)[above of=L] {$V \otimes \mathcal{O}_{\mathbb{P}^{2}}$};
\node(Z)[below of=D] {$\oplus$};
\node(G)[below of=Z] {$D \otimes \mathcal{O}_{\mathbb{P}^{2}}$};
\node(Q)[right of=D] {$$};
\node(2)[right of=Q] {$$};
\node(3)[right of=2] {$$};
\node(4)[right of=3] {$$};
\node(7)[right of=4] {$$};
\node(5)[right of=7] {$$};
\node(8)[right of=5] {$$};
\node(10)[right of=8] {$$};
\node(9)[right of=10] {$$};
\node(11)[right of=9] {$$};
\node(12)[right of=11] {$$};
\node(13)[right of=12] {$$};
\node(14)[right of=13] {$$};
\node(15)[right of=14] {$$};
\node(16)[right of=15] {$$};
\node(17)[right of=16] {$$};
\node(18)[right of=17] {$$};
\node(22)[right of=18] {$$};
\node(21)[right of=22] {$$};
\node(20)[right of=21] {$$};
\node(19)[right of=20] {$V \otimes \mathcal{O}_{\mathbb{P}^{2}}(1)$};
\draw[->] (D) --(19) node[below,midway] {$\begin{bmatrix}
    -(z_{0}\bar{x}-tz_{2})  & z_{0}x-t^{-1}z_{1} & z_{0}p   \\

\end{bmatrix}$};
\draw[->] (P) --(D) node[below,midway] {$\begin{bmatrix}
    z_{0}x-t^{-1}z_{1}      \\
    z_{0}\bar{x}-tz_{2}       \\
    z_{0}q
\end{bmatrix}$};
\end{tikzpicture}

Let us emphasize that the trivialization of $tE_{(x,\bar{x},p,q)}$ at infinity is 
induced by 
\begin{equation*}
D \otimes \mathcal{O}_{\ell_\infty} \xrightarrow{\eta(t^{-1})}D \otimes \mathcal{O}_{\ell_\infty} \xrightarrow{s \mapsto (s,0,0)} (D \oplus V \oplus V) \otimes \mathcal{O}_{\ell_\infty}.
\end{equation*}

The following commutative diagram gives the isomorphism between the monad for $tE_{(x,\bar{x},p,q)}$ and the one for the quadruple
$(tx,t^{-1}\bar{x},p\eta(t)^{-1},\eta(t)q )$:

\begin{tikzpicture}[node distance=0.4cm, auto]
\node (P) {$V \otimes \mathcal{O}_{\mathbb{P}^{2}}(-1)$};
\node(A)[right of=P] {$$};
\node(K)[right of=A] {$$};
\node(J)[right of=K] {$$};
\node(H)[right of=J] {$$};
\node(S)[right of=H] {$$};
\node(X)[right of=S] {$$};
\node(C)[right of=X] {$$};
\node(V)[right of=C] {$V \otimes \mathcal{O}_{\mathbb{P}^{2}}$};
\node(B)[right of=V] {$$};
\node(N)[right of=B] {$$};
\node(M)[right of=N] {$\oplus$};
\node(1)[right of=M] {$$};
\node(D)[right of=1] {$$};
\node(2)[right of=D] {$V \otimes \mathcal{O}_{\mathbb{P}^{2}}$};
\node(3)[right of=2] {$$};
\node(4)[right of=3] {$$};
\node(5)[right of=4] {$\oplus$};
\node(6)[right of=5] {$$};
\node(7)[right of=6] {$$};
\node(8)[right of=7] {$D \otimes \mathcal{O}_{\mathbb{P}^{2}}$};
\node(9)[right of=8] {$$};
\node(10)[right of=9] {$$};
\node(11)[right of=10] {$$};
\node(12)[right of=11] {$$};
\node(13)[right of=12] {$$};
\node(14)[right of=13] {$$};
\node(15)[right of=14] {$$};
\node(16)[right of=15] {$$};
\node(17)[right of=16] {$$};
\node(18)[right of=17] {$V \otimes \mathcal{O}_{\mathbb{P}^{2}}(1)$};
\draw[->] (8) --(18) node[below,midway] {$$};
\draw[->] (P) --(V) node[below,midway] {$$};

\node(df)[below of=P] {$$};
\node(qr)[below of=df] {$$};
\node(qe)[below of=qr] {$$};
\node(qw)[below of=qe] {$$};
\node (P1)[below of=qw] {$V \otimes \mathcal{O}_{\mathbb{P}^{2}}(-1)$};
\node(A1)[right of=P1] {$$};
\node(K1)[right of=A1] {$$};
\node(J1)[right of=K1] {$$};
\node(H1)[right of=J1] {$$};
\node(S1)[right of=H1] {$$};
\node(X1)[right of=S1] {$$};
\node(C1)[right of=X1] {$$};
\node(V1)[right of=C1] {$V \otimes \mathcal{O}_{\mathbb{P}^{2}}$};
\node(B1)[right of=V1] {$$};
\node(N1)[right of=B1] {$$};
\node(M1)[right of=N1] {$\oplus$};
\node(20)[right of=M1] {$$};
\node(D1)[right of=20] {$$};
\node(21)[right of=D1] {$V \otimes \mathcal{O}_{\mathbb{P}^{2}}$};
\node(22)[right of=21] {$$};
\node(23)[right of=22] {$$};
\node(24)[right of=23] {$\oplus$};
\node(26)[right of=24] {$$};
\node(25)[right of=26] {$$};
\node(27)[right of=25] {$D \otimes \mathcal{O}_{\mathcal{P}^{2}}$};
\node(28)[right of=27] {$$};
\node(29)[right of=28] {$$};
\node(30)[right of=29] {$$};
\node(31)[right of=30] {$$};
\node(32)[right of=31] {$$};
\node(33)[right of=32] {$$};
\node(34)[right of=33] {$$};
\node(35)[right of=34] {$$};
\node(36)[right of=35] {$$};
\node(37)[right of=36] {$V \otimes \mathcal{O}_{\mathcal{P}^{2}}(1)$};
\draw[->] (27) --(37) node[below,midway] {$$};
\draw[->] (P1) --(V1) node[below,midway] {$$};

\draw[->] (P) --(P1) node[right,midway] {$Id$};
\draw[->] (V) --(V1) node[right,midway] {$t$};
\draw[->] (2) --(21) node[right,midway] {$t^{-1}$};
\draw[->] (8) --(27) node[right,midway] {$\eta(t)$};
\draw[->] (18) --(37) node[right,midway] {$Id$};

\end{tikzpicture}

\quad

 The morphism on cohomology induces the isomorphism between corresponding vector bundles with trivializations. \end{proof}

\sss{An isomorphism between $\coprod\limits_{v, \sum {v_{i}}=u} \mathfrak{M}_{\zero}^{reg}(v,d)$ and $\operatorname{Bun}_{GL_{m},-\lambda}^{u}$ }\qquad

Let $\coprod\limits_{v, \sum {v_{i}}=u} \mathfrak{M}_{\zero}^{reg}(v,d)$ be the disjoint union of quiver varieties of type A with vertices numbered by integers with a framing of dimension $d$ and $\sum {v_{i}}=u$. We define a morphism (see \cite[Section 3]{N02}): 
$$\tilde{\Theta}_{v,d}\colon \mathfrak{M}_{\zero}^{reg}(v,d) \rightarrow \mathfrak{M}_{\zero}^{reg}(V,D)\simeq \operatorname{Bun}_{GL_{m}}^{u}(\mathbb{A}^{2}),$$ $$\tilde{\Theta}_{v,d}(x_{i},\bar{x}_{i},p_{i},q_{i})=(\oplus x_{i},\oplus \bar{x}_{i},\oplus p_{i},\oplus q_{i}).$$

Morphisms $\tilde{\Theta}_{v,d}$ for different $v$ induce the morphism $\tilde{\Theta}_{d}\colon \coprod\limits_{v, \sum {v_{i}}=u} \mathfrak{M}_{\zero}^{reg}(v,d) \rightarrow \operatorname{Bun}_{GL_{m}}^{u}(\mathbb{A}^{2})$.

\sss{\bf Lemma}
\label{117}
$\tilde{\Theta}_{d}$ induces an isomorphism from $\coprod\limits_{v, \sum {v_{i}}=u} \mathfrak{M}_{\zero}^{reg}(v,d)$ to $\operatorname{Bun}_{GL_{m},-\lambda}^{u}$ where $\lambda$ is as in subsection~\ref{111} (see also section~\ref{data}).

\begin{proof} We describe the inverse morphism.
 Let $(x,\bar{x},p,q)$ be a fixed point under the $\mathbb{C}^{\ast}$-action on $\mathfrak{M}_{\zero}^{reg}(V,D)$ corresponding to $-\lambda$. Then using
lemma~\ref{115} we conclude that for every $t \in \mathbb{C}^{\ast}$ there exists
 $\rho_{V}(t) \in GL(V)$ such that

\quad

\begin{equation}
\label{73}
(tx,t^{-1}\bar{x},p \lambda(t),\lambda(t)^{-1}q ) =
(\rho_{V}(t)x\rho_{V}(t)^{-1},\rho_{V}(t)\bar{x}\rho_{V}(t)^{-1},\rho_{V}(t)p,q\rho_{V}(t)^{-1}).
\end{equation}

\quad

Note that $\rho_{V}(t)$ is uniquely determined by $t$ since the $GL(V)$-action on stable and costable quadruples is free. In particular $\rho_{V}$ defines a cocharacter of $GL(V)$. We decompose
$V$ into the direct sum $\bigoplus_i V_{i}$ where $V_{i}$ is the $t^{i}$-eigenspace of $\rho_{V}$ and similarly decompose $D$ into the direct sum $\bigoplus_i D_{i}$ with respect to $\lambda$ (recall that $D_i$ is the $t^{i}$-eigenspace of $\la$). The condition~(\ref{73})
implies that $\forall i \subset \mathbb{Z}$, $x(V_{i}) \subset V_{i+1}$, $\bar{x}(V_{i}) \subset V_{i-1}$, $p(D_{i}) \subset V_{i}$,
$q(V_{i}) \subset D_{i}$. So $(x,\bar{x},p,q)$ defines a point of a type $A$ quiver variety with vertices being numbered by integers such that
$\sum\limits_{i=-\infty}^{+\infty}v_{i}=u$, and the framing equal to $d$. The inverse map is constructed.
\end{proof}

\sss{An isomorphism between $\mathfrak{M}_{\zero}^{reg}(v,d)$ and $\operatorname{Bun}_{GL_{m},-\lambda}^{-\mu, \sum {v_{i}}} $}
\label{119}

\sss{\bf Lemma}
$\tilde{\Theta}_{d}$ induces the isomorphism  \cite[Section 4]{N02} \begin{equation*}
\Theta\colon \mathfrak{M}_{\zero}^{reg}(v,d) \iso \operatorname{Bun}_{GL_{m},-\lambda}^{-\mu, \sum {v_{i}}}
\end{equation*}
where  $\lambda, \mu$ are as in subsection~\ref{111} (see also section~\ref{data}).

\begin{proof}
 It is enough to prove that $\mathbb{C}^{\ast}$ acts on the fiber at the origin of $E \in \Theta(\mathfrak{M}_{\zero}^{reg}(v,d))$ by $\mu$. Recall the cocharater $\rho_V \colon \mathbb{C}^\ast \ra GL(V)$ that acts with eigenvalue $t^{i}$ on the space $V_{i}$.

\quad

   Let $(x,\bar{x},p,q):=\Theta(x_{i},\bar{x}_{i},p_{i},q_{i})$ and
 $E_{(x,\bar{x},p,q)}$ be the corresponding vector bundle. The bundle $tE_{(x,\bar{x},p,q)}$ is the middle cohomology of the following monad:

\quad

\begin{tikzpicture}[node distance=0.4cm, auto]
\node (P) {$V \otimes \mathcal{O}_{\mathbb{P}^{2}}(-1)$};
\node(A)[right of=P] {$$};
\node(K)[right of=A] {$$};
\node(J)[right of=K] {$$};
\node(H)[right of=J] {$$};
\node(S)[right of=H] {$$};
\node(X)[right of=S] {$$};
\node(C)[right of=X] {$$};
\node(V)[right of=C] {$$};
\node(B)[right of=V] {$$};
\node(N)[right of=B] {$$};
\node(M)[right of=N] {$$};
\node(1)[right of=M] {$$};
\node(D)[right of=1] {$V \otimes \mathcal{O}_{\mathbb{P}^{2}}$};
\node(L)[above of=D] {$\oplus$};
\node(F)[above of=L] {$V \otimes \mathcal{O}_{\mathbb{P}^{2}}$};
\node(Z)[below of=D] {$\oplus$};
\node(G)[below of=Z] {$D \otimes \mathcal{O}_{\mathbb{P}^{2}}$};
\node(2)[right of=D] {$$};
\node(3)[right of=2] {$$};
\node(4)[right of=3] {$$};
\node(4)[right of=4] {$$};
\node(5)[right of=4] {$$};
\node(6)[right of=5] {$$};
\node(7)[right of=6] {$$};
\node(8)[right of=7] {$$};
\node(9)[right of=8] {$$};
\node(10)[right of=9] {$$};
\node(11)[right of=10] {$$};
\node(12)[right of=11] {$$};
\node(14)[right of=12] {$$};
\node(13)[right of=14] {$$};
\node(15)[right of=13] {$$};
\node(16)[right of=15] {$$};
\node(17)[right of=16] {$$};
\node(18)[right of=17] {$$};
\node(19)[right of=18] {$$};
\node(20)[right of=19] {$$};
\node(21)[right of=20] {$V \otimes \mathcal{O}_{\mathbb{P}^{2}}(1)$};
\draw[->] (D) --(21) node[below,midway] {$\begin{bmatrix}
    -(z_{0}\bar{x}-tz_{2})  & z_{0}x-t^{-1}z_{1} & z_{0}p   \\

\end{bmatrix}$};
\draw[->] (P) --(D) node[below,midway] {$\begin{bmatrix}
    z_{0}x-t^{-1}z_{1}      \\
    z_{0}\bar{x}-tz_{2}       \\
    z_{0}q
\end{bmatrix}$};
\end{tikzpicture}

 The  map $(\rho_V(t),t^{-1}\rho_V(t) \oplus t\rho_V(t) \oplus \la(t),\rho_V(t))$ provides an
 isomorphism from $E_{(x,\bar{x},p,q)}$ to $tE_{(x,\bar{x},p,q)}$. In particular the inverse to this map induces a $\mathbb{C}^{\ast}$-module structure on the fiber at the origin of the  monad corresponding to $E_{(x,\bar{x},p,q)}$:

\quad

\begin{tikzcd}
V \ar[
    ,loop
    ,out=60
    ,in=120
    ,distance=2.5em
    ]{}{\rho_V^{-1}} \ar{rr}{} && V\oplus V\oplus D \ar[
    ,loop
    ,out=27
    ,in=153
    ,distance=7.5em
    ]{}{t^{-1}\rho_V^{-1}\oplus t\rho_V^{-1}\oplus \la^{-1}} \ar{rr} && V  \ar[
    ,loop
    ,out=60
    ,in=120
    ,distance=2.5em
    ]{}{\rho_V^{-1}}
\end{tikzcd}

Now we can calculate the action of $\mathbb{C}^{\ast}$ on the cohomology of this complex. The cocharacter
corresponding to this action is equal to the difference of the cocharacters on $V \oplus V \oplus D$ and the double cocharacter on $V$.
It follows that the desired cocharacter acts with an eigenvalue $t^{-i}$ on a subspace of dimension
$v_{i-1}+v_{i+1}+d_{i}-2v_{i}$. So this cocharacter is $-\mu$.
\end{proof}

\sss{An isomorphism $\operatorname{Bun}_{GL_{m},-\lambda}^{-\mu}(\mathbb{A}^{2}/\mathbb{G}_{m}) 
\simeq (L^{< 0}GL_{m} \cdot \mathsf{L}_\lambda \cap L^{\geq 0}GL_{m} \cdot \mathsf{L}_\mu)$.}
\label{1110}

\qquad

Recall the construction of the isomorphism (see \cite[theorem 5.2]{BF1})
$$
\Upsilon\colon \operatorname{Bun}_{GL_{m},-\lambda}^{-\mu}(\mathbb{A}^{2}/\mathbb{G}_{m}) \iso (L^{< 0}GL_{m} \cdot \mathsf{L}_\lambda \cap L^{\geq 0}GL_{m} \cdot \mathsf{L}_\mu).
$$

Recall the space $\operatorname{Bun}^u_{GL_m}(\mathbb{A}^2)$ of $m$-dimensional vector bundles on $\mathbb{P}^{2}$ with trivialization on the line at infinity and second Chern class equal to $u$. We will also denote this space by $\operatorname{Bun}^u_{GL_m}(\mathbb{P}^2;\ell_\infty)$. 
Let $\operatorname{Bun}^u_{GL_m}(\mathbb{P}^1 \times \mathbb{P}^1;\mathbb{P}^1 \times \{\infty\} \cup \{\infty\} \times \mathbb{P}^1)$ be the moduli space of $m$-dimensional  vector bundles on $\mathbb{P}^{1}\times\mathbb{P}^{1}$ with trivialization on the union
$(\mathbb{P}^{1}\times \{\infty\}) \cup (\{\infty\} \times \mathbb{P}^{1})$ of lines $\mathbb{P}^{1}\times \{\infty\},\, \{\infty\} \times \mathbb{P}^{1} \subset \mathbb{P}^1 \times \mathbb{P}^1$ and second Chern class $u$. 
We claim that there exists an isomorphism $\operatorname{Bun}^u_{GL_m}(\mathbb{P}^2;\ell_\infty) \simeq \operatorname{Bun}^u_{GL_m}(\mathbb{P}^1 \times \mathbb{P}^1;\mathbb{P}^1 \times \{\infty\} \cup \{\infty\} \times \mathbb{P}^1)$ (see~\cite[section~3.6]{Hen} and references therein).

The identification works as follows. Pick two distinct points $p_1,\, p_2 \in \ell_{\infty}$. Consider the blowup $\pi_1\colon \operatorname{Bl}_{p_1,p_2}(\mathbb{P}^2) \ra \mathbb{P}^2$ of $\mathbb{P}^2$ at the points $p_1,p_2$.  Consider also the blowup $\pi_2\colon \operatorname{Bl}_{\{\infty\} \times \{\infty\}}(\mathbb{P}^1 \times \mathbb{P}^1) \ra \mathbb{P}^1 \times \mathbb{P}^1$ of $\mathbb{P}^1 \times \mathbb{P}^1$ at the point $\{\infty\} \times \{\infty\}$. It is a classical result (see for example~\cite[example~7.22]{Ha}) that we have an isomorphism $\operatorname{Bl}_{p_1,p_2}(\mathbb{P}^2) \simeq \operatorname{Bl}_{\{\infty\} \times \{\infty\}}(\mathbb{P}^1 \times \mathbb{P}^1)$ that restricts to the isomorphism $\pi_1^{-1}(\ell_\infty) \simeq \pi_2^{-1}(\mathbb{P}^1 \times \{\infty\} \cup \{\infty\} \times \mathbb{P}^1)$. 
We set $Z:=\operatorname{Bl}_{p_1,p_2}(\mathbb{P}^2) \simeq \operatorname{Bl}_p(\mathbb{P}^1 \times \mathbb{P}^1)$, $Q:=\pi_1^{-1}(\ell_\infty) \simeq \pi_2^{-1}(\mathbb{P}^1 \times \{\infty\} \cup \{\infty\} \times \mathbb{P}^1)$.

So we obtain the following diagram 
\begin{equation*}
\begin{tikzcd}
 & Z  \arrow[dr, shift left=1, "\pi_2"] \arrow[dl, shift right=1, "\pi_1"] & 
\\\mathbb{P}^2 &&  \mathbb{P}^1 \times \mathbb{P}^1 
\end{tikzcd}
\end{equation*}

Let us denote by $\operatorname{Bun}_{GL_m}^u(Z;Q)$ the moduli space of $m$-dimensional  vector bundles $E$ on $Z$ with fixed trivialization  $D \otimes \mathcal{O}_Q \iso E|_Q$ and second Chern class equal to $u$.
Then the functors $\pi_1^*,\,\pi_2^*$ induce isomorphisms
\begin{equation*}
\operatorname{Bun}^u_{GL_m}(\mathbb{P}^2;\ell_\infty) \simeq  \operatorname{Bun}^u_{GL_m}(Z;Q) \simeq \operatorname{Bun}^u_{GL_m}(\mathbb{P}^1 \times \mathbb{P}^1,\mathbb{P}^1 \times \{\infty\} \cup \{\infty\} \times \mathbb{P}^1).
\end{equation*}
The composition of the isomorphisms above gives us the desired identification
\begin{equation*}
\operatorname{Bun}_{GL_m}^u(\mathbb{A}^2)=\operatorname{Bun}^u_{GL_m}(\mathbb{P}^2;\ell_\infty) \simeq \operatorname{Bun}^u_{GL_m}(\mathbb{P}^1 \times \mathbb{P}^1;\mathbb{P}^1 \times \{\infty\} \cup \{\infty\} \times \mathbb{P}^1).
\end{equation*}

The morphism $\Upsilon$  is constructed as follows (see~\cite[theorem 5.2]{BF1}). Pick a bundle $E \in \operatorname{Bun}_{GL_{m},-\lambda}^{-\mu}(\mathbb{A}^{2}/\mathbb{G}_{m})$ that we consider as a bundle on $\mathbb{P}^1 \times \mathbb{P}^1$ with a trivialization at $\mathbb{P}^1 \times \{\infty\} \cup \{\infty\} \times \mathbb{P}^1$. We claim that $E$ is trivial being restricted to $\mathbb{P}^{1} \times (\mathbb P^{1} \setminus \{0\})$. 
To see this note that $E$ is
trivial being restricted to the line $\mathbb{P}^{1} \times \{\infty\}$. The set of $p \in \mathbb{P}^1$ such that $E|_{\mathbb{P}^1 \times \{p\}}$ is nontrivial is closed and proper, hence, finite. 
So the number of horizontal ``jumping lines" of $E$ (i.e. lines $\ell$ such that the restriction $E|_{\ell}$ is nontrivial)
has to be finite. Using invariance of $E$ under $\mathbb{C}^{\ast}$-action we conclude that the only jumping line must be $\mathbb P^{1} \times \{0\}$ i.e. $E|_{\mathbb{P}^{1} \times (\mathbb P^{1} \setminus \{0\})}$ is trivial. 

Note now that we have a fixed trivialization of $E$ on the line $\mathbb{P}^{1} \times \{\infty\}$ that we can uniquely extend to the trivialization of $E|_{\mathbb P^{1} \times (\mathbb P^{1} \setminus \{0\})}$ (explicit formula for this trivialization appears in the proof of lemma~\ref{1112} see (\ref{triv_tau_formula!})). Then we restrict $E$ together with the trivialization above to the line
$\{1\} \times \mathbb{P}^{1}$ and get a point of the affine Grassmannian $\mathcal{G}$. Finally, we apply $z^{-\lambda}$ to this point and obtain the desired point $\Upsilon(E)$ of the slice $(L^{<0}GL_{m} \cdot \mathsf{L}_\lambda \cap L^{\geq 0}GL_{m} \cdot \mathsf{L}_\mu)$. It follows from~\cite[theorem 5.2]{BF1} that the map $\Upsilon$ that we described induces the isomorphism between $\operatorname{Bun}^{-\mu}_{GL_m,-\la}(\mathbb{A}^2/\mathbb{G}_m)$ and $(L^{<0}GL_m \cdot \mathsf{L}_\la \cap L^{\geq}GL_m \cdot \mathsf{L}_\mu)$\footnote{There is a typo in the proof of \cite[theorem 5.2]{BF1}: indeed recall that the loop rotation action in \cite[Section 2.1]{BF1} is defined by $g(s) \mapsto g(as)$ and  in the proof of~\cite[lemma 5.3]{BF1} it is claimed that if we take $x \in \mathcal{W}^{\overline{\la}}_{G,\overline{\mu}}$ and define $g(z)=x^{z^{-1}}$ then it satisfies $g(0) \in (\operatorname{Gr}^{\overline{\la}}_G)^{\mathbb{C}^*}$. This is not true in general since
\begin{equation*}
g(0)=\underset{\tau \ra 0}{\operatorname{lim}}\,x^{\tau^{-1}}=\underset{\tau \ra \infty}{\operatorname{lim}}\,x^{\tau}=s^{\overline{\mu}},
\end{equation*}
where the last equality holds by \cite[equation (5.2)]{BF1}. Note that the condition $s^{\overline{\mu}} \in (\operatorname{Gr}^{\overline{\la}}_G)^{\mathbb{C}^*}$ is equivalent to $\mu \in W\la$ that is not true in general (here $W$ is the Weyl group of $G$).
Claims above become true if we invert the loop rotation action i.e. if $a \in \mathbb{C}^\times$ acts via $g(s) \mapsto g(a^{-1}s)$ (then we should define $g(z)$ above as $x^{z}$).}

\sss{An isomorphism between $\mathfrak{M}_{\zero}^{reg}(v,d)$ and $(L^{< 0}GL_{m} \cdot \mathsf{L}_\lambda \cap L^{\geq 0}GL_{m} \cdot \mathsf{L}_\mu)$.} Composing the isomorphisms from subsections~\ref{119} and~\ref{1110} we obtain the isomorphism 
$$
\Upsilon \circ \Theta: \mathfrak{M}_{\zero}^{reg}(v,d) \iso (L^{< 0}GL_{m} \cdot \mathsf{L}_\lambda \cap L^{\geq 0}GL_{m} \cdot \mathsf{L}_\mu).
$$

\sss{\bf Lemma}
\label{1112}
The isomorphism $\Upsilon \circ \Theta$ is given by the following formula: $$(x_{i},\bar{x}_{i},p_{i},q_{i}) \mapsto (1+\sum_{}^\infty z^{-l}q \bar{x}^{r}x^{l}p)\mathsf{L}_\la,$$ where $(x,\bar{x},p,q):=(\oplus x_{i},\oplus \bar{x}_{i},\oplus p_{i},\oplus q_{i})$.

\begin{proof}(Similar to \cite[proposition 4.8]{Hen}) Take $(x_{i},\bar{x}_{i},p_{i},q_{i}) \in \mathfrak{M}_{\zero}^{reg}(v,d)$.

Let $E_{(x,\bar{x},p,q)} := \Theta(x_{i},\bar{x}_{i},p_{i},q_{i})$.
The vector bundle $E_{(x,\bar{x},p,q)}$ can be described as the middle cohomology of the following monad (see \cite[subsection 2.4]{BF2}):

\quad

\begin{tikzpicture}[node distance=0.4cm, auto]
\node (P) {$V \otimes \mathcal{O}_{\mathbb{P}^{1} \times \mathbb{P}^{1}}(-1,-1)$};
\node(A)[right of=P] {$$};
\node(K)[right of=A] {$$};
\node(J)[right of=K] {$$};
\node(H)[right of=J] {$$};
\node(S)[right of=H] {$$};
\node(X)[right of=S] {$$};
\node(C)[right of=X] {$$};
\node(V)[right of=C] {$$};
\node(B)[right of=V] {$$};
\node(N)[right of=B] {$$};
\node(M)[right of=N] {$$};
\node(1)[right of=M] {$$};
\node(D)[right of=1] {$V \otimes \mathcal{O}_{\mathbb{P}^{1} \times \mathbb{P}^{1}}(-1,0)$};
\node(L)[above of=D] {$\oplus$};
\node(F)[above of=L] {$V \otimes \mathcal{O}_{\mathbb{P}^{1} \times \mathbb{P}^{1}}(0,-1)$};
\node(Z)[below of=D] {$\oplus$};
\node(G)[below of=Z] {$D \otimes \mathcal{O}_{\mathbb{P}^{1} \times \mathbb{P}^{1}}$};
\node(2)[right of=D] {$$};
\node(3)[right of=2] {$$};
\node(4)[right of=3] {$$};
\node(4)[right of=4] {$$};
\node(5)[right of=4] {$$};
\node(6)[right of=5] {$$};
\node(7)[right of=6] {$$};
\node(8)[right of=7] {$$};
\node(9)[right of=8] {$$};
\node(10)[right of=9] {$$};
\node(11)[right of=10] {$$};
\node(12)[right of=11] {$$};
\node(21)[right of=12] {$V \otimes \mathcal{O}_{\mathbb{P}^{1} \times \mathbb{P}^{1}}$};
\draw[->] (D) --(21) node[below,midway] {${\bf{b}}$};
\draw[->] (P) --(D) node[below,midway] {${\bf{a}}$};
\end{tikzpicture}

$${\bf{a}} = \begin{bmatrix}
    sx-y\\
    h\bar{x}-z\\
    shq
\end{bmatrix} , {\bf{b}} = \begin{bmatrix}
    -(h\bar{x}-z),&sx-y,& p\\

\end{bmatrix},$$
where $([y:s],[z:h])$ are coordinates on $\mathbb{P}^{1} \times \mathbb{P}^{1}$.
Set $(\infty,\infty)$ := $([1:0],[1:0])$, $(0,0):=([0:1],[0:1])$.
Let us describe the trivialization of
   $E_{(x,\bar{x},p,q)}$ restricted to $\mathbb{P}^{1} \times (\mathbb P^{1} \setminus \{0\})$.
For this it suffices to construct a map $D\otimes\mathcal{O}_{\mathbb{P}^{1} \times (\mathbb P^{1} \setminus \{0\})} \rightarrow
\operatorname{Ker({\bf{b}})}|_{\mathbb{P}^{1} \times (\mathbb P^{1} \setminus \{0\})}$ transversal to $\operatorname{Im({\bf{a}})}|_{\mathbb{P}^{1} \times (\mathbb P^{1} \setminus \{0\})}$. It is easy to see that the map:

\begin{equation}\label{triv_tau_formula!}
\begin{tikzpicture}[node distance=0.6cm, auto]
\node (P) {$D \otimes \mathcal{O}_{\mathbb{P}^{1} \times (\mathbb P^{1} \setminus \{0\})}$};
\node(A)[right of=P] {$$};
\node(K)[right of=A] {$$};
\node(J)[right of=K] {$$};
\node(H)[right of=J] {$$};
\node(S)[right of=H] {$$};
\node(X)[right of=S] {$$};
\node(C)[right of=X] {$$};
\node(V)[right of=C] {$$};
\node(B)[right of=V] {$$};
\node(N)[right of=B] {$$};
\node(M)[right of=N] {$$};
\node(1)[right of=M] {$$};
\node(D)[right of=1] {$V \otimes \mathcal{O}_{\mathbb{P}^{1} \times(\mathbb P^{1} \setminus \{0\})}(-1,0)$};
\node(L)[above of=D] {$\oplus$};
\node(F)[above of=L] {$V \otimes \mathcal{O}_{\mathbb{P}^{1} \times (\mathbb P^{1} \setminus \{0\})}(0,-1)$};
\node(Z)[below of=D] {$\oplus$};
\node(G)[below of=Z] {$D \otimes \mathcal{O}_{\mathbb{P}^{1} \times (\mathbb P^{1} \setminus \{0\})}$};
\draw[->] (P) --(D) node[below,midway] {$\tau_{1}$ = $\begin{bmatrix}
    (h\bar{x}-z)^{-1}p\\
    0\\
    \operatorname{Id}
\end{bmatrix}$ };
\end{tikzpicture}
\end{equation}
satisfies the requirement.

 Note that $\tau_{1}$ is well defined
 because $h\bar{x}$ is nilpotent ($\bar{x}$ = $\oplus \bar{x}_{i}$, and $\bar{x}_{i}$ sends $V_{i}$ to $V_{i-1}$, so that $\oplus \bar{x}_{i}$ acts
 nilpotently on $\bigoplus V_{i}$), hence $h\bar{x}-z$ is invertible being restricted to $\mathbb{P}^{1} \times (\mathbb P^{1} \setminus \{0\})$ (since $z \neq 0$ on $\mathbb{P}^{1} \times (\mathbb P^{1} \setminus \{0\})$ and $h\bar{x}$ is nilpotent). We denote by $\bar{\tau}_1$ the trivialization induced by $\tau_1$.

For the same reasons the map:

$$\begin{tikzpicture}[node distance=0.6cm, auto]
\node (P) {$D \otimes \mathcal{O}_{(\mathbb{P}^{1} \setminus \{0\}) \times \mathbb P^{1}}$};
\node(A)[right of=P] {$$};
\node(K)[right of=A] {$$};
\node(J)[right of=K] {$$};
\node(H)[right of=J] {$$};
\node(S)[right of=H] {$$};
\node(X)[right of=S] {$$};
\node(C)[right of=X] {$$};
\node(V)[right of=C] {$$};
\node(B)[right of=V] {$$};
\node(N)[right of=B] {$$};
\node(M)[right of=N] {$$};
\node(1)[right of=M] {$$};
\node(D)[right of=1] {$V \otimes \mathcal{O}_{(\mathbb{P}^{1} \setminus \{0\}) \times \mathbb P^{1}}(-1,0)$};
\node(L)[above of=D] {$\oplus$};
\node(F)[above of=L] {$V \otimes \mathcal{O}_{(\mathbb{P}^{1} \setminus \{0\}) \times \mathbb P^{1}}(0,-1)$};
\node(Z)[below of=D] {$\oplus$};
\node(G)[below of=Z] {$D \otimes \mathcal{O}_{(\mathbb{P}^{1} \setminus \{0\}) \times \mathbb P^{1}}$};
\draw[->] (P) --(D) node[below,midway] {$\tau_{2}$ = $\begin{bmatrix}
    0\\
    (y-sx)^{-1}p\\
    \operatorname{Id}
\end{bmatrix}$ };
\end{tikzpicture}$$

induces the trivialization of $E_{(x,\bar{x},p,q)}$ restricted to $(\mathbb{P}^{1} \setminus \{0\}) \times \mathbb{P}^1$ to be denoted $\bar{\tau}_2$. Note that these two trivializations agree at the point $(\infty,\infty)$ and extend the trivialization of
$E_{(x,\bar{x},p,q)}$ restricted to $\{\infty\} \times \mathbb{P}^1$, $\mathbb{P}^1 \times \{\infty\}$ respectively. Now we can construct $\Upsilon(E_{(x,\bar{x},p,q)})$. To this end we have to calculate the transition function $(\bar{\tau}_{1}^{-1} \circ \bar{\tau}_{2})_{[1:1] \times (\mathbb{P}^{1}\setminus \{ 0,\infty
\})}$ that is the point of $\mathcal{G}$ corresponding to $(E_{(x,\bar{x},p,q)})|_{[1:1] \times \mathbb{P}^1}$ and the trivialization induced by 
$$
\tau_{1}\colon D \otimes \mathcal{O}_{\mathbb{P}^1 \times (\mathbb{P}^1 \setminus \{0\}) } \iso (E_{(x,\bar{x},p,q)})|_{\mathbb{P}^1 \times (\mathbb{P}^1 \setminus \{0\})}.
$$

Let us compute $\bar{\tau}_{1}^{-1} \circ \bar{\tau}_{2}$ on the fiber at a point $([y_0:s_0],[z_0:h_0])=p$. On the fiber of $\tau_{1}, \tau_{2}$ at $p$ we have the following morphisms:

$$\begin{tikzpicture}[node distance=0.3cm, auto]
\node (P) {$D$};
\node(G)[right of=P] {$$};
\node(H)[right of=G] {$$};
\node(J)[right of=H] {$$};
\node(K)[right of=J] {$$};
\node(L)[right of=K] {$$};
\node(Q)[right of=L] {$$};
\node(W)[right of=Q] {$V$};
\node(E)[right of=W] {$\oplus$};
\node(D)[right of=E] {$V$};
\node(21)[right of=D] {$\oplus$};
\node(T)[right of=21] {$D$};
\node(Y)[right of=T] {$$};
\node(U)[right of=Y] {$$};
\node(I)[right of=U] {$$};
\node(O)[right of=I] {$$};
\node(22)[right of=O] {$$};
\node(23)[right of=22] {$$};
\node(B)[right of=23] {$D$};
\draw[->] (P) --(W) node[below,midway] {$(\tau_{1})|_{p}$};
\draw[->] (B) --(T) node[below,midway] {$(\tau_{2})|_{p}$};
\end{tikzpicture}$$
They induce isomorphisms:
$$\begin{tikzpicture}[node distance=0.7cm, auto]
\node (P) {$D$};
\node(G)[right of=P] {$$};
\node(H)[right of=G] {$$};
\node(J)[right of=H] {$$};
\node(K)[right of=J] {$$};
\node(L)[right of=K] {$$};
\node(Q)[right of=L] {$$};
\node(D)[right of=Q] {$(\operatorname{Ker({\bf{b}})}/\operatorname{Im({\bf{a}})})|_{p}$};
\node(Y)[right of=D] {$$};
\node(U)[right of=Y] {$$};
\node(I)[right of=U] {$$};
\node(O)[right of=I] {$$};
\node(22)[right of=O] {$$};
\node(23)[right of=22] {$$};
\node(B)[right of=23] {$D$};
\draw[->] (P) --(D) node[below,midway] {$(\bar{\tau}_1)|_p$};
\draw[->] (B) --(D) node[below,midway] {$(\bar{\tau}_2)|_p$};

\end{tikzpicture}$$

 For a vector $w \in D$ we want to find 
 $(\bar{\tau}_{1}^{-1})|_p \circ (\bar{\tau}_{2})|_p(w)$ i.e. 
 the vector $\tilde{w} \in D$ such that $(\tau_{1})|_p(\tilde{w}) -
 (\tau_{2})|_p(w) \in \operatorname{Im({\bf{a}})}$. It means that there exists a vector $u\in V$ such that $(\tau_{1})|_p(\tilde{w}) - (\tau_{2})|_p(w) =
 {\bf{a}}|_p(u)$. We obtain the system of equations:

\quad
\begin{equation}
\label{74}
\begin{cases}
(h_0\bar{x}-z_0)^{-1}p(\tilde{w}) = s_0x(u)-y_0u \\
(s_0x-y_0)^{-1}p(w) = h_0\bar{x}(u)-z_0u\\
\tilde{w}-w = s_0h_0q(u)
\end{cases}
\Rightarrow
\end{equation}

\quad

\begin{equation}
\begin{cases}
u = (h_0\bar{x}-z_0)^{-1}(s_0x-y_0)^{-1}p(w)\\
\tilde{w} = w + s_0h_0q(u)
\end{cases}
\end{equation}

\quad

Hence $\tilde{w} = w + s_0h_0q(h_0\bar{x}-z_0)^{-1} (s_0x-y_0)^{-1} p(w)$.

So

$$(\bar{\tau}_{2}^{-1} \circ \bar{\tau}_{1})_{[1:1] \times (\mathbb{P}^{1}\setminus \{ 0,\infty\})} = (1 +
q (\bar{x}-z)^{-1}(x-1)^{-1}p).$$

Thus we obtained the point $1 + z^{-1}\sum\limits_{r,l=0}^{\infty}z^{-r}q\bar{x}^{r}x^lp$ of $\mathcal{G}$.
It remains to multiply it by $z^{-\lambda}$. 
We conclude that

   $$\Upsilon \circ \Theta(x_{i},\bar{x}_{i},p_{i},q_{i}) = z^{-\lambda}(1 +
  z^{-1}\sum\limits_{r,l=0}^{\infty}z^{-r}q \bar{x}^{r}x^{l}p)=(1+\sum_{}^\infty z^{-l}q \bar{x}^{r}x^{l}p)\mathsf{L}_\la,$$
 the last equality holds since for $w \in D_i$ we have $q\bar{x}^rx^lp(w) \in D_{i+l-r}$ so
 \begin{equation*}
 z^{-\la}(z^{-r}q\bar{x}^rx^lp)z^{\la}(w)=z^{-l} q\bar{x}^rx^lp(w).
 \end{equation*}
  \end{proof}

 \srem{} Note that from~(\ref{74}) it follows that $w = \tilde{w} - s_0h_0q(s_0x-y_0)^{-1}(h_0\bar{x}-z_0)^{-1}p(\tilde{w})$. So conjugating by $z^\la$ we get the following nontrivial equation (compare with~\cite[proof of proposition 4.10]{Hen}):
  \begin{equation}
\label{76}
(1 +
  z^{-1}\sum\limits_{r,l=0}^{\infty}z^{-l}q\bar{x}^{r}x^{l}p)^{-1}=(1 -
  z^{-1}\sum\limits_{r,l=0}^{\infty}z^{-l}qx^{l}\bar{x}^{n}p).\end{equation}

\sus{Proof of Theorem~\ref{102}}
\label{12}

\sss{\bf Lemma}
\label{121}
The isomorphism $\psi \circ \phi$ restricted to $\mathfrak{M}_{\zero}^{reg}(v,d)$ induces an isomorphism between $\mathfrak{M}_{\zero}^{reg}(v,d)$ and $(L^{< 0}GL_{m} \cdot \mathsf{L}_\lambda \cap L^{\geq 0}GL_{m} \cdot \mathsf{L}_\mu)$ and is given by the formula:  \begin{equation}(x_{i},\bar{x}_{i},p_{i},q_{i}) \mapsto (1 +
  z^{-1}\sum\limits_{r,l =0}^{\infty}z^{-l}q\bar{x}^{r}x^{l}p)\mathsf{L}_\la\end{equation} where $(x,\bar{x},p,q):=(\oplus x_{i},\oplus \bar{x}_{i},\oplus p_{i},\oplus q_{i})$.

\begin{proof} In Subsection~\ref{reg_part} we proved that the map $$\Upsilon \circ \Theta:(x_{i},\bar{x}_{i},p_{i},q_{i}) \mapsto (1 +
  z^{-1}\sum\limits_{r,l =0}^{\infty}z^{-l}q\bar{x}^{r}x^{l}p)\mathsf{L}_\la$$ is an isomorphism between $\mathfrak{M}_{\zero}^{reg}(v,d)$ and $(L^{< 0}GL_{m} \cdot \mathsf{L}_\lambda \cap L^{\geq 0}GL_{m} \cdot \mathsf{L}_\mu)$. Recall that $\mathcal{G}$ is the moduli space of lattices $L\subset D(K)$. Set 
  \begin{equation*}
  L:= (1 +
  z^{-1}\sum\limits_{r,l =0}^{\infty}z^{-l}q\bar{x}^{r}x^{l}p)\mathsf{L}_\la \in (L^{<0}GL_m \cdot \mathsf{L}_\la \cap L^{\geq 0}GL_m \cdot \mathsf{L}_\mu).	
  \end{equation*}

 According to the proof of proposition in subsection~\ref{Isomorphism of a pronilpotent slice to z and a neighborhood of sL} 
this lattice $L$ is uniquely determined by a $\C$-linear map $f\colon \mathsf{L}_{\la} \rightarrow 
 \mathcal{I}_{\la}$ (here $ 
 \mathcal{I}_\la=z^{-1}\mathrm{U}_\la[z^{-1}]$, $\mathsf{U}_{\la}=\oplus_i z^{-i}D_i$ compare with Subsection~\ref{UU si})
 and $f$ is uniquely determined by $f_{1}\colon \mathsf{L}_{\la} \rightarrow \mathsf{U}_{\la}$.

 Note that $f$ can be constructed as follows. Denote the projection of $\mathsf{L}_{\la}\oplus \mathcal{I}_{\la}
 $ to $\mathsf{L}_{\la}$ (resp. 
 $\mathcal{I}_\la$) along 
 $\mathcal{I}_\la$
 (resp. $\mathsf{L}_{\la}$) by $\pi_{\la}$ (resp. $\pi_{\la}^{-}$). Note that $\pi_{\la}$ induces the isomorphism $\pi\colon L \iso \mathsf{L}_{\la}$. Then $$f:=\pi_{\la}^{-} \circ \pi^{-1}.$$
We have the following commutative diagram:



\begin{equation} \label{78}
\begin{tikzcd}
\mathsf{L}_{\la} \arrow[rr, "1+z^{-1}\sum\limits_{r,l=0}^{\infty}z^{-l}q\bar{x}^{r}x^{l}p"] \arrow[ddrr] && \mathsf{L}_{\la} \oplus 
\mathcal{I}_\la
\arrow[dd, shift left=1, "\pi_{\la}"] \arrow[rr, "\pi_{\la}^{-}"] && 
\mathcal{I}_\la
\\
&&&&
\\
&& \mathsf{L}_{\la}\arrow[uu, shift left=1, "\pi^{-1}"]  \arrow[uurr, "f"] &&
\end{tikzcd}
\end{equation}

Let $f=\sum\limits_{k=1}^{\infty} z^{-k}f_{k}$, $f_k \in \operatorname{Hom}(\mathsf{L}_\la,\mathsf{U}_\la)$ as in the proof of Proposition in Subsection~\ref{Isomorphism of a pronilpotent slice to z and a neighborhood of sL}. For a vector $z^{-h}e$, $e \in D_{j'}$,

\begin{equation*}
\pi^{-1}(z^{-h}e)=z^{-h}e+\sum\limits_{k=1}^{\infty} z^{-k}f_{k}(z^{-h}e)=z^{-h}e+\sum\limits_{j} z^{-j-1}w_{j}+\sum\limits_{k=2}^{\infty} z^{-k}f_{k}(z^{-h}e)
\end{equation*}
for some $w_{j} \in D_{j}$ (that we want to compute).

 Recall that by~(\ref{76}) the map $$1 - z^{-1}\sum\limits_{r,l=0}^{\infty}z^{-l}qx^{l}\bar{x}^{r}p\colon \mathsf{L}_{\la}\oplus 
 \mathcal{I}_\la
 \rightarrow \mathsf{L}_{\la}\oplus 
 \mathcal{I}_\la
 $$ is inverse to the map 
 $$
 1 +
  z^{-1}\sum\limits_{r,l=0}^{\infty}z^{-l}q\bar{x}^{r}x^{l}p \colon \mathsf{L}_{\la}\oplus 
 \mathcal{I}_\la
  \rightarrow \mathsf{L}_{\la}\oplus 
 \mathcal{I}_\la
  .$$
It follows from the commutativity of the diagram~(\ref{78}) that 

   \begin{equation}\label{commut_diag!}
   \Xi_{h,j'} := (1 - z^{-1}\sum\limits_{r,l=0}^{\infty}z^{-l}qx^{l}\bar{x}^{r}p)(z^{-h}e+\sum\limits_{j} z^{-j-1}w_{j}+\sum\limits_{k=2}^{\infty} z^{-k}f_{k}(z^{-h}e)) \in \mathsf{L}_{\la}.
   \end{equation}
   
Note that we have two gradings on $D(\mathsf{K})$. One is by the degree of $z$ and the other comes from the decomposition $D=\oplus_i D_i$.  
Straightforward  computation shows that $(-j-1,j)$-component of the vector $\Xi_{h,j'}$ is $z^{-j-1}(w_j-qx^{j-h}\bar{x}^{j'-h}p(e))$ (to prove this we observe that the sum of degrees of components of vectors $z^{-k}f_k(z^{-h}e)$ is equal to $-k$ that is less than $-1$ for $k>1$ and also note that the operator $z^{-l-1}qx^l\bar{x}^rp$ shifts the sum of degrees by $-r-1$). Using~(\ref{commut_diag!}) together with the fact that vectors of $\mathsf{L}_\la$ do not have any components of degree $(-j-1,j)$ we conclude that 

 $$w_{j}=q_{j}x_{j-1}\dots x_{h}\bar{x}_{h}\dots \bar{x}_{j'-1}p_{j'}(e_{j'}).$$
 From that and (\ref{phi}) it follows that the maps $f_{1}$ for $L$ and for  $\psi \circ \phi(x_{i},\bar{x}_{i},p_{i},q_{i})$ are the same. So

 $$\psi \circ \phi(x_{i},\bar{x}_{i},p_{i},q_{i})=L=\eta \circ \Theta(x_{i},\bar{x}_{i},p_{i},q_{i}).$$ \end{proof}

\subsubsection{Proof of Theorem~\ref{102}} \label{301}
\begin{proof} By Lemma~\ref{121}, the morphism $\psi \circ \phi$ restricted to the dense open subvariety $\mathfrak{M}_{0}^{reg}(v,d) \subset \mathfrak{M}_{0}(v,d)$ is given by the formula $$(x_{i},\bar{x}_{i},p_{i},q_{i}) \mapsto (1 +
  z^{-1}\sum\limits_{r,l =0}^{\infty}z^{-l}q\bar{x}^{r}x^{l}p)\mathsf{L}_\la.$$

Now continuity of the map $\psi \circ \phi$ implies Theorem~\ref{102}.
\end{proof}


\begin{thebibliography}{MVi2}

\bibitem[BGK]{BGK}
V.~ Baranovsky, V.~ Ginzburg, and A.~ Kuznetsov,
\emph{Wilson's Grassmannian and a noncommutative quadric},
Int. Math. Res. Not. 2003, no. 21, 1155--1197.

\bibitem[BD]{BD}
A.~ Beilinson and V.~ Drinfeld,
\emph{Quantization of Hitchin's integrable system
and Hecke eigensheaves}, available at http://math.uchicago.edu/mitya/langlands.html.


\bibitem[BF1]{BF1}
A.~Braverman and M.~Finkelberg,
\emph{Pursuing the double affine Grassmannian I: transversal slices via instantons on $A_k$-singularities},
Duke Math. J. 152 (2010), no. 2, 175-206.

\bibitem[BF2]{BF2}
A.~Braverman and M.~Finkelberg,
\emph{Pursuing the double affine Grassmannian III: Convolution with affine Zastava},
Mosc. Math. J. 13 (2013), no. 2, 233-265.

\bibitem[BG]{BG}
A.~Braverman and D.~Gaitsgory \emph{On Ginzburg's Lagrangian construction
of representations of ${\rm
GL}(n)$}, Math. Res. Lett. 6 (1999), no. 2, 195--201.

\bibitem[BGV]{BGV}
A. ~Braverman, D.~ Gaitsgory, and M. Vybornov,
\emph{Relation between two geometrically defined bases in representations of $GL_n$},
preprint 2004, math.RT/0411252.




\bibitem[CG]{CG}
N.~ Chriss and V.~Ginzburg,
\emph{Representation theory and complex geometry},
Birkh\" auser, Boston, 1997.

\bibitem[CB]{CB}
W.~ Crawley-Boevey, \emph{Normality of Marsden-Weinstein reductions for representations of quivers},
Math. Ann. 325 (2003), no. 1, 55--79.





\bibitem[Ha]{Ha} 
J. Harris, 
\emph{Algebraic Geometry: A First Course}, book, Springer, 1992.




\bibitem[Hen]{Hen}
A. ~Henderson,
\emph{Involutions on the affine Grassmannian and moduli spaces of principal bundles},
preprint 2015, arXiv:1512.04254.

\bibitem[Ho]{Ho}
R.~Howe,
\emph{Perspectives on invariant theory: Schur duality,
multiplicity-free actions and beyond},
The Schur lectures (1992) (Tel Aviv), 1--182, Israel Math.
Conf. Proc., 8, Bar-Ilan Univ., Ramat Gan, 1995.

\bibitem[IMW]{IMW}
Mee Seong Im, Chun-Ju Lai, Arik Wilber,
\emph{Springer fibers via quiver varieties using Maffei-Nakajima isomorphism}.
arXiv:2009.08778.

\bibitem[KP]{KP}
H. ~ Kraft and C.~ Procesi,
{Closures of conjugacy classes of matrices are normal},
Invent. Math. 53 (1979), no. 3, 227--247.

\bibitem[L1]{L81}
G.~Lusztig, Green polynomials and singularities of unipotent classes, Adv. in Math. 42 (1981),
no. 2, 169--178.

\bibitem[L2]{Lu}
Singularities, character formulas and a q-analog of weight multiplicities, Asterisque, 101–-102 (1983), 208–-229.

\bibitem[L3]{L98}
G.~Lusztig, On quiver varieties, Adv. in Math. 136 (1998), 141--182.

\bibitem[L4]{L00}
G.~Lusztig, Semicanonical bases arising from enveloping algebras, Adv. Math. 151 (2000), no. 2,
129–-139





\bibitem[Maf]{M}
A.~Maffei, \emph{Quiver varieties of type A},
Comment. Math. Helv. 80 (2005), 1--27.

\bibitem[Mal]{Mal}
A.~Malkin,
\emph{Tensor product varieties and crystals: GL case},
Trans. Amer. Math. Soc. 354 (2002), no. 2, 675--704.

\bibitem[Mir]{Mir}
I.~Mirkovi\'c, \emph{Character sheaves on reductive Lie algebras},
Mosc. Math. J. 4 (2004), no.4, 897--910, 981.

\bibitem[MVi1]{MV1}
I.~ Mirkovi\'c, and K.~ Vilonen,
\emph{Perverse sheaves on affine Grassmannians and Langlands duality},
Math. Res. Lett. 7 (2000), no. 1, 13--24.

\bibitem[MVi2]{MV2}
I.~ Mirkovi\'c, and K.~ Vilonen,
\emph{Geometric Langlands duality and representations of algebraic groups over commutative rings}, Annals Mathematics 166.1 (2007): 95--143.



\bibitem[MVy]{MVyb}
I. ~ Mirkovi\' c and M.~ Vybornov,
\emph{On quiver varieties and affine Grassmannians of type A},
C. R. Acad. Sci. Paris, Ser. I  (2003) 336 (3) 207--212.

\bibitem[MV]{MV}
I. ~ Mirkovi\' c and M.~ Vybornov,
{\em Quiver varieties and Beilinson-Drinfeld Grassmannians of type A}.
arXiv:0712.4160. 




\bibitem[N1]{N94}
H.~Nakajima,  \emph{Instantons on ALE spaces, quiver varieties,
and Kac-Moody algebras},
Duke Math. J. 76 (1994), no. 2,365--416.

\bibitem[N2]{N98}
H.~Nakajima, \emph{Quiver varieties and Kac-Moody algebras},
Duke Math. J. 91 (1998), no. 3, 515--560.


\bibitem[N3]{N01b}
H. ~Nakajima,
\emph{Quiver varieties and tensor products},
Invent. Math. 146 (2001), no. 2, 399--449.

\bibitem[N4]{N99}
H. ~Nakajima,
\emph{Lectures on Hilbert schemes of points on surfaces},
University Lecture Series 18, American Mathematical Society, Providence, RI, 1999.

\bibitem[N5]{N02}
H. ~Nakajima,
\emph{Geometric construction of representations of affine algebras},
Proceedings of the International Congress of Mathematicians, Vol. I (Beijing, 2002), 423--438.


\bibitem[Sav]{Sav}
A. ~Savage,
\emph{On two geometric constructions of $U({\mathfrak{sl}}_n)$ and its representations},
J. Algebra 305 (2006), no. 2, 664--686.



\bibitem[Sp]{Sp}
N.~ Spaltenstein, \emph{Classes unipotentes et sous-groupes de Borel},
Lecture Notes in Mathematics, 946. Springer-Verlag, Berlin-New York, 1982.



\bibitem[TZ]{TZ}
James Tao, Yifei Zhao,\
{\em Extensions by $K_2$ and factorization line bundles},
matharXiv:1901.08760,


\bibitem[Wa1]{Wa1}
W. ~ Wang, \emph{Lectures at Yale University}, 1998/1999.

\bibitem[Wa2]{Wa2}
W. ~ Wang, \emph{Lagrangian Construction of the $(gl_n, gl_m)$-Duality},
Commun. Contemp. Math. 3, (2001), 201--214.

\bibitem[W]{W}
G. ~ Wilson, \emph{Collisions of Calogero-Moser particles and an
adelic Grassmannian. With an appendix by I.~G.~ Macdonald},
Invent. Math. 133 (1998), no. 1, 1--41.


\end{thebibliography}
\end{document}